\documentclass{iopart}
\usepackage{diagbox}%
\usepackage{iopams}
\eqnobysec
\usepackage{graphicx}
\usepackage{epstopdf}
\usepackage{color}
\usepackage{amsfonts}
\usepackage{amssymb}
 \expandafter\let\csname equation*\endcsname\relax
 \expandafter\let\csname endequation*\endcsname\relax
\usepackage{amsmath}
\usepackage{amsxtra}
\usepackage{amsthm}
\usepackage{latexsym}
\usepackage{multirow}

\newtheorem{lemma}{Lemma}[section]
\newtheorem{theorem}[lemma]{Theorem}

\newtheorem{definition}[lemma]{Definition}

\begin{document}
\jl{5}
\title[]
{Elastic-net regularization for nonlinear electrical impedance tomography with a splitting approach}

\author{Jing Wang$^{1}$, Bo Han$^{2}$ and Wei Wang$^3$}

\medskip

\address{$^1$School of Mathematical Sciences, Heilongjiang University, Harbin, Heilongjiang 150080, China (jingwangmath@hlju.edu.cn)}
\address{$^2$Department of Mathematics, Harbin Institute of Technology, Harbin, Heilongjiang 150001, China (bohan@hit.edu.cn)}
\address{$^3$College of Mathematics, Physics and Information Engineering, Jiaxing University, Zhejiang 314001, China (weiwang\_math@126.com) }





\begin{abstract}
Image reconstruction of EIT mathematically is a typical nonlinear and severely ill-posed inverse problem. Appropriate priors or penalties are required to enable the reconstruction.
The commonly used $l_{2}$-norm can enforce the stability to preserve local smoothness, and the current $l_{1}$-norm can enforce the sparsity to preserve sharp edges.
Considering the fact that $l_{2}$-norm penalty always makes the solution overly smooth and $l_{1}$-norm penalty always makes the solution too sparse, elastic-net regularization approach with a convex combination term of $l_{1}$-norm and $l_{2}$-norm emerges for fully nonlinear EIT inverse problems.
Our aim is to combine the strength of both terms: sparsity in the transform domain and smoothness in the physical domain, in an attempt to improve the reconstruction resolution and robustness to noise.
Nonlinearity and non-smoothness of the generated composite minimization problem make it challenging to find an efficient numerical solution. Then we develop one simple and fast numerical optimization scheme based on the split Bregman technique for the joint penalties regularization. The method is validated using simulated data for some typical conductivity distributions. Results indicate that the proposed inversion model with an appropriate parameter choice achieves an efficient regularization solution and enhances the quality of the reconstructed image.
\end{abstract}




%
\oddsidemargin 15mm
\evensidemargin 15mm
\topmargin  0 pt


%
\section{Introduction}
Electrical Impedance Tomography (EIT) is a noninvasive imaging technique aiming at visualizing the electrical conductivity distribution inside the interested object from the injected excitation currents and measured boundary voltages at surface electrodes.
With the advantages of low cost, high speed and lack of radiation, the main potential applications for EIT are, e.g., medical
imaging, geophysical exploration, and non-destructive testing of materials, see the references \cite{Cheney:1999,Borcea:2002,Lionheart:2005}.

The process of generating an image from the known currents and measured voltages data is called the image reconstruction, which is a highly nonlinear ill-posed problem.
So EIT imaging problem suffers from inherent low spatial resolution and instability.
In order to achieve desired images, EIT algorithms need to incorporate regularization to overcome the extreme sensitivity to modeling errors and measurement noise.
The penalty regularization technique is one most commonly used, whose basic idea is that some prior constraints are added up to the least-squares fitting term to preserve the stability of the inverse solution.
Different type of penalties have been employed for nonlinear EIT imaging.
$l_{2}$-norm regularization imposing smoothness priors is one of the most successful methods \cite{Tikhonov:2002,Tikhonov:2003}.
This regularization improves the stability at the cost of imaging resolution, so that the reconstructed images exhibit overly smooth features with fuzzy edges.
This implies that jumps and edges cannot be nicely reconstructed by $l_{2}$-norm regularization.
Total variation (TV) regularization has obvious edge-preserving and scale-dependent properties \cite{TV:2003}.
Applications of TV regularization for EIT imaging \cite{Levelset_TV:2005,TV:2007,TV:2010} is verified to work well for reconstructing large piecewise constant discontinuous conductivity distributions with sharp interfaces between the background and inclusions.
Furthermore, one hybrid regularization combining $l_{2}$-norm and TV prior is proposed to improve the reconstruction resolution and robustness to noise \cite{Liu:2013}.
In recent years, an alternative for penalty term is $l_{1}$-norm on the expansion coefficients of the solution with respect to a certain basis/frame.
It is a convex relaxation of the NP-hard $l_{0}$-norm sparsity metric, and is popularly used to enforce the sparsity.
It is shown that $l_{1}$-norm regularization has been receiving considerable interest for solving inverse problems \cite{IST:2004}, and can efficiently capture the small-scale details of the solution \cite{Loris:2007,CS:2012}.
In light of this, $l_{1}$-norm regularization has gained much popularity for nonlinear EIT imaging to improve the reconstruction quality \cite{BorsicL1_PDIPM:2012,Jin_L12:2012,Maass:2014,Jin_L1:2012}.

Note that $l_{1}$-norm regularized minimum often requires inverting potentially ill-conditioned operator and thus leads to numerical problems, then one remedy method named by 'elastic-net regularization' is proposed in an ill-posed problem framework \cite{Jin_Elastic:2009}.
Inspired by this, in order to alleviate the numerical problem and obtain a meaningful stable solution, we further regularize the sparsity-promoting $l_{1}$-norm regularization by adding the stable $l_{2}$-norm constraint term.
Moreover, considering the fact that $l_{2}$-norm always makes the solution excessively smooth and $l_{1}$-norm penalty always makes the solution too sparse, therefore this paper investigates an elastic-net regularization approach with stronger convex combination of $l_{1}$-norm and $l_{2}$-norm for fully nonlinear EIT inverse problem.
Our aim is to combine the strength of both terms: sparsity in the transformed domain and smoothness in the physical domain, in an attempt to improve the reconstruction resolution and robustness to noise.
This approach is more suited to reconstruct the smooth and edge regions for small spatial profiles.
Due to the non-differentiable $l_{1}$-norm prior, one simple and fast numerical scheme based on splitting Bregman iteration technique is developed efficiently, which allows separating the Gauss-Newton minimization from the proposed cost function.
Numerical simulations with synthetic data are implemented to validate the feasibility and effectiveness of the proposed inversion model.

The contents of this paper are organized as follows. Section 2 briefly
describes the basic mathematical model for the forward problem and the
inverse problem in EIT. In section 3, we design a compound regularized
inversion model and discuss its well-posedness. In addition, an efficient numerical implementation based on split Bregman technique is provided for handling the joint $l_{1}$-norm and $l_{2}$-norm penalties.
In section 4, numerical simulations with synthetic data are performed to
demonstrate that our proposed method has good performance for nonlinear EIT imaging problem.
Finally, we draw conclusions in section 5.

\section{Mathematical formulation for EIT}
We shall recapitulate the mathematical formulation for EIT problem in this section.
Since inverse imaging problems are characterized by the property of being sensitive to measurement noise and modeling errors, it is very important to model the measurement procedure accurately.
The complete electrode model (CEM) \cite{CEM,Somersalo:1992} is currently the most accurate mathematical model for reproducing EIT experimental data, as it takes into account the presence of the discrete electrodes and effects of the contact impedance, and can match the measurement precision of the experiment well.
In this study, we use the CEM for the forward simulation, which can be stated as follows:
\begin{equation}
\label{forward}
\left\{
\begin{array}{rl}
\nabla\cdot(\sigma\nabla u)&=0,~\textup{in}~\Omega,\cr
\sigma\frac{\partial u}{\partial n}&=0,~\textup{on}~\partial\Omega\setminus\Gamma_{e},\cr
\int_{e_{l}}\sigma\frac{\partial u}{\partial n}dS&=I_{l},~\textup{for}~l=1,2,\ldots,L,\cr
u+z_{l}\sigma\frac{\partial u}{\partial n}&=U_{l},~\textup{on}~e_{l},~l=1,2,\ldots,L,
\end{array}
\right.
\end{equation}
where $\sigma$ is the electrical conductivity, $u$ is the electrical potential inside the object domain $\Omega$, also referred to as the inner potential, $z_{l}$ is the effective contact impedance between the $l$\textup{th} electrode and tissue, $L$ is the number of electrodes in the system, $e_{l}$ denotes the area of the $l$\textup{th} electrode, $\Gamma_{e}=\cup_{l=1}^{L}e_{l}\subset\partial\Omega$, and $n$ is the unit outward normal to $\partial\Omega$ and $\frac{\partial u}{\partial n}=\nabla u\cdot n$.
$U_l$ and $I_l$ are the electrical potential and electrical current on the $l$\textup{th} electrode, respectively.
Additionally, in order to ensure the solvability of the elliptic boundary value problem (\ref{forward}) and the uniqueness of the solution, the following two other conditions should be satisfied:
\begin{equation}
\sum_{l=1}^{L}I_{l}=0,~~\sum_{l=1}^{L}U_{l}=0.
\end{equation}
We denote vectors consisting of the electrical potentials
and currents on the electrodes by $U$ and $I$, i.e. $U=(U_{1},U_{2},\ldots,U_{L})\in \mathbb{R}^{L}$ and $I=(I_{1},I_{2},\ldots,I_{L})\in \mathbb{R}^{L}$, respectively.

EIT is composed of forward problem and inverse problem.
The forward problem is to find potentials $u$ in $\Omega$ and $U$ on the electrodes from known conductivity $\sigma$, injected currents $I$, and contact impedances $\{z_{l}\}$.
Hence, the solution of the forward problem is equivalent to solving the boundary value problem \eqref{forward}.
Mathematically the existence and uniqueness of the forward solution for problem \eqref{forward} has been proven in \cite{Somersalo:1992}.
Since the forward solution cannot be solved analytically for complex geometries, we have to resort to the numerical technique to find an approximate solution.
In this paper, the finite element method (FEM) is employed for the discretization of the forward problem, and we refer to \cite{Woo:1994,Vauhkonen:1997} for more details on the forward solution and FEM approximation of the CEM.

In practice, the conductivity $\sigma$ is unknown. What we know is merely all pairs of injected current data $I$ and resulting voltage data $U$ on the electrodes.
The inverse problem seeks to reconstruct the unknown conductivity $\sigma$ from these measured data on the electrodes.
With the aid of FEM discretization, the relation between the applied
electrode current data $I$ and the computed electrode voltage data $U(\sigma)$ can be defined by
\begin{equation}
U(\sigma)=F(\sigma;I),
\end{equation}
where $F$ is the solution operator (also called the forward operator), which depends linearly on $I$ but nonlinearly on $\sigma$.
For a fixed current vector $I$, we can view $F(\sigma;I)$ as
a function of $\sigma$ only.
The measured data is unavoidably noisy, thus using an additive Gaussian model for the measurement error, the observation model response can be formulated as the nonlinear operator equation:
\begin{equation}
U^{\delta}=F(\sigma)+e,
\end{equation}
where $U^{\delta}$ denotes the noisy potential measurements, $e$ denotes the Gaussian distributed noise error with zero mean in the measurements, and the admissible set $\mathcal{A}$ is given by
$$\mathcal{A}=\{\sigma\in L^{\infty}(\Omega):\lambda\leq\sigma\leq\lambda^{-1}~a.e.~in~ \Omega,~\sigma|_\Gamma=\sigma_0|_\Gamma\},$$
for some fixed constant $\lambda\in(0,1)$.
The condition $\sigma|_\Gamma=\sigma_0|_\Gamma$ reflects the fact that inclusions are located inside the object.

Therefore, the EIT inverse problem mathematically is to find an approximation to the physical conductivity distribution $\sigma^{*}$ from the finite set of noisy potential data $U^{\delta}$ of voltages $U(\sigma^{*})$ (or denoted by $U$), also called as image reconstruction.
This constitutes a nonlinear and severely ill-posed problem in the sense that small errors in the data can lead to very large deviations in the solutions.

\section{Regularization and iterative realization}
In order to reconstruct the conductivity stably, some sort of regularization is required.
In this section we devise the elastic-net regularization, and provide one effective algorithm for minimizing the elastic-net functional.

\subsection{Elastic-net regularization}
In this study, image reconstruction for nonlinear EIT inverse problem is based on a traditional nonlinear least-squares approach $\min\limits_{\sigma}\big{\|}F(\sigma)-U^{\delta}\big{\|}^{2}$.
As this forms an ill-posed problem with a high condition number, the regularized least-squares need to be applied.
As mentioned before, the regularizer we seek to in this paper contains two parts: a sparsity-promoting $l_{1}$-norm penalty term denoted as $S_1(\sigma)$ and a smoothing penalty term denoted as $S_2(\sigma)$.

It is noteworthy that the sparsity in $l_{1}$-norm framework here mainly refers to two types: the space-domain sparseness and the transform-domain sparseness.
Let $\sigma_{0}$ be the background conductivity.
When $\sigma_{0}$ is known and the inhomogeneity $\sigma-\sigma_{0}$ has a certain sparsity degree in the space domain (i.e., most elements of $\sigma-\sigma_{0}$ are zero), we say that the conductivity distributions have the spatial sparsity, and $S_1(\sigma)$ can be defined by $l_{1}$-norm of the elements of the inhomogeneity as follows
$$S_1(\sigma)=\big{\|}\sigma-\sigma_{0}\big{\|}_{1}
=\sum_{k}\big{|}[\sigma-\sigma_{0}]_k\big{|}.$$
When $\sigma_0$ is unknown, with the aid of the sparsity with respect to a preassigned basis/frame $\Phi:=\{\phi_{k}\}_{k\in \mathbb{N}}$, i.e., there are only finitely many (or a small number of) nonzero expansion coefficients $[\Phi\sigma]_k:=\langle\sigma,\phi_k\rangle$, $S_1(\sigma)$ can be defined by
$l_{1}$-norm of the transform coefficients as follows $$S_1(\sigma)=\big{\|}\Phi\sigma\big{\|}_{1}
=\sum_k\big{|}[\Phi\sigma]_k\big{|}.$$
For the smoothing penalty term, given an a priori guess $\sigma_{ref}$, we utilize the standard $l_{2}$-norm $$S_2(\sigma)=\big{\|}R(\sigma-\sigma_{ref})\big{\|}_{2}^{2},$$
where $R$ is the regularization matrix.

Moreover, we inspired by the 'elastic-net' \cite{Jin_Elastic:2009} introduce the hybrid regularizer
\begin{equation}\label{regularizer1}
\Theta_\beta(\sigma)=(1-\beta)S_1(\sigma)+\beta S_2(\sigma),
\end{equation}
which is a convex combination of $l_{1}$-norm and $l_{2}$-norm, so called 'elastic-net regularizer'.
Consequently, we in this paper consider such variational regularization model
\begin{eqnarray}\label{L1L2}
\mathcal{J}_{\alpha,\beta}(\sigma)
=\big{\|}F(\sigma)-U^{\delta}\big{\|}^{2}
+\alpha\Theta_{\beta}(\sigma),
\end{eqnarray}
where $\alpha>0$ is a regularization parameter.
For $\beta=1$ model \eqref{L1L2} reduces to $l_{2}$-norm regularization; while for $\beta=0$ model \eqref{L1L2} reduces to $l_{1}$-norm regularization.
We here assign a non-zero value for $\beta$ (i.e., $0<\beta<1$), since we want to enforce the structural properties of sparsity for all iterations.
This will typically relieves overly smooth behavior while well preserving sharp edges.

\subsection{Well-posedness for elastic-net regularization}
We denote by $\sigma_{\alpha,\beta}^{\delta}$ the minimizer of the functional $\mathcal{J}_{\alpha,\beta}(\sigma)$ with noisy data $U^\delta$ and parameters $(\alpha,\beta)$.
Throughout this paper, we will need the following definition of $\Theta_{\beta}$-minimizing solutions, as well as some properties of the forward operator $F(\sigma)$ for details and proofs we refer to \cite{Jin_Analysis:2012}.
\begin{definition}
An element $\sigma^{\dag}$ is said to be a $\Theta_{\beta}$-minimizing solution to $F(\sigma)=U$ if $F(\sigma^{\dag})=U$ and
$$\Theta_{\beta}(\sigma^{\dag})=\min_{\sigma\in\mathcal{A}}
\{\Theta_{\beta}(\sigma):F(\sigma)=U\}.$$
\end{definition}
\noindent Such an exact solution $\sigma^{\dag}$ may take a different form if the value of $\beta$ varies.
\begin{lemma}
Let $I\in \mathbb{R}_{\diamond}^{L}$ be a fixed current vector. For any $\sigma, \sigma+\delta\sigma\in \mathcal{A}$, the forward operator $F(\sigma)$ is uniformly bounded and continuous, and is also Fr$\acute{e}$chlet differentiable with respect to $\sigma$.
Moreover, the Fr$\acute{e}$chlet derivative $F'(\sigma)$ is bounded and Lipschitz continuous, and the following estimations are satisfied:
\begin{displaymath}
\begin{aligned}
&\|F(\sigma+\delta\sigma)-F(\sigma)\|_{H}\leq C\|\delta\sigma\|_{L^{\infty}(\Omega)},\\
&\big{\|}F'(\sigma+\delta\sigma)-F'(\sigma)\big{\|}_{\mathcal{L}(L^{\infty}(\Omega),~H)}\leq L\big{\|}\delta\sigma\big{\|}_{L^{\infty}(\Omega)},\\
&\|F(\sigma+\delta\sigma)-F(\sigma)-F'(\sigma)\delta\sigma\|_{L^{2}(\partial\Omega)}\leq \frac{L}{2}\|\delta\sigma\|^{2}_{L^{\infty}(\Omega)}.
\end{aligned}
\end{displaymath}
\end{lemma}
We next address the well-posedness of the regularization functional \eqref{L1L2}.
More precisely, by the results in Lemma 3.2, we have the following existence and stability results.
The proofs of the next results are standard and analogous to those in \cite{Jin_Elastic:2009}, and are thus omitted.

\begin{theorem}
Assuming that both parameters $\alpha>0$ and $0<\beta<1$ for any given measurement data $U^{\delta}$.
The problem (\ref{L1L2}) is well-posed and consistent, i.e.,

(i) There exists at least one minimizer $\sigma_{\alpha,\beta}^{\delta}$ to the functional $\mathcal{J}_{\alpha,\beta}(\sigma)$ on the admissible set $\mathcal{A}$.

(ii) Let $\{U^{n}\}\subset L^{2}(\partial\Omega)$ be a sequence of noisy data converging to $U^{\delta}$, and let $\sigma_{\alpha,\beta}^{n}$ be the minimizer of the problems:
$$\mathcal{J}_{\alpha,\beta}(\sigma;U^{n})=\big{\|}F(\sigma)-U^{n}\big{\|}_{L^{2}(\partial\Omega)}^{2}
+\alpha \Theta_{\beta}(\sigma).$$
Here, to show explicitly the dependence of $\mathcal{J}_{\alpha,\beta}(\sigma)$ on the measurement data, we use $\mathcal{J}_{\alpha,\beta}(\sigma;U^{n})$ instead of $\mathcal{J}_{\alpha,\beta}(\sigma)$.
Then a subsequence $\{\sigma_{\alpha,\beta}^{n_1}\}$ of $\{\sigma_{\alpha,\beta}^{n}\}$ and $\tilde{\sigma}\in\mathcal{A}$ exist such that $\{\sigma_{\alpha,\beta}^{n_1}\}$ converges to $\tilde{\sigma}$, and
\begin{eqnarray}\label{RM1}
\lim_{n_{1}\rightarrow\infty}\Theta_{\beta}
(\sigma_{\alpha,\beta}^{n_{1}})=\Theta_{\beta}(\tilde{\sigma}).
\end{eqnarray}
Furthermore, $\tilde{\sigma}$ is a solution to $\mathcal{J}_{\alpha,\beta}(\sigma)$.

(iii) We take a sequence of noise level $\{\delta_{n}\}$.
Suppose that $\|U-U^{\delta_{n}}\|\leq\delta_{n}$ and the regularization parameter $\alpha_{n}:=\alpha(\delta_{n})$ satisfy
$$\lim_{\delta_{n}\rightarrow0}\alpha_{n}=0~and~
\lim_{\delta_{n}\rightarrow0}\frac{\delta_{n}^{2}}{\alpha_{n}}=0.$$
Then the sequence of corresponding minimizers $\{\sigma_{\alpha_{n},\beta}^{n}\}$ has a subsequence converging to an $\Theta_{\beta}$-minimizing solution $\sigma^{\dag}$ as $\delta_{n}\rightarrow0$.
\end{theorem}

\subsection{Iterative algorithm realization}
Having described the elastic-net functional, we now proceed to the algorithmic part of minimizing the functional.
As the forward operator $F$ depends nonlinearly on the parameter $\sigma$, it is almost impossible to find an accurate solution by any non-iterative algorithm with a simplified linear model.
In order to obtain reasonable reconstructions, we have to resort to the iterative method by updating the solution multiple times.
For the nonlinear inversion model, we usually follow the standard local linearized iterative approach \cite{Engl_Book:1996}.
Using an iterative framework, let $F(\sigma)$ be the first order Taylor expansion at $\sigma^{k}$, i.e.,
$$F(\sigma)\approx F(\sigma^{k})+F'(\sigma^{k})(\sigma-\sigma^{k}),$$
our problem becomes that for given $\sigma^{k}$, to
solve
\begin{equation}
\min_{\delta\sigma}\bigg{\{}\big{\|}F'(\sigma^{k})\delta\sigma+F(\sigma^{k})-U^{\delta}\big{\|}^{2}+\alpha_{k}
\Theta_\beta(\sigma^{k}+\delta\sigma)\bigg{\}},
\end{equation}
set $\sigma^{k+1}=\sigma^{k}+\delta\sigma$, $\alpha_{k+1}<\alpha_{k}$ and repeat the iteration
until an expected number of iterations or a satisfactorily low discrepancy is achieved.

However, the presence of $l_{1}$-norm results in non-smooth objective functional, and so the minimizer can no longer be explicitly computed.
A review of optimization algorithms for $l_{1}$-norm minimization problem is offered by Loris \cite{OptimizaitonReview:2009}.
Thereinto, the split Bregman method (SBM) \cite{SplitBregman:2009} on account of the fast computational efficiency and the satisfactory
numerical performance has received wide applications.

\subsubsection{SBM for the transform domain}
In the transform domain framework, assuming that conductivity distributions have the sparsity representation with respect to the chosen basis $\Phi$, the minimizing inversion model becomes
\begin{eqnarray}\label{L1L2S}
\min\limits_{\delta\sigma}\bigg{\{}\big{\|}F'(\sigma^{k})\delta\sigma+F(\sigma^{k})-U^{\delta}\big{\|}^{2}
&+\alpha_{k}\bigg{[}(1-\beta)\big{\|}\Phi(\sigma^{k}+\delta\sigma)\big{\|}_{1}\cr
&+\beta\big{\|}R(\sigma^{k}+\delta\sigma-\sigma_{ref})\big{\|}_{2}^{2}\bigg{]}\bigg{\}}.
\end{eqnarray}
According to the basic idea of Bregman iteration, we introduce a new auxiliary variable $d$ and enforce an equality constraint $d=\Phi(\sigma^{k}+\delta\sigma)$, which leads to the following optimization subproblems:
\begin{equation}\label{sub-problems}
\begin{array}{ll}
&(\delta\sigma,d^{k+1})=\arg\min\limits_{\delta\sigma,d}\bigg{\{}\big{\|}F'(\sigma^{k})\delta\sigma+F(\sigma^{k})-
U^{\delta}\big{\|}^{2}\cr
&\qquad\qquad\qquad\qquad\qquad
+\alpha_{k}\bigg{[}(1-\beta)\big{\|}d\big{\|}_{1}+\beta\big{\|}R(\sigma^{k}+\delta\sigma-\sigma_{ref})
\big{\|}^{2}\bigg{]}\cr
&\qquad\qquad\qquad\qquad\qquad
+\mu\big{\|}d-\Phi(\sigma^{k}+\delta\sigma)-b_{d}^{k}\big{\|}^{2}\bigg{\}},\cr
&b_{d}^{k+1}=b_{d}^{k}+\Phi(\sigma^{k}+\delta\sigma)-d^{k+1},
\end{array}
\end{equation}
where $\mu>0$ is the relaxation factor associated with the equality
constraint.

By virtue of the alternating direction iteration optimization techniques, we next address the solution of subproblems (\ref{sub-problems}).
Minimizing (\ref{sub-problems}) with respect to $\delta\sigma$ for given $d$ and $b_{d}^{k}$ leads to the following standard least square problem
\begin{eqnarray}\label{sub1}
\min\limits_{\delta\sigma}\bigg{\{}\big{\|}F'(\sigma^{k})\delta\sigma+F(\sigma^{k})-
U^{\delta}\big{\|}^{2}&+\alpha_{k}\beta\big{\|}R(\sigma^{k}+\delta\sigma-\sigma_{ref})
\big{\|}^{2}\cr
&+\mu\big{\|}d-\Phi(\sigma^{k}+\delta\sigma)-b_{d}^{k}\big{\|}^{2}\bigg{\}},
\end{eqnarray}
where we ignore the constant terms independent of $\delta\sigma$.
To obtain the minimum of equation (\ref{sub1}), solving the first-order optimality condition with respect to $\delta\sigma$ yields the Gauss-Newton update
\begin{eqnarray}\label{GN-update}
\big{[}F'(\sigma^{k})^{T}F'(\sigma^{k})&+\alpha_{k}\beta R^{T}R+\mu
I\big{]}\delta\sigma=F'(\sigma^{k})^{T}(U^{\delta}-F(\sigma^{k}))\cr
&-\alpha_{k}\beta R^{T}R(\sigma^{k}-\sigma_{ref})-\mu(\sigma^{k}+\Phi^{T}(b_{d}^{k}-d)),
\end{eqnarray}
which can be efficiently computed using the conjugate gradient (CG) method.
The minimization with respect to $d$ for given
$\sigma^{k+1}=\sigma^{k}+\delta\sigma$ and $b_{d}^{k}$ is as follows:
\begin{equation}\label{sub2}
\min\limits_{d}\bigg{\{}\alpha_{k}(1-\beta)\big{\|}d\big{\|}_{1}+
\mu\big{\|}d-\Phi\sigma^{k+1}-b_{d}^{k}\big{\|}^{2}\bigg{\}}.
\end{equation}
This sub-problem is the $l_1$-norm denoising problem, and the minimizer can be solved quickly with the soft shrinkage/thresholding operator, which can be formulated by
\begin{equation}\label{sub2_solution}
d^{k+1}=\mathrm{shrinkage}(\Phi\sigma^{k+1}+b_{d}^{k},\eta),
\end{equation}
where $\eta=\frac{\alpha_{k}(1-\beta)}{2\mu}$, $\mathrm{shrinkage}(x,t)=\mathrm{sign}(x)\mathrm{max}(|x|-t,0)$ for $x$ and $t$.

Finally, we summarize the complete procedure for solving the optimization problem (\ref{L1L2S}) in Algorithm 1.
$\alpha_{0}$ is the chosen initial regularization parameter and $q_{\alpha}$ is the step factor used to adjust the regularization parameter adaptively.
Another problem is the stopping criterion for the iterative reconstruction.
We here consider the Morozov's discrepancy principle \cite{Morozov} for the Gauss-Newton iteration.
It states that the residual shall not be smaller than the error in the measurement, i.e., when the data residual falls below the noise level.
Thus, $\alpha$ should be chosen such that
\begin{equation}
\big{\|}F(\sigma_{\alpha})-U^{\delta}\big{\|}\approx\big{\|}U-U^{\delta}\big{\|},
\end{equation}
where $\sigma_{\alpha}$ denotes the reconstructed conductivity distribution with regularization parameter $\alpha$.

\begin{flushleft}
\begin{tabular}{l}
\hline
\textbf{Algorithm~1.} SBM for the transform domain\\
\hline
Choose $\alpha_{0},~q_{\alpha},~\beta,~\mu$ and $\tau$\\
Set $k=0,~\sigma^{k}=\sigma_{ref},~d^{0}=0,~b_{d}^{0}=0$\\
\bf{repeat}\\
\qquad Set $j=0,~\eta=\alpha_{k}(1-\beta)/(2\mu$)\\
\qquad \bf{repeat}\\
\qquad \qquad Compute $\delta\sigma^{j+1}$ as minimizer of
(\ref{sub1}) by solving (\ref{GN-update})\\
\qquad \qquad Compute $d^{j+1}$ as minimizer of (\ref{sub2}) using shrinkage operator, i.e. \\
\qquad \qquad \qquad \qquad \quad $d^{j+1}=\mathrm{shrinkage}(\Phi(\sigma^{k}+\delta\sigma^{j+1})+b_{d}^{j},\eta)$\\
\qquad \qquad Set $b_{d}^{j+1}=b_{d}^{j}+\Phi(\sigma^{k}+\delta\sigma^{j+1})-d^{j+1}$\\
\qquad \qquad update $j\leftarrow j+1$\\
\qquad \bf{until} $\big{\|}\delta\sigma^{j}\big{\|}<\tau$\\
\qquad Set $\sigma^{k+1}=\sigma^{k}+\delta\sigma^{j}$\\
\qquad Update $\alpha_{k+1}=q_{\alpha}\alpha_{k},~k\leftarrow k+1$\\
\bf{until} $\big{\|}F(\sigma^{k})-U^{\delta}\big{\|}<\big{\|}U-U^{\delta}\big{\|}$\\
\hline
\end{tabular}
\end{flushleft}

\subsubsection{SBM for the space domain}
In the space domain framework, assuming that conductivity distributions have the spatial sparsity, that is, most elements for the inhomogeneity $\sigma-\sigma_{0}$ are zero, the corresponding inversion model becomes
\begin{eqnarray}\label{L1L2SS}
\min\limits_{\delta\sigma}\bigg{\{}\big{\|}F'(\sigma^{k})\delta\sigma+F(\sigma^{k})-U^{\delta}\big{\|}^{2}
&+\alpha_{k}\bigg{[}(1-\beta)\big{\|}\sigma^{k}+\delta\sigma-\sigma_{0}\big{\|}_{1}\cr
&+\beta\big{\|}R(\sigma^{k}+\delta\sigma-\sigma_{ref})\big{\|}_{2}^{2}\bigg{]}\bigg{\}}.
\end{eqnarray}

Analogous to Algorithm 1, we by introducing the auxiliary variable $d$ and enforcing the equality constraint $d=\sigma^{k}+\delta\sigma-\sigma_{0}$ obtain three-stage iterations:
\begin{displaymath}
\begin{array}{lll}
\textrm{Step1}:&\delta\sigma=\arg\min\limits_{\delta\sigma}\bigg{\{}\big{\|}
F'(\sigma^{k})\delta\sigma+F(\sigma^{k})-
U^{\delta}\big{\|}^{2}+\alpha_{k}\beta\big{\|}R(\sigma^{k}+\delta\sigma
-\sigma_{ref})\big{\|}^{2}\cr
&\quad\qquad\qquad\qquad\qquad\qquad\qquad +\mu\big{\|}d-(\sigma^{k}+\delta\sigma-\sigma_{0})-b_{d}^{k}\big{\|}^{2}\bigg{\}},\cr
\textrm{Step2}:&d^{k+1}=\arg\min\limits_{d}\bigg{\{}\alpha_{k}(1-\beta)\big{\|}d\big{\|}_{1}+
\mu\big{\|}d-(\sigma^{k}+\delta\sigma-\sigma_{0})-b_{d}^{k}\big{\|}^{2}\bigg{\}},\cr
\textrm{Step3}:&b_{d}^{k+1}=b_{d}^{k}+(\sigma^{k}+\delta\sigma-\sigma_{0})-d^{k+1},
\end{array}
\end{displaymath}
The variational equation for Step1 is given by
\begin{eqnarray}\label{GN-update2}
\big{[}F'(\sigma^{k})^{T}F'(\sigma^{k})&+\alpha_{k}\beta R^{T}R+\mu I\big{]}\delta\sigma=F'(\sigma^{k})^{T}(U^{\delta}-F(\sigma^{k}))\cr
&-\alpha_{k}\beta
R^{T}R(\sigma^{k}-\sigma_{ref})-\mu(\sigma^{k}-\sigma_{0}+b_{d}^{k}-d).
\end{eqnarray}
Step2 can be solved efficiently through the shrinkage operator as follows
\begin{equation}\label{sub2_solution2}
d^{k+1}=\mathrm{shrinkage}(\sigma^{k}+\delta\sigma-\sigma_{0}+b_{d}^{k},\eta).
\end{equation}
where $\eta=\frac{\alpha_{k}(1-\beta)}{\tilde{2\mu}}$.
The full procedure for solving the optimization problem (\ref{L1L2SS}) is summarized in Algorithm 2.

\begin{flushleft}
\begin{tabular}{l}
\hline
\textbf{Algorithm~2.} SBM for the space domain\\
\hline
Choose $\alpha_{0},~q_{\alpha},~\beta,~\mu$ and $\tau$\\
Set $k=0,~\sigma^{k}=\sigma_{ref},~d^{0}=0,~b_{d}^{0}=0$\\
\bf{repeat}\\
\qquad Set $j=0,~\eta=\alpha_{k}(1-\beta)/(2\mu)$\\
\qquad \bf{repeat}\\
\qquad \qquad Compute $\delta\sigma^{j+1}$ by solving the variational equation (\ref{GN-update2})\\
\qquad \qquad Compute $d^{j+1}$ using the shrinkage operator (\ref{sub2_solution2}), i.e.\\
\qquad \qquad \qquad \qquad
$d^{j+1}=\mathrm{shrinkage}(\sigma^{k}+\delta\sigma^{j+1}-\sigma_{0}+b_{d}^{k},\eta)$\\
\qquad \qquad Set $b_{d}^{j+1}=b_{d}^{j}+(\sigma^{k}+\delta\sigma^{j+1}-\sigma_{0})-d^{j+1}$\\
\qquad \qquad update $j\leftarrow j+1$\\
\qquad \bf{until} $\big{\|}\delta\sigma^{j}\big{\|}<\tau$\\
\qquad Set $\sigma^{k+1}=\sigma^{k}+\delta\sigma^{j}$\\
\qquad Update $\alpha_{k+1}=q_{\alpha}\alpha_{k},~k\leftarrow k+1$\\
\bf{until} $\big{\|}F(\sigma^{k})-U^{\delta}\big{\|}<\big{\|}U-U^{\delta}\big{\|}$\\
\hline
\end{tabular}
\end{flushleft}

\section{Results from numerical simulations}
In order to evaluate the feasibility and effectiveness of the proposed elastic-net regularization approach, we tested it and compared it with TV regularization, $l_{1}$-norm regularization and $l_{2}$-norm regularization using numerical simulations with synthetic data.
We selected a 16-electrode numbered clockwise system for circular domain and used constant injection current between adjacent electrodes and adjacent voltage measurement between all other electrodes.

All the simulations are performed on two dimensional numerical phantoms, as this allows an easier visualization of the results and faster simulation.
The simulated data were generated using the finite element forward model based on the two-dimensional EIDORS demo version developed by Vauhkonen \cite{Vauhkonen:2001}.
To avoid committing "the inverse crime" \cite{InverseCrime}, two different finite elements meshes were used for the forward and inverse solvers, shown in Figure 1.
For the forward simulations we used the FEM with fine meshes of 1049 nodes and 1968 elements.
For the inverse computations we used the FEM with coarse meshes of 279 nodes and 492 elements to reduce the computational burden.
The small green strips on the perimeter of the meshes show the position of the ideal electrodes.
In the reconstructions, contact impedances under the electrodes were assumed to be known, and are set to be 0.05 for all electrodes.

\begin{figure}[htbp]
\centering
\includegraphics[width=1.65in]{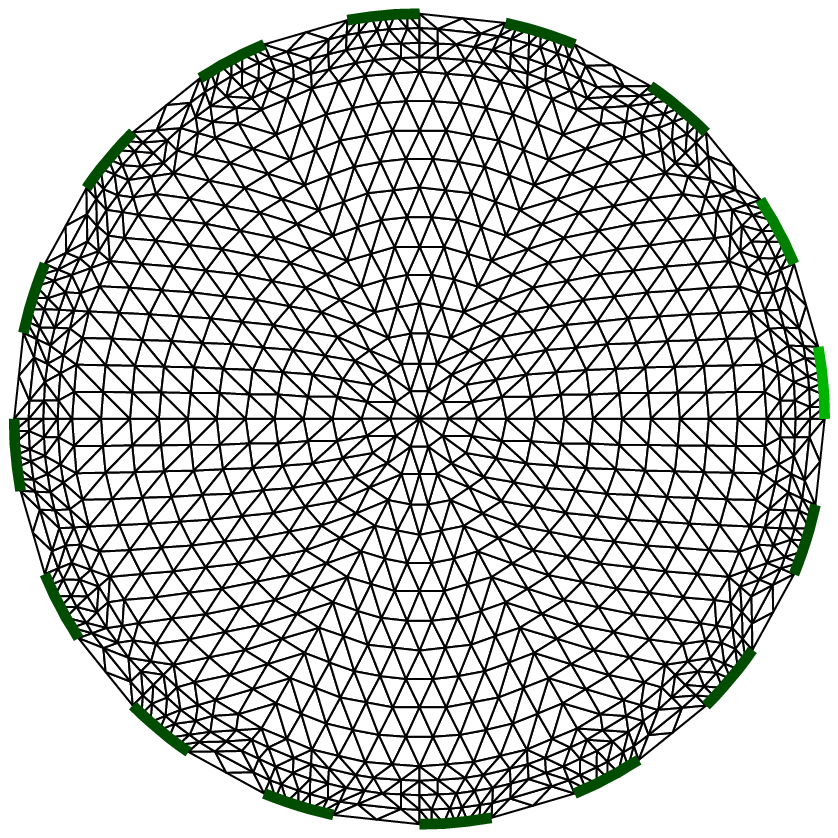}
\includegraphics[width=1.65in]{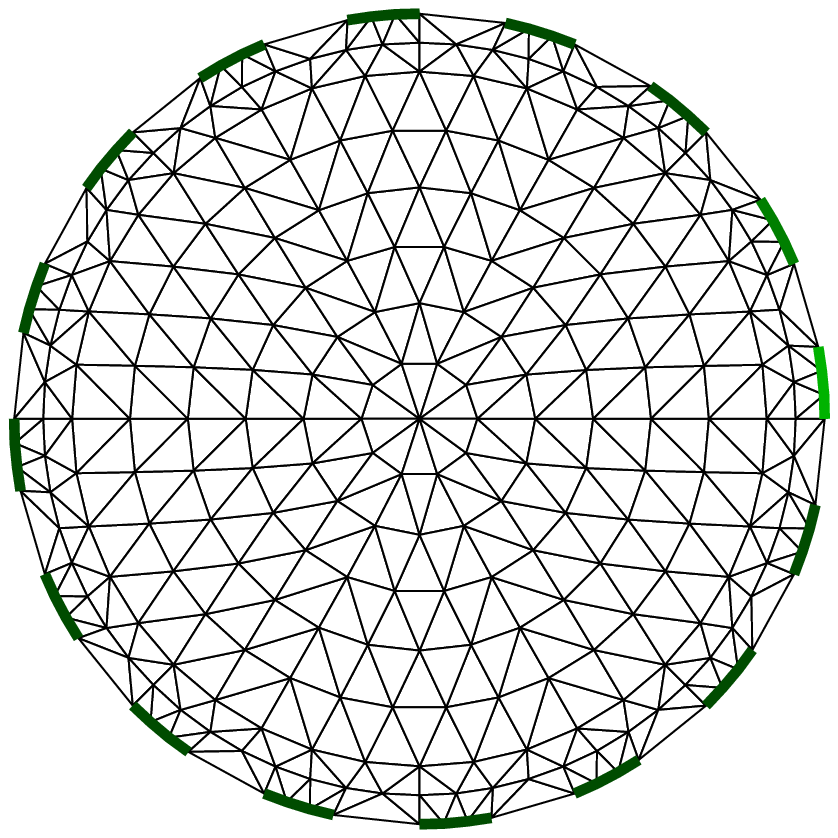}
\caption{The FEM meshes for numerical and experimental studies. Left: forward fine mesh with 1968 triangulation elements and 1049 nodes; Right: inverse coarse mesh with 492 triangulation elements and 279 nodes.}
\label{Figure 1}
\end{figure}

\subsection{Simulation phantoms}
The conductivity distribution $\sigma$ is assumed to be piecewise constant, homogeneous in $\Omega$ except for several inclusions $\{B_{j}\}_{j=1}^{p}$ ($p$ is a number), where $\sigma$ is assumed to be different from the background:
\begin{displaymath}
\sigma(r) = \left\{ \begin{array}{ll}
\sigma_{0}, & r\in\Omega\backslash\cup_{j}B_{j},\\
\sigma_{j}, & r\in B_{j}.
\end{array} \right.
\end{displaymath}
The boundary of the inclusion $B_{j}$ is denoted by $\partial B_{j}$.
The sharp edges refers to the boundaries of all the inclusions, i.e. $\cup_{j}\partial B_{j}$, while the smooth parts is $\Omega\backslash\cup_{j}\partial B_{j}$.
We here consider the following three simulation circular phantoms with small inclusions to test our proposed method, as depicted in Figure 2.
Phantom A consists of two inclusions ($p=2$), Phantom B consists of three inclusions ($p=3$), and Phantom C concerns five inclusions ($p=5$).
The resistivity values are 4 $\Omega m$ for the background, 1 $\Omega m$ for the blue color parts and 8 $\Omega m$ for the red color parts.

\begin{figure}[htbp]
\centering
\includegraphics[width=1.5in]{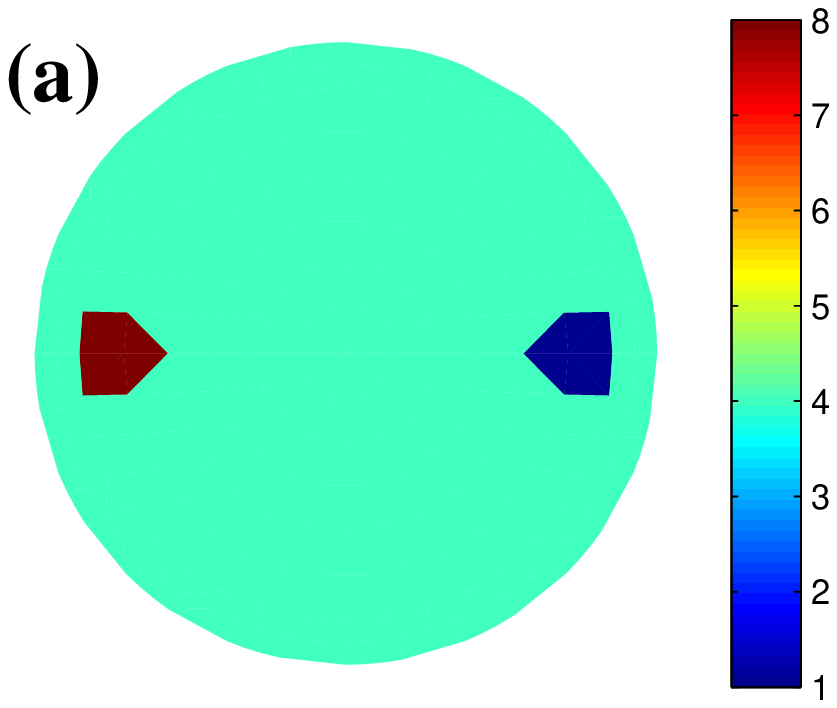}
\includegraphics[width=1.5in]{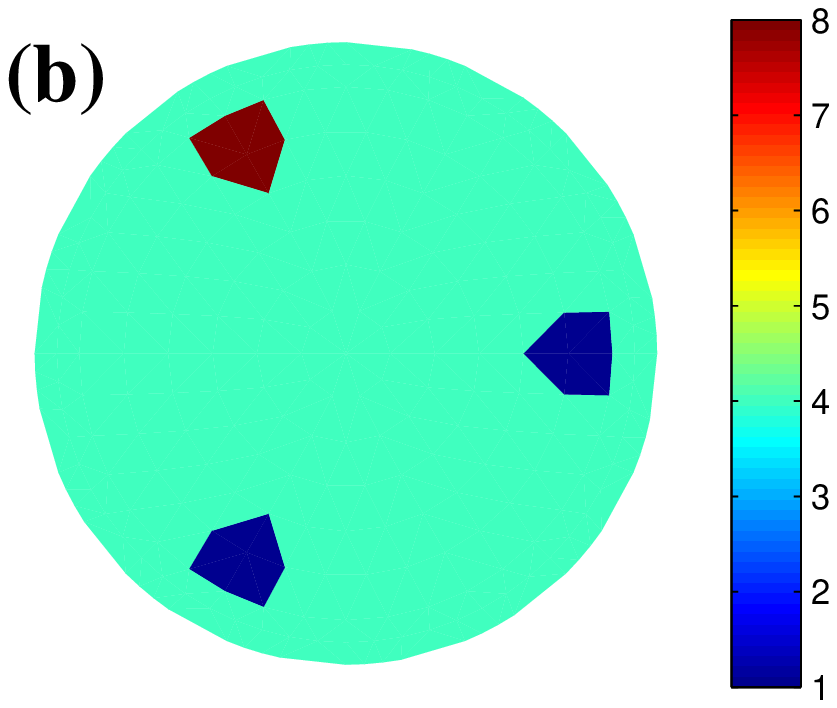}
\includegraphics[width=1.5in]{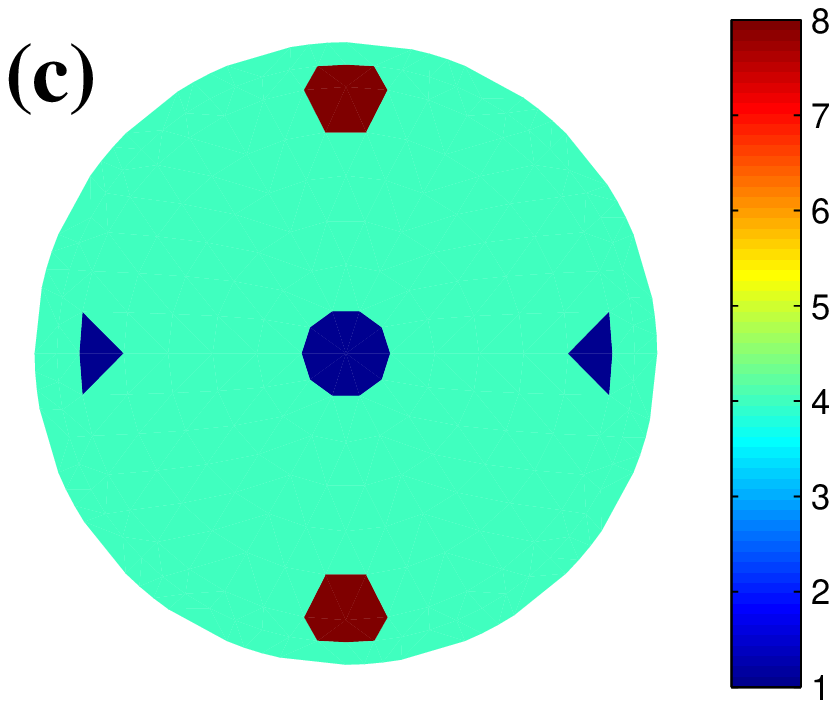}
\caption{Three simulation phantoms. (a) Phantom A.  (b) Phantom B.   (c) Phantom C.}
\label{Figure 2}
\end{figure}

The conductivity $\sigma$ for these three phantoms are apparently not sparse in the space domain.
We have to seek some sort of sparsity representation method.
We can resort to the sparsity representation in wavelet transform domain.
Here take $\Phi$ to be Haar wavelet.
Figure 3 displays the vertical bar graph of the wavelet expansion coefficients.
It can be seen that most coefficients are around zero, that is, there exist many small coefficients.
Beyond that, let $\sigma_0$ be known, and we can use the sparsity of the inhomogeneity ($\sigma-\sigma_0$) in the space domain.
Figure 4 shows the vertical bar graph of the elements for inhomogeneities.
It is obvious that most elements for $\sigma-\sigma_0$ are zero.
Hence both Haar wavelet transformation and inhomogeneity can effectively provide one sparsity representation.
But Figure 4 shows better sparsity for the phantoms.

\begin{figure}[htbp]
\centering
\includegraphics[width=1.65in]{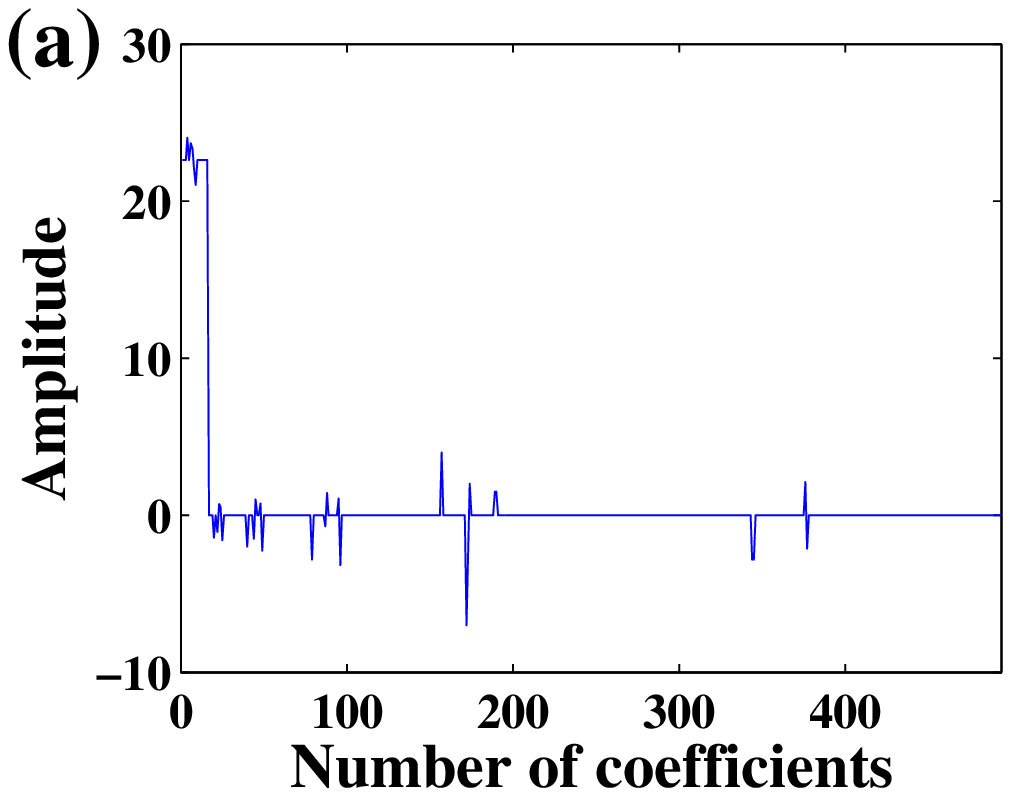}
\includegraphics[width=1.65in]{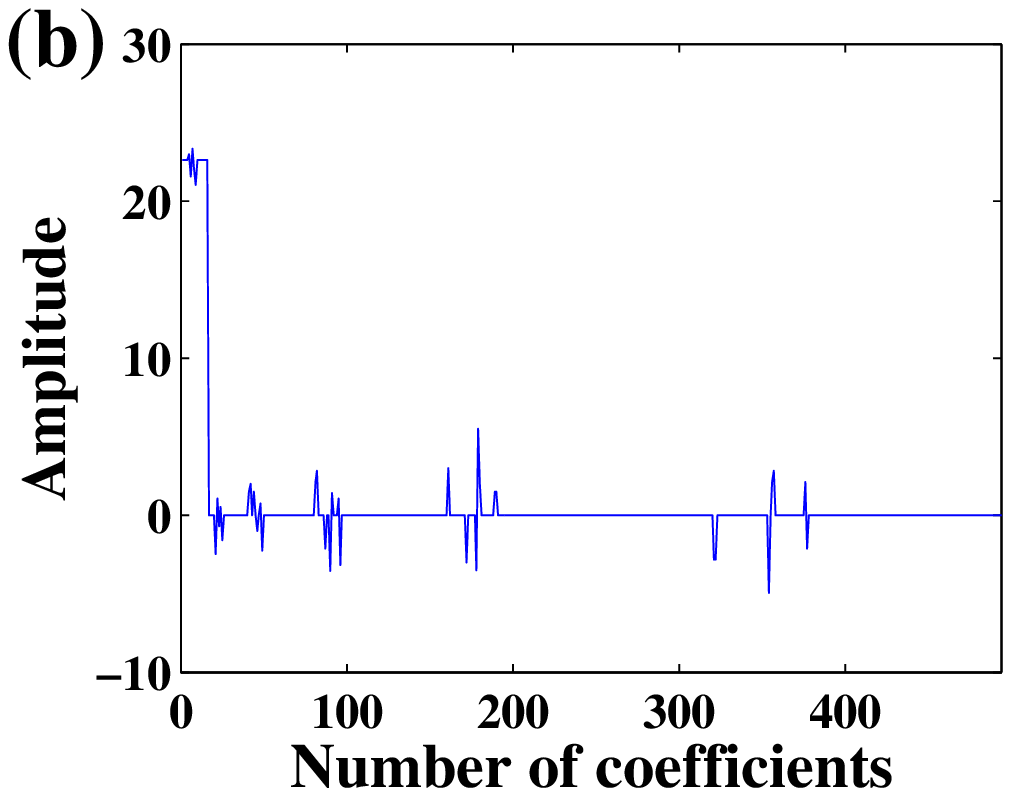}
\includegraphics[width=1.65in]{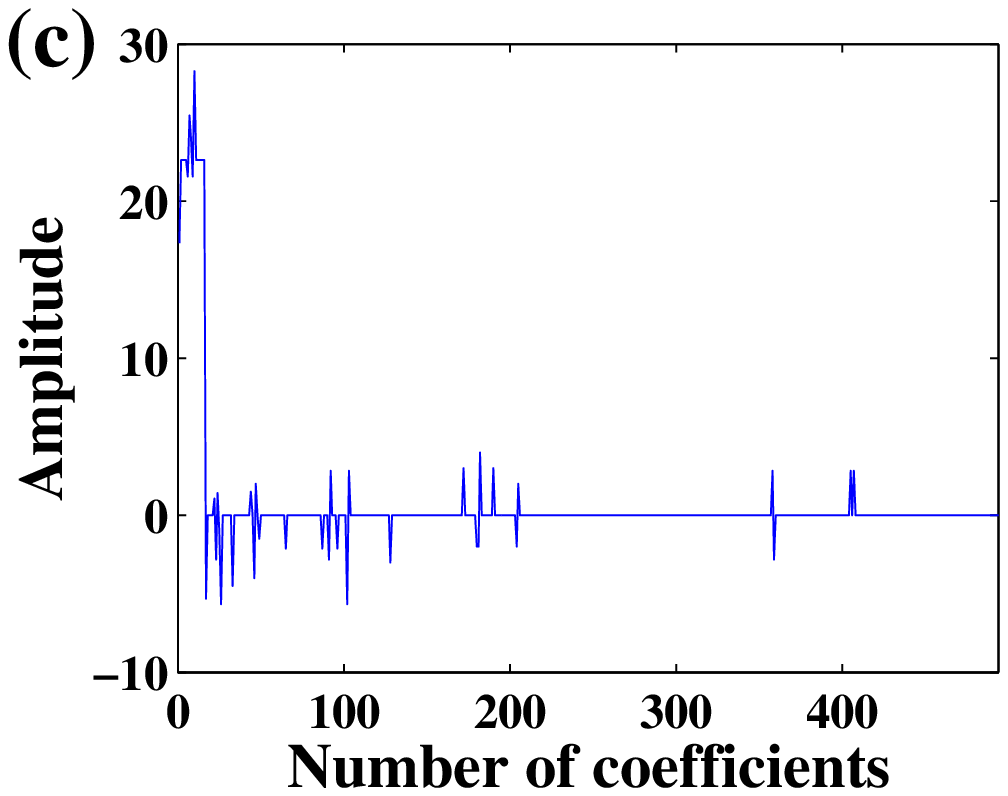}
\caption{Bar graph of the wavelet coefficients. (a) for Phantom A.  (b) for Phantom B.   (c) for Phantom C.}
\label{Figure 3}
\end{figure}

\begin{figure}[htbp]
\centering
\includegraphics[width=1.65in]{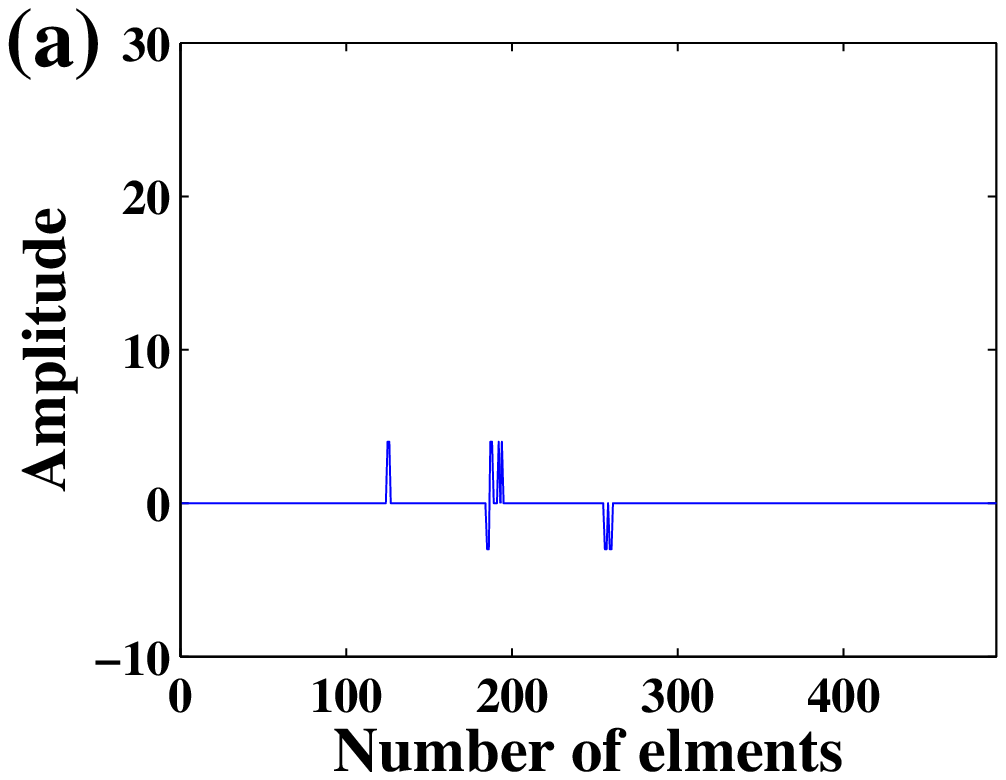}
\includegraphics[width=1.65in]{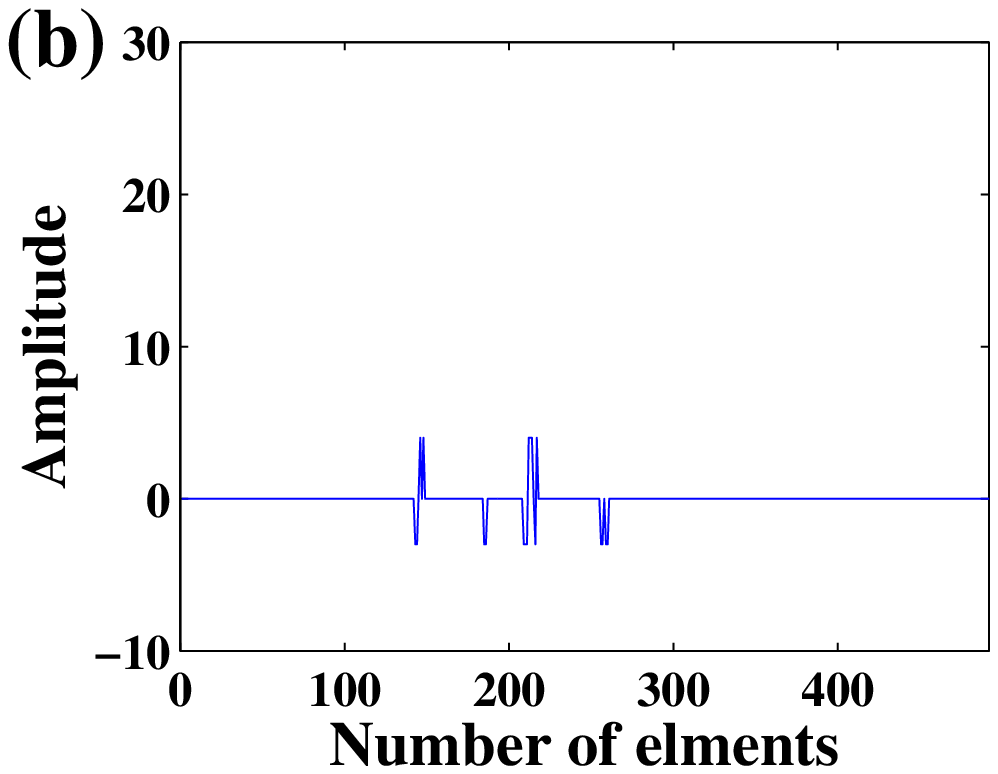}
\includegraphics[width=1.65in]{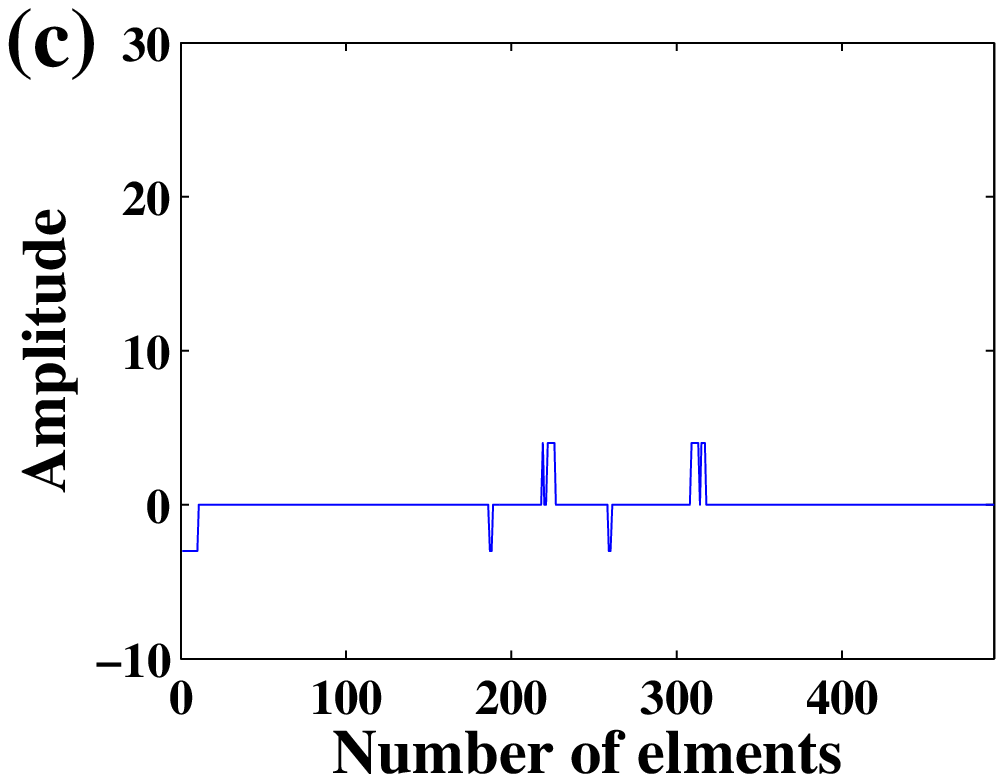}
\caption{Bar graph of the elements for inhomogeneities. (a) for Phantom A.  (b) for Phantom B.   (c) for Phantom C.}
\label{Figure 4}
\end{figure}

\subsection{Choices of parameters}
Selection of parameters is very important because an inappropriate parameter can lead to a meaningless result.
In this part, the strategies of choosing parameters for our approach are discussed.
We test it by the noise-free data later.
In addition, we need a distance function as an indicator of reconstruction error to compare the reconstruction performance more quantitatively.
The relative error criteria is used, i.e.
$$\textrm{RE}=\frac{\big{\|}\sigma-\sigma^{*}\big{\|^{2}}}{\big{\|}\sigma^{*}\big{\|}^{2}},$$
where $\sigma$ is the calculated conductivity and $\sigma^{*}$ is the true one.
Smaller error will give a better performance.

The regularization matrix (i.e. weighting operator) $R$ defines a suitable smoothing norm for the reconstruction problem, which is selected as an identity matrix in the standard form.
We here consider a version of 2D difference that includes smoothness assumptions.
The choice of regularization parameter $\alpha$ has attracted quite some interest.
The basic problem is that this parameter balances the goodness of misfit and the stability of the inversion.
In this paper, a geometric sequence $\alpha_{k}=q^{k}\alpha_{0}$ is used for the regularization parameter.
The choice of the initial value $\alpha_{0}$ and the factor $q$ depends on the inversion model and the noise level.
The Gauss-Newton iteration is stopped by the Morozov's discrepancy principle, i.e., $\big{\|}F(\sigma)-U^{\delta}\big{\|}<\big{\|}U-U^{\delta}\big{\|}$.
The inner loop is stopped after 10 iterations or when the update $\delta\sigma$ is less than $\tau=1\times10^{-2}$.

Additionally, the total iterations and the computational time for both Algorithm 1 and Algorithm 2 depend on the relaxation factor $\mu$.
A higher value for $\mu$ will lead to slower convergence of the split Bregman technique, since the threshold level $\eta$ in (\ref{sub2_solution}), which plays the key role in for promoting the sparsity of the reconstruction, is dependent on $\mu$.
Choosing a proper relaxation factor $\mu$ usually results to reasonable convergence of both CG and the split Bregman method.
In addition, it can be seen from the proposed inversion model (\ref{L1L2}) that parameter $\beta$ is the key parameter, which balances the influence of $l_{1}$-norm and $l_{2}$-norm.
We find that the effect of $l_{2}$-norm penalty is weakened gradually but the effect of $l_{1}$-norm penalty is strengthened gradually along with the decrease of $\beta$.
We will discuss the influence of parameter $\mu$ and $\beta$ for the reconstruction performance later.

\subsection{Noise-free case}
The numerical performance of Algorithm 1 is first tested by the noise-free data.
Table 1 lists the selected algorithmic parameters for the test phantoms.
Fixing $\beta=0.1$, we first investigate how the relaxation factor $\mu$ affects the performance of Algorithm 1.
Figure 5 shows the reconstruction errors curves against iterations with various parameters $\mu$.
The results from Figure 5 indicate that it can result in good convergence with appropriate parameter $\mu$.
We here take $\mu=1\times10^{-10}$.

\begin{table}[htbp]\footnotesize
\centering
\caption{Parameters for Algorithm 1 with noise-free case.}
\centering
\begin{tabular*}{11cm}{@{\extracolsep{\fill}}cccc}
\hline
Parameters & Phantom A & Phantom B  & Phantom C\\
\hline
$\alpha_{0}$  & $1\times10^{-6}$   & $1\times10^{-6}$  & $1\times10^{-7}$\\
$q_{\alpha}$  & 0.6   & 0.8  & 0.5\\
\hline
\end{tabular*}
\end{table}

\begin{figure}[htbp]
\centering
\includegraphics[width=1.65in]{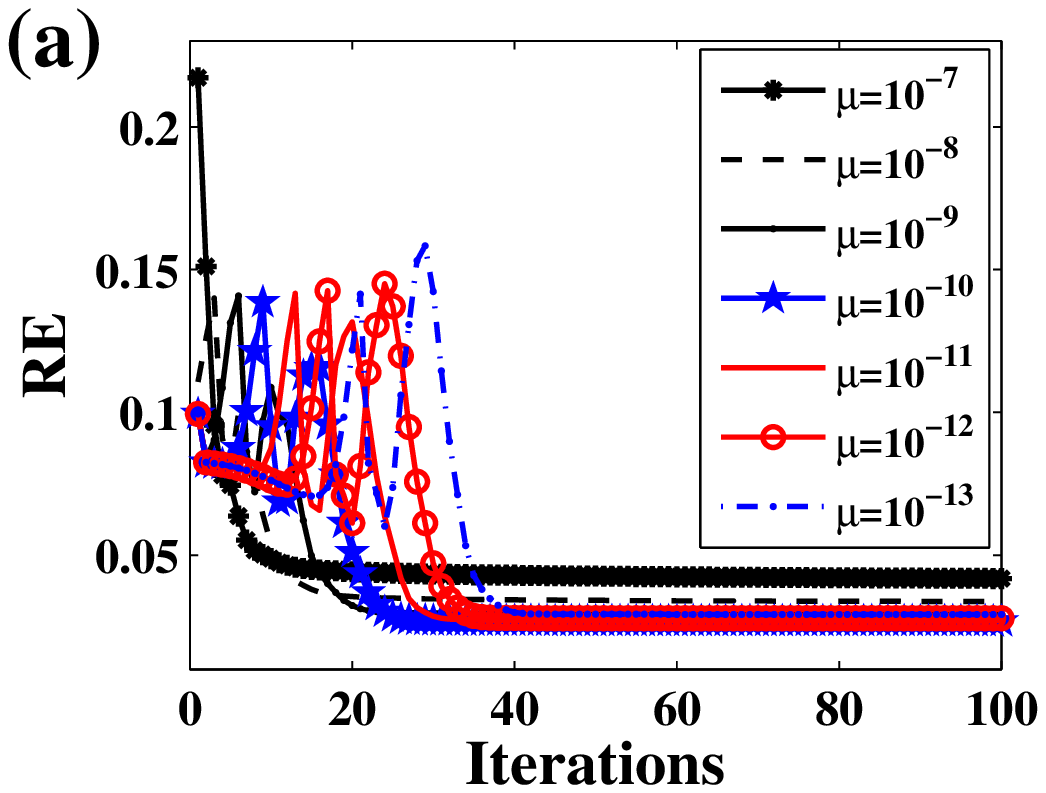}
\includegraphics[width=1.65in]{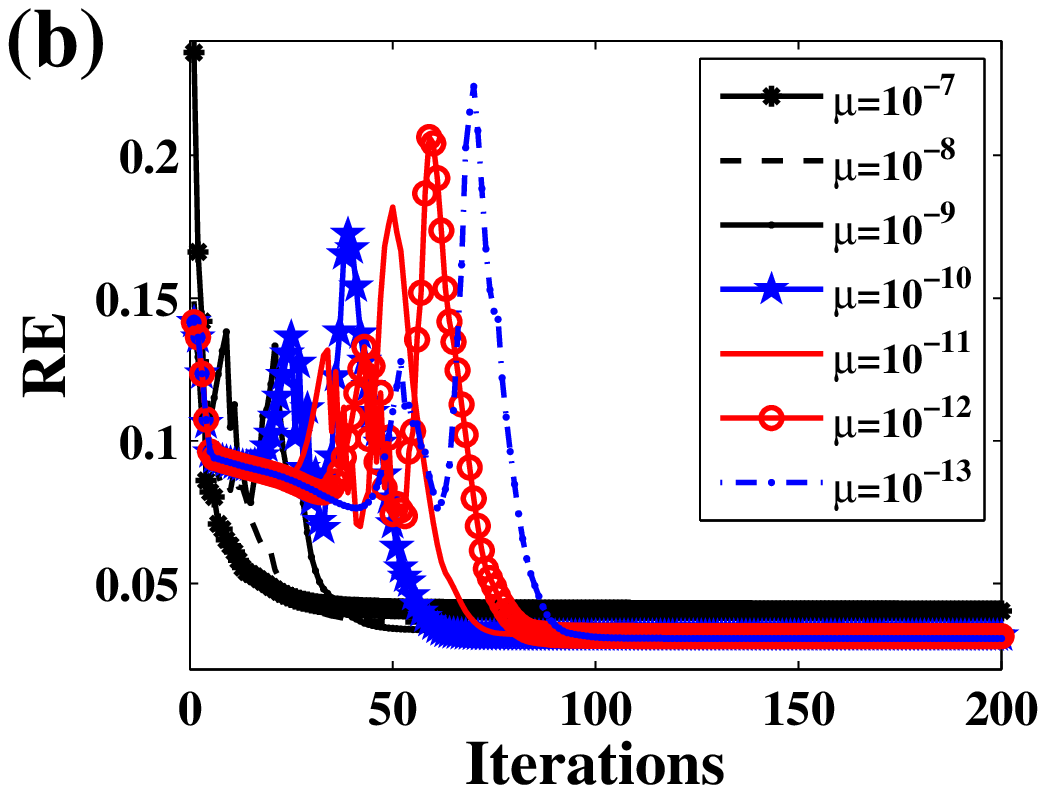}
\includegraphics[width=1.65in]{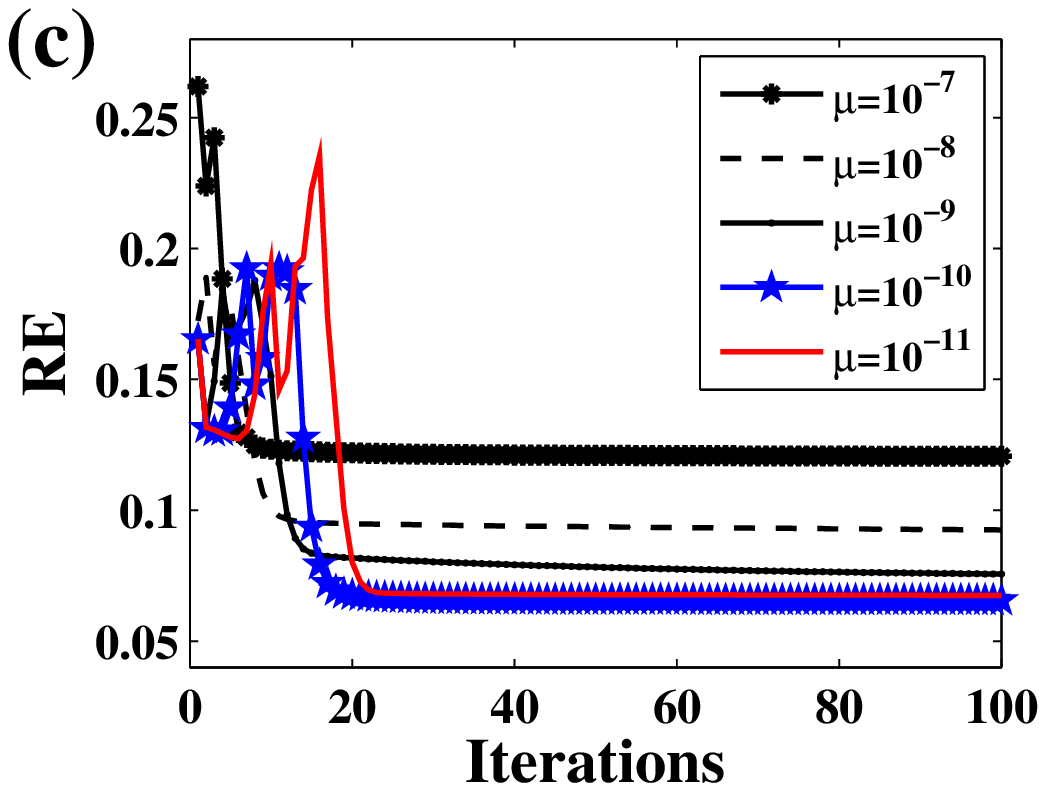}
\caption{Reconstruction errors against iterations with various relaxation factors $\mu$ as fixing $\beta=0.1$. (a) for Phantom A.  (b) for Phantom B.   (c) for Phantom C.}
\label{Figure 5}
\end{figure}

We next investigate the influence of parameter $\beta$ for the regularizer performance.
Figure 6(a)-(c) illustrates RE curves with different parameter $\beta$.
It can be found from Figure 6(a)-(c) that Algorithm 1 is convergent.
Figure 6(d) depicts the change of RE along with $\beta$.
It indicates that the smaller parameter $\beta$ is the less RE for phantoms with the noise-free data considered in this work.
Figure 7(a)-(c) shows reconstructed images with different $\beta$ for Phantom A, Phantom B, and Phantom C, respectively.
It is clearly visible that the smaller the parameter $\beta$ is the clearer for the edges of the inclusions.
And the artifacts can be effectively decreased with smaller $\beta$.
It further verifies that $l_{1}$-norm penalty can improve the imaging quality.
The corresponding errors for reconstructions in Figure 7(a)-(c) are listed in Table 2.
As can be expected, numerical results above indicate that Algorithm 1 shows a satisfactory numerical performance.
Let background $\sigma_0$ be known.
Simulations for Algorithm 2 will yield similar results, so here not listed in detail.
We then in Figure 8 list the reconstructed images with both Algorithm 1 and Algorithm 2.
And corresponding errors for reconstructions displayed in Figure 8 sre summarized in Table 3.
It can be seen form Figure 8 and Table 3 that Algorithm 2 is superior to Algorithm 1 with shaper edges and less artifacts.
As we are expected, Algorithm 2 obtains better results because it can be found from Figure 3 and Figure 4 that Figure 4 used for Algorithm 2 shows better sparsity apparently.

\begin{figure}[htbp]
\centering
\includegraphics[width=1.8in]{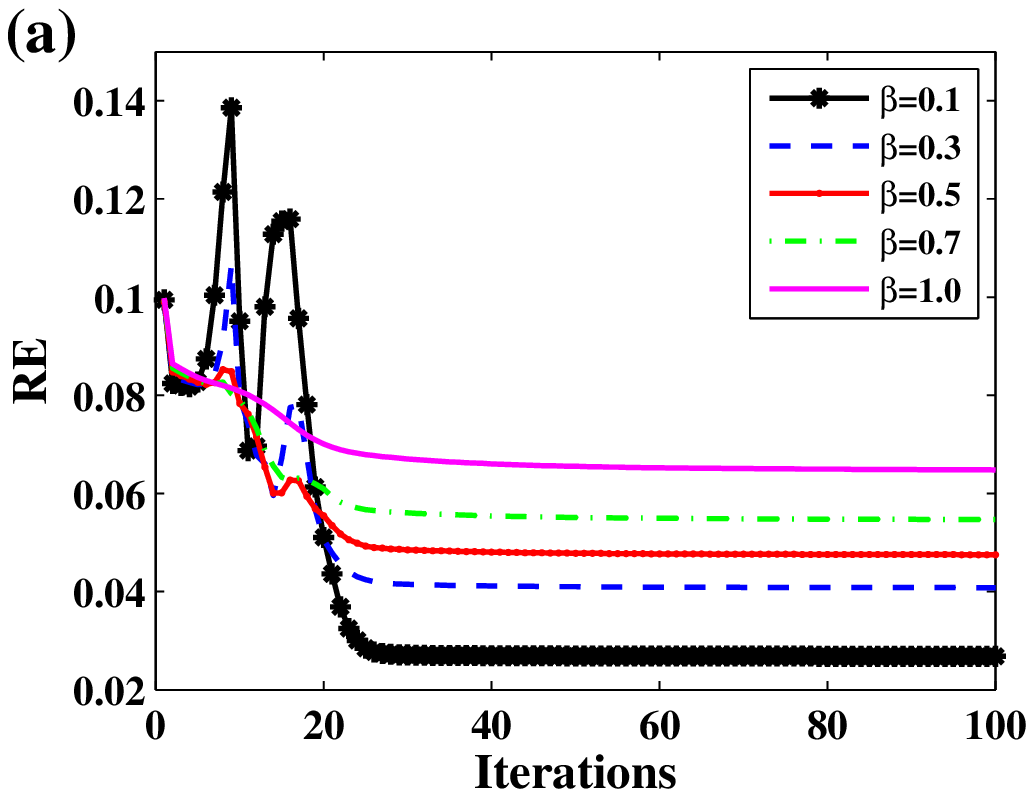}
\includegraphics[width=1.8in]{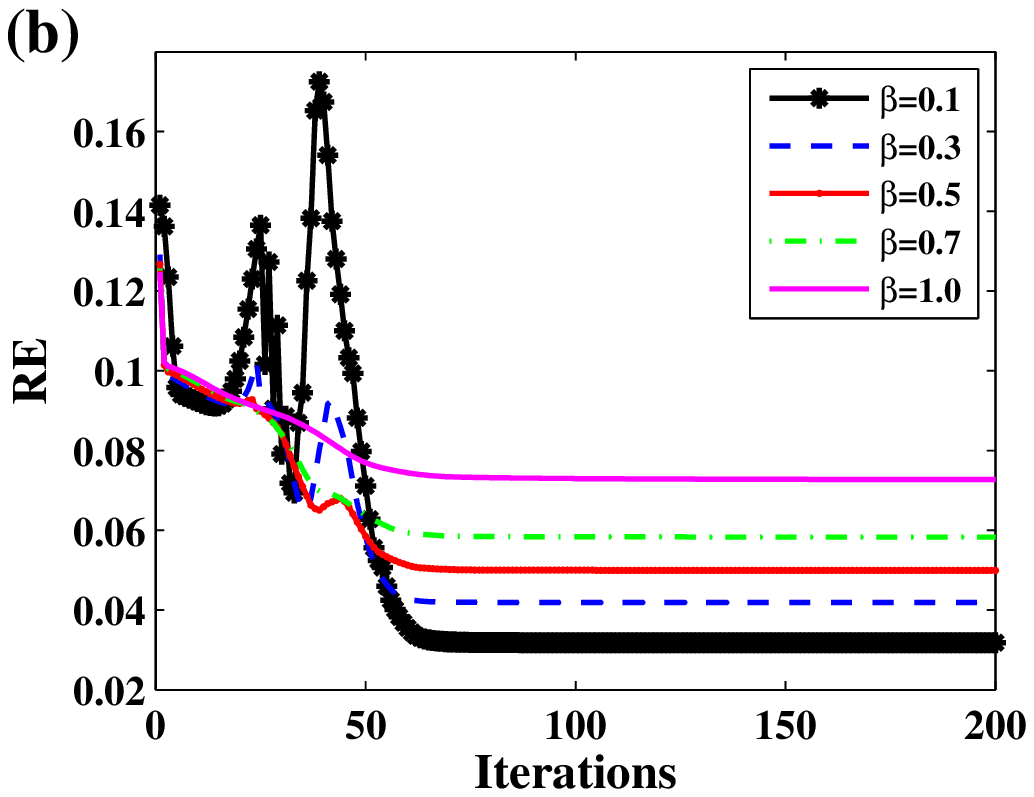}
\includegraphics[width=1.8in]{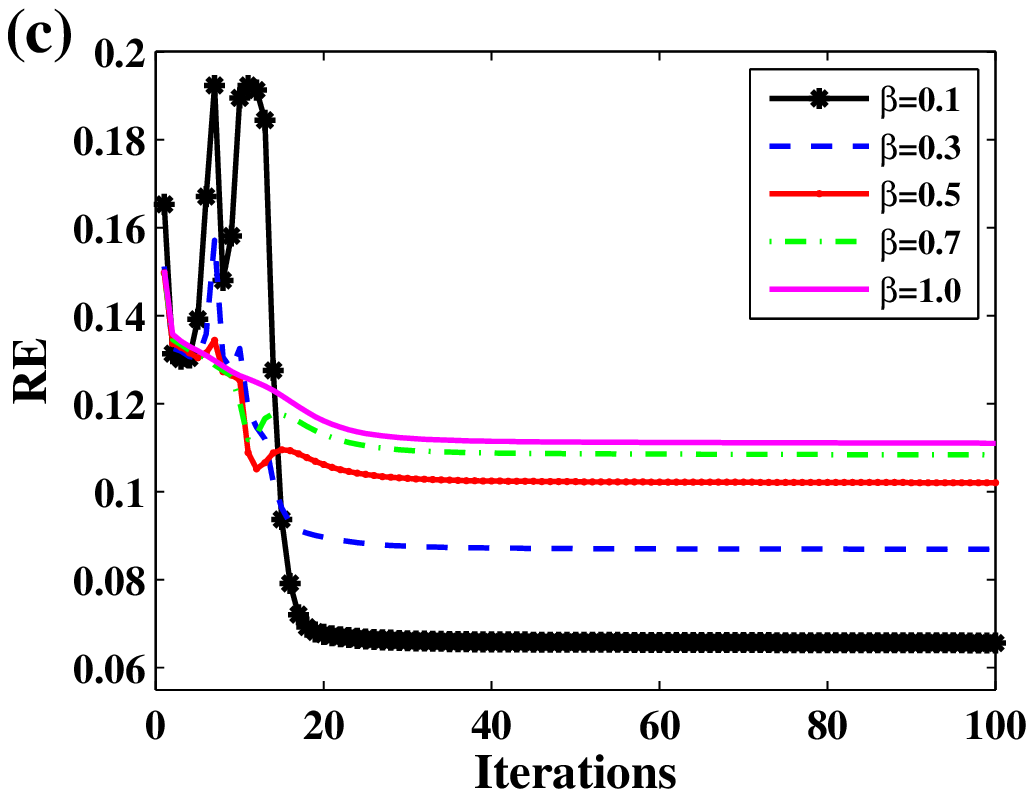}
\includegraphics[width=1.8in]{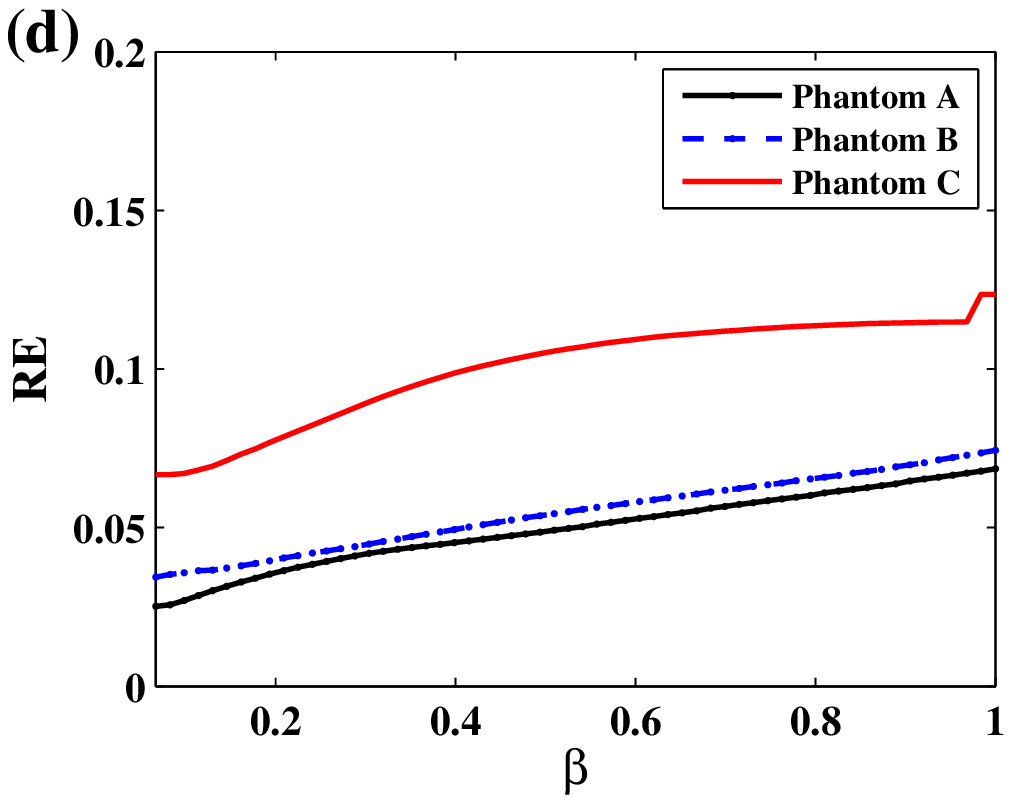}
\caption{Reconstructed image errors against iterations with different parameters $\beta$ for (a) Phantom A, (b) Phantom B, (c) Phantom C, respectively.
 (d) Relative errors against parameters $\beta$.}
\label{Figure 6}
\end{figure}

\begin{figure}[htbp]
\centering
\begin{minipage}[b]{0.18\textwidth}
\centering
\includegraphics[width=1.15in]{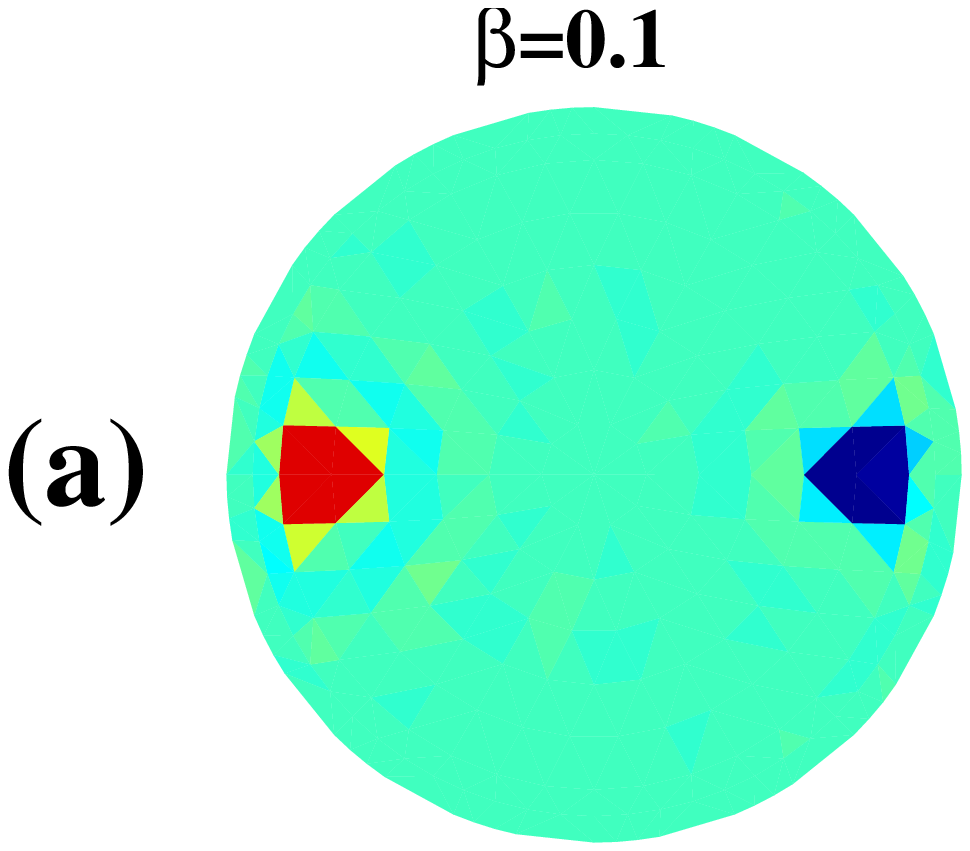}
\end{minipage}
\begin{minipage}[b]{0.18\textwidth}
\centering
\includegraphics[width=1.15in]{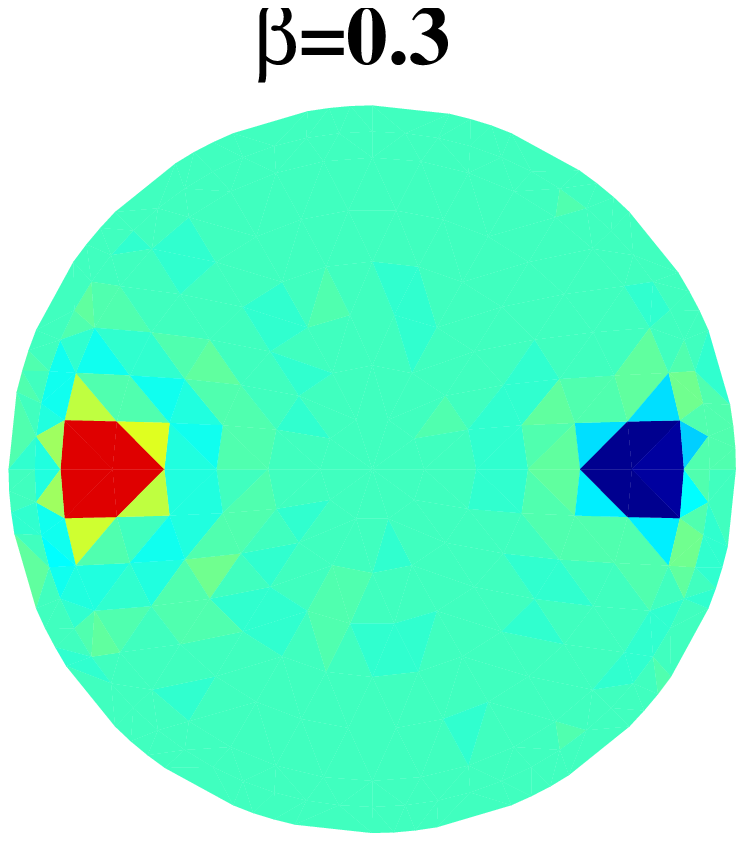}
\end{minipage}
\begin{minipage}[b]{0.18\textwidth}
\centering
\includegraphics[width=1.15in]{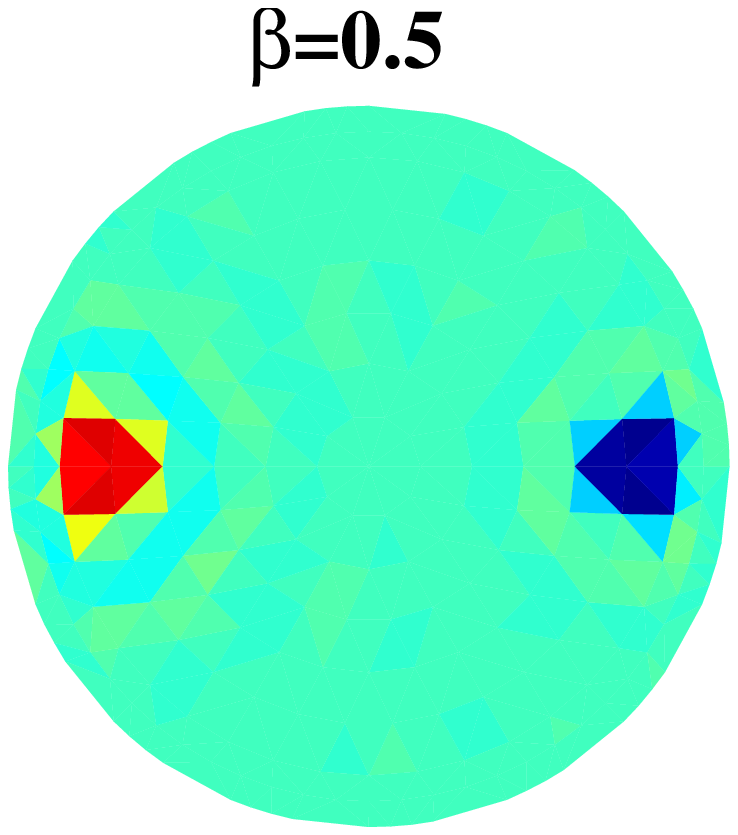}
\end{minipage}
\begin{minipage}[b]{0.18\textwidth}
\centering
\includegraphics[width=1.15in]{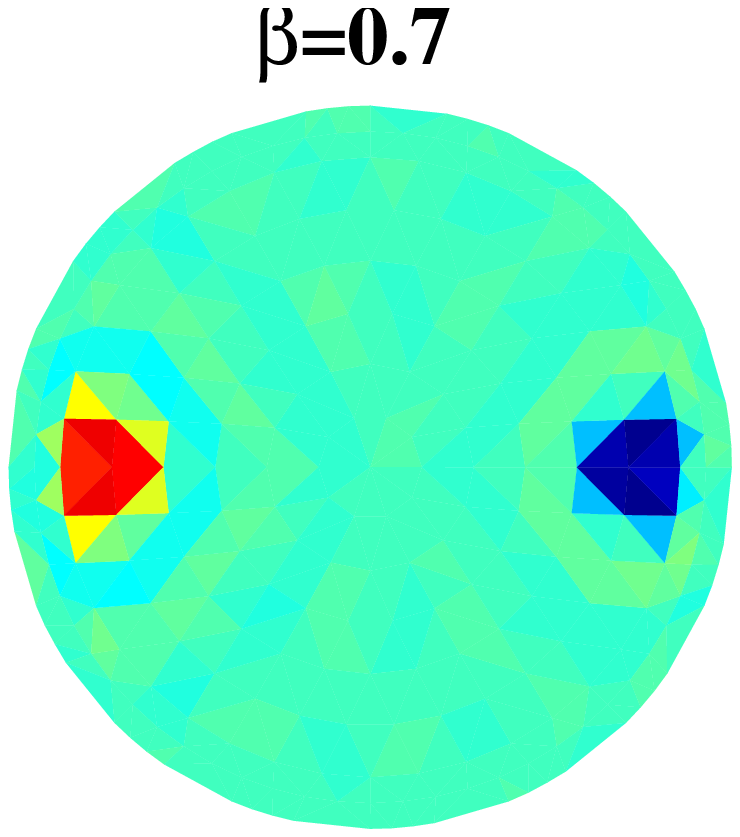}
\end{minipage}
\begin{minipage}[b]{0.18\textwidth}
\centering
\includegraphics[width=1.15in]{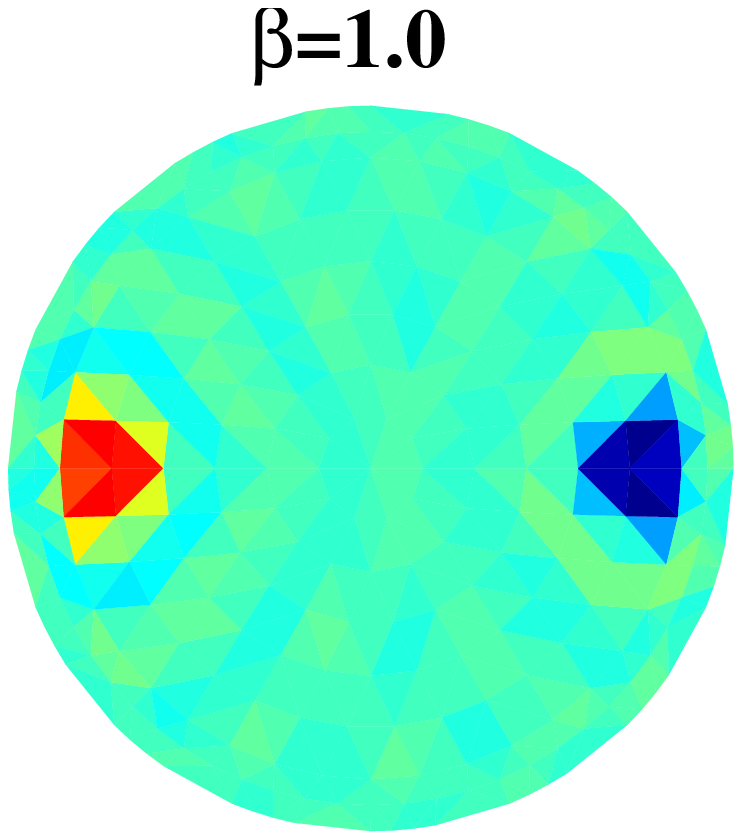}
\end{minipage}
\begin{minipage}[b]{0.18\textwidth}
\centering
\includegraphics[width=1.15in]{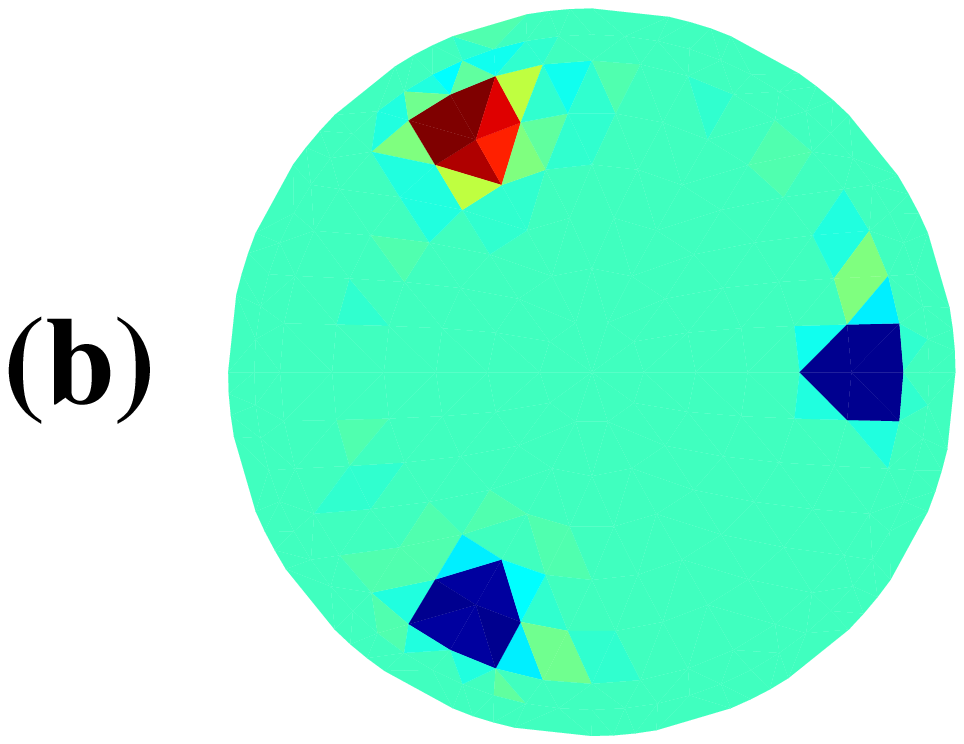}
\end{minipage}
\begin{minipage}[b]{0.18\textwidth}
\centering
\includegraphics[width=1.15in]{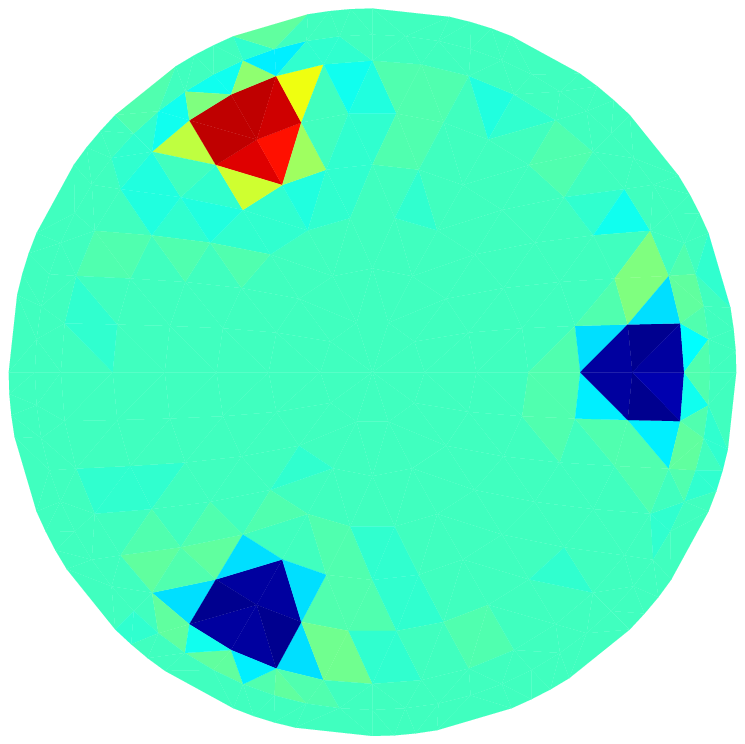}
\end{minipage}
\begin{minipage}[b]{0.18\textwidth}
\centering
\includegraphics[width=1.15in]{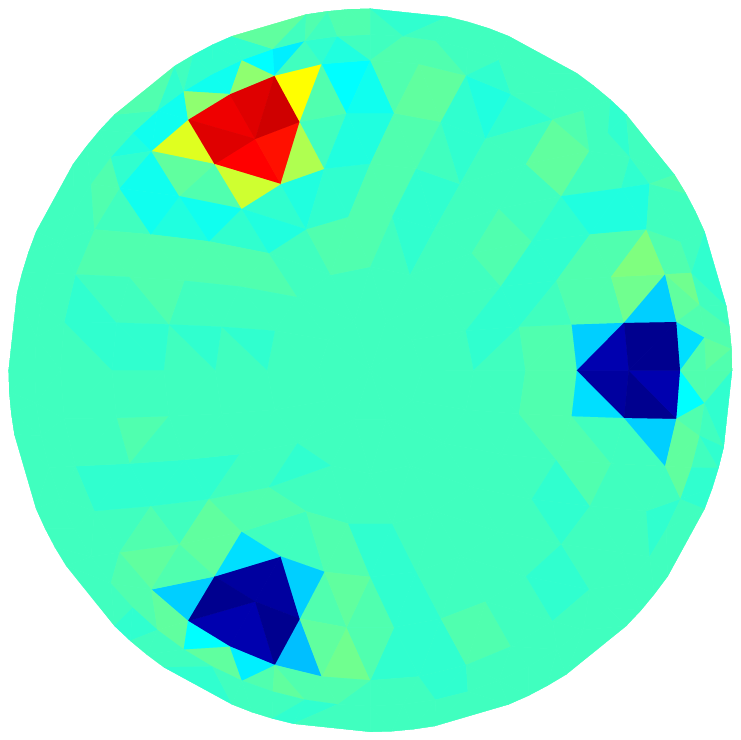}
\end{minipage}
\begin{minipage}[b]{0.18\textwidth}
\centering
\includegraphics[width=1.15in]{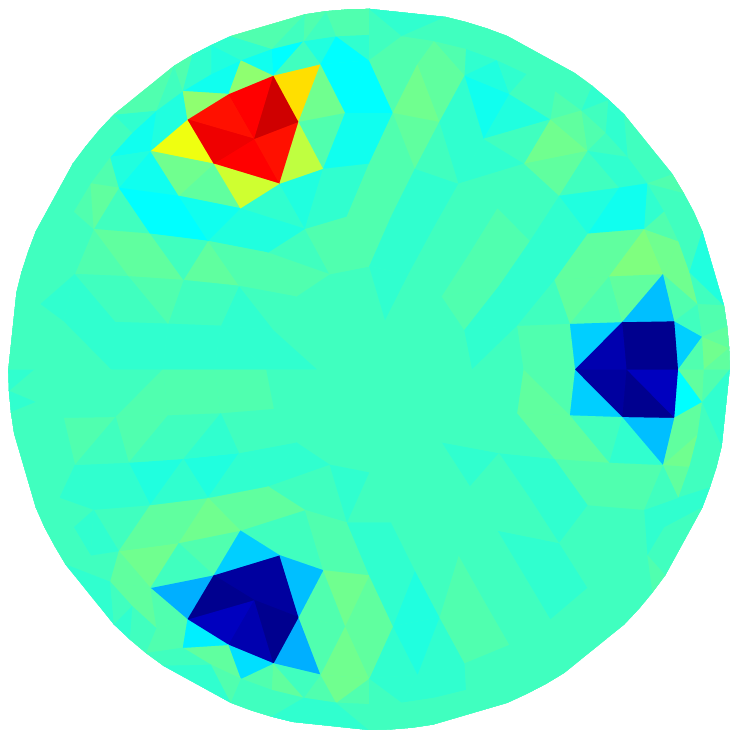}
\end{minipage}
\begin{minipage}[b]{0.18\textwidth}
\centering
\includegraphics[width=1.15in]{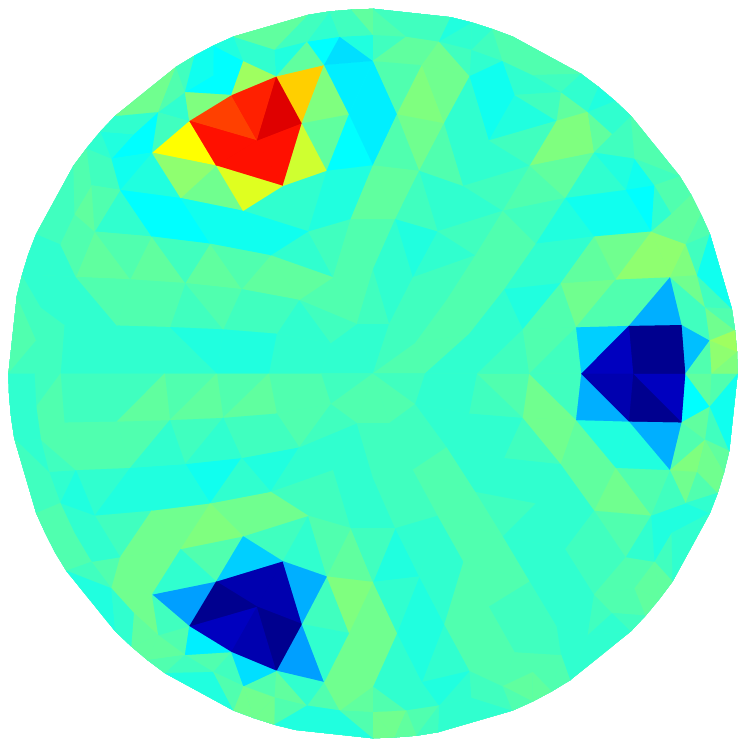}
\end{minipage}
\begin{minipage}[b]{0.18\textwidth}
\centering
\includegraphics[width=1.15in]{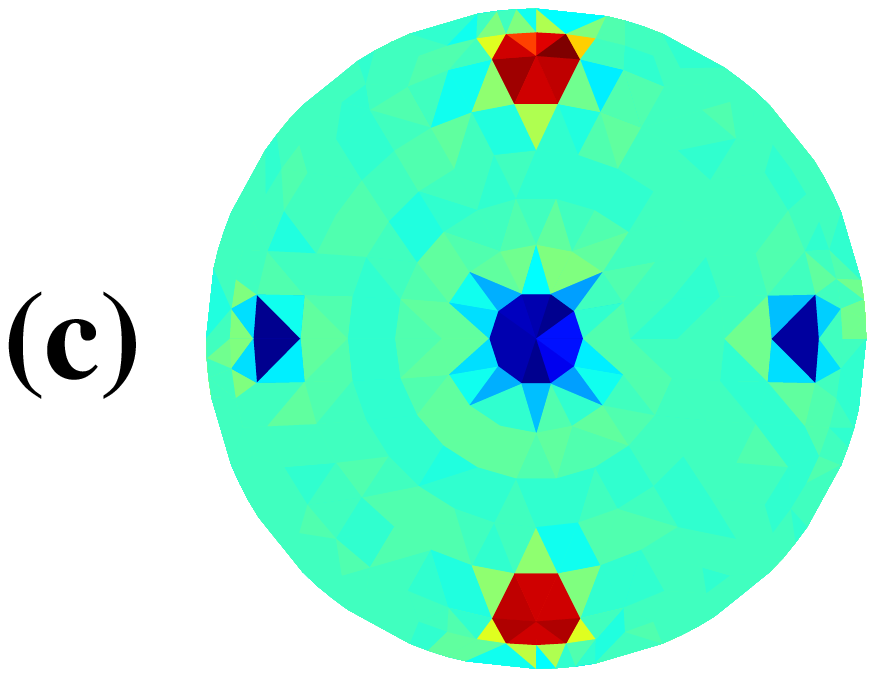}
\end{minipage}
\begin{minipage}[b]{0.18\textwidth}
\centering
\includegraphics[width=1.15in]{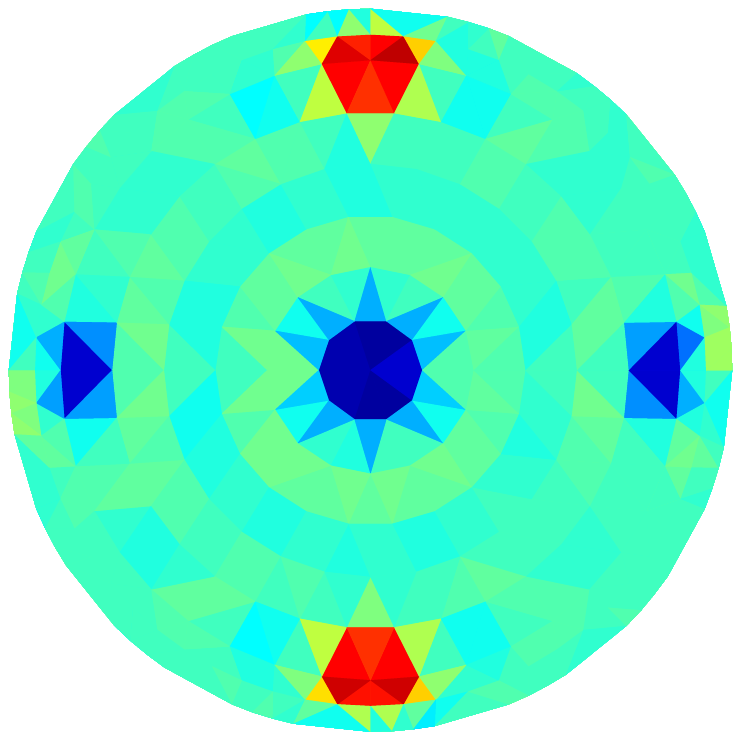}
\end{minipage}
\begin{minipage}[b]{0.18\textwidth}
\centering
\includegraphics[width=1.15in]{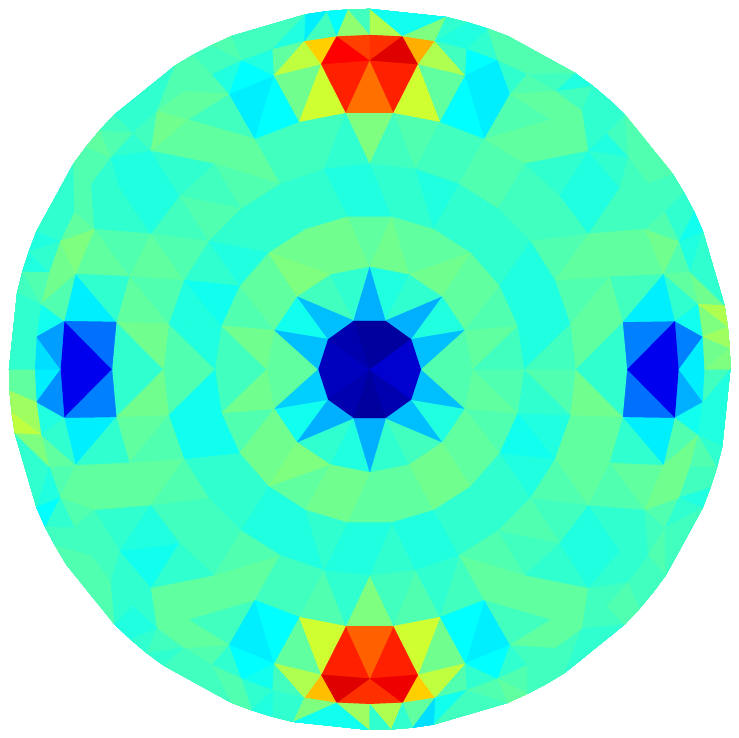}
\end{minipage}
\begin{minipage}[b]{0.18\textwidth}
\centering
\includegraphics[width=1.15in]{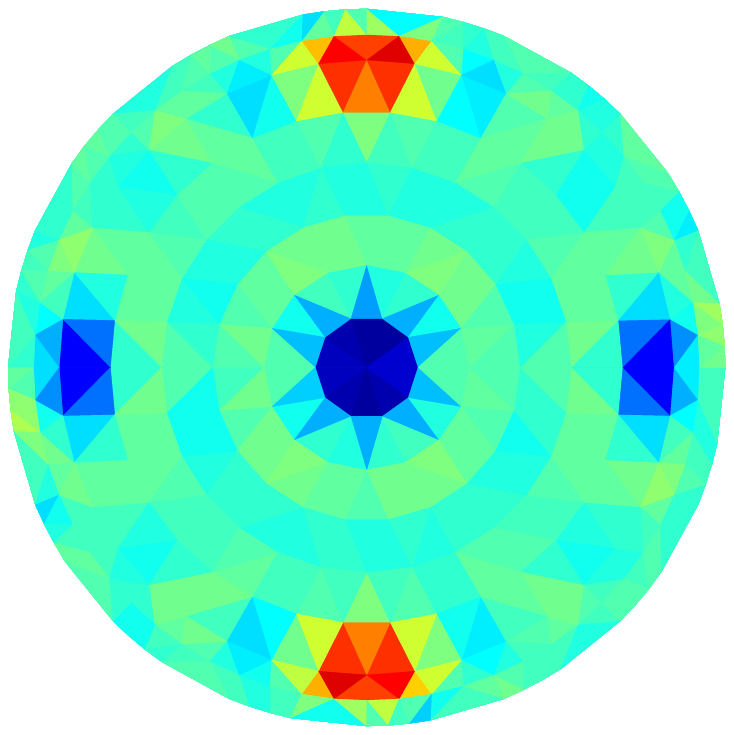}
\end{minipage}
\begin{minipage}[b]{0.18\textwidth}
\centering
\includegraphics[width=1.15in]{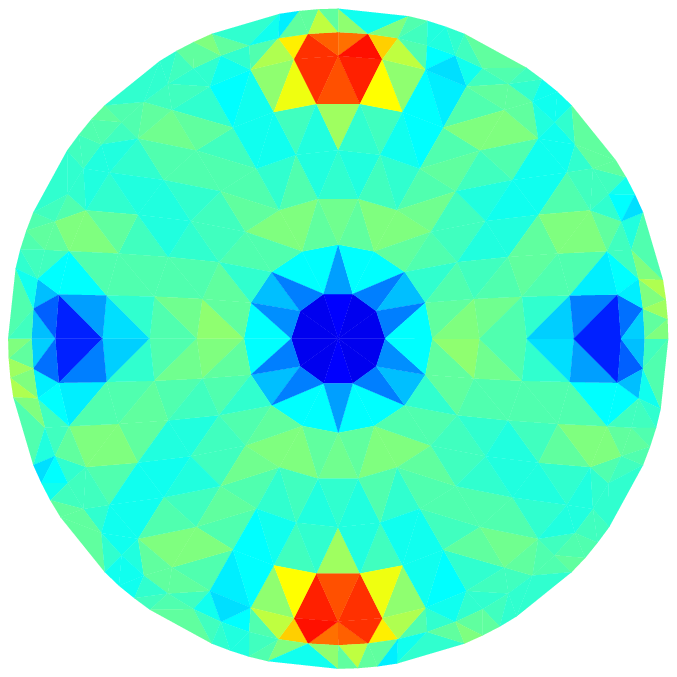}
\end{minipage}
\caption{Reconstructions with noise-free case at various parameters $\beta$ (from left column to right column with $\beta=0.1$, 0.3, 0.5, 0.7, 1.0, respectively). (All displayed by the unified colorbars from 1 to 8).}
\label{Figure 7}
\end{figure}

\begin{table}[htbp]\footnotesize
\centering
\caption{Comparison of image errors for the results displayed in Figure 7.}
\centering
\begin{tabular*}{12cm}{@{\extracolsep{\fill}}cccccc}
\hline
Parameters & $\beta=0.1$ & $\beta=0.3$   & $\beta=0.5$ &$\beta=0.7$  &$\beta=1.0$  \\
\hline
Phantom A  & 0.0272   & 0.0417  & 0.0488  & 0.0567 &0.0685\\
Phantom B  & 0.0324   & 0.0424  & 0.0500  & 0.0659 &0.0884\\
Phantom C  & 0.0670   & 0.0893  & 0.1056  & 0.1123 &0.1240\\
\hline
\end{tabular*}
\end{table}

\begin{figure}[htbp]
\centering
\begin{minipage}[b]{0.25\textwidth}
\centering
\includegraphics[width=1.5in]{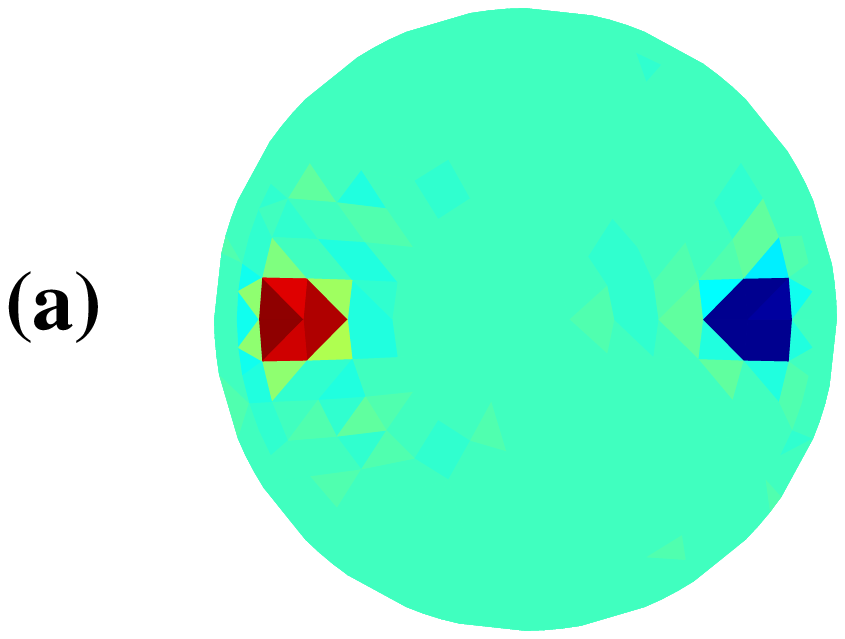}
\end{minipage}
\begin{minipage}[b]{0.25\textwidth}
\centering
\includegraphics[width=1.5in]{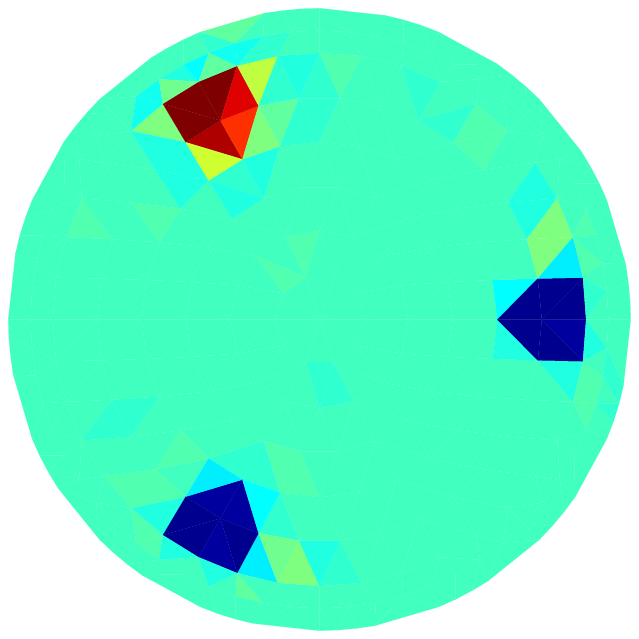}
\end{minipage}
\begin{minipage}[b]{0.25\textwidth}
\centering
\includegraphics[width=1.5in]{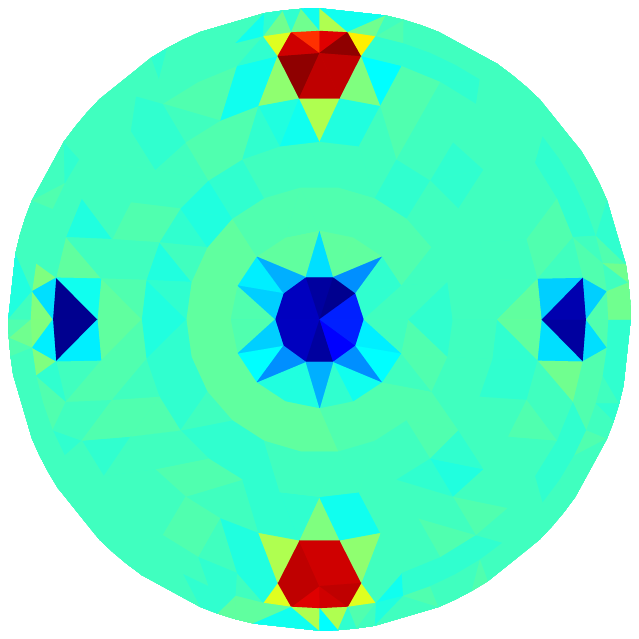}
\end{minipage}
\begin{minipage}[b]{0.25\textwidth}
\centering
\includegraphics[width=1.5in]{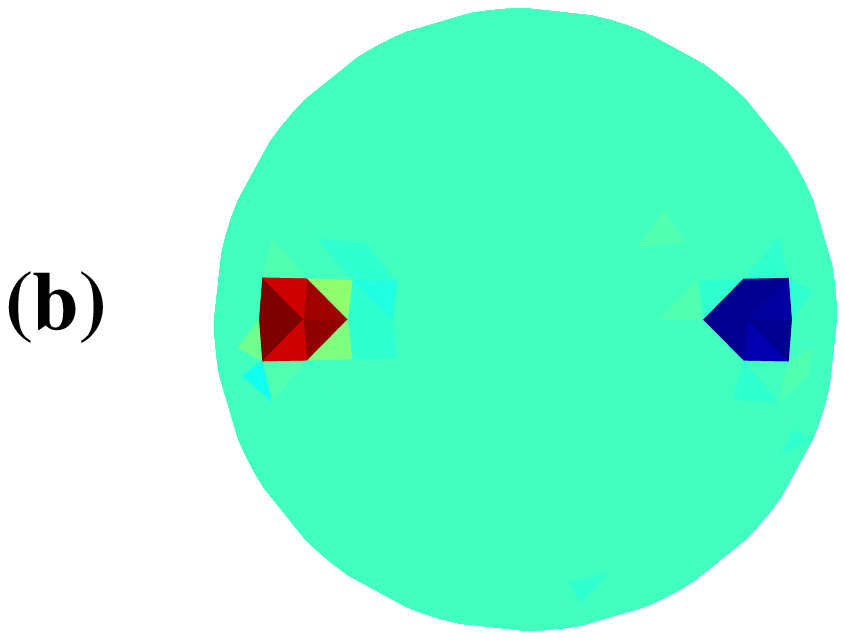}
\end{minipage}
\begin{minipage}[b]{0.25\textwidth}
\centering
\includegraphics[width=1.5in]{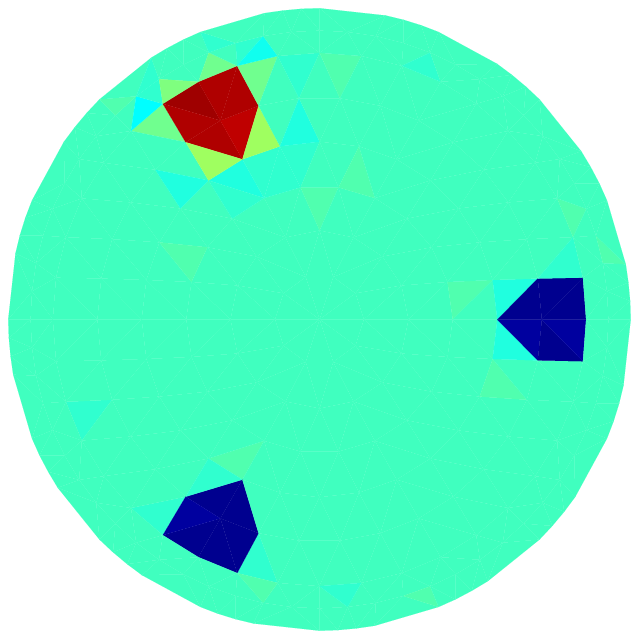}
\end{minipage}
\begin{minipage}[b]{0.25\textwidth}
\centering
\includegraphics[width=1.5in]{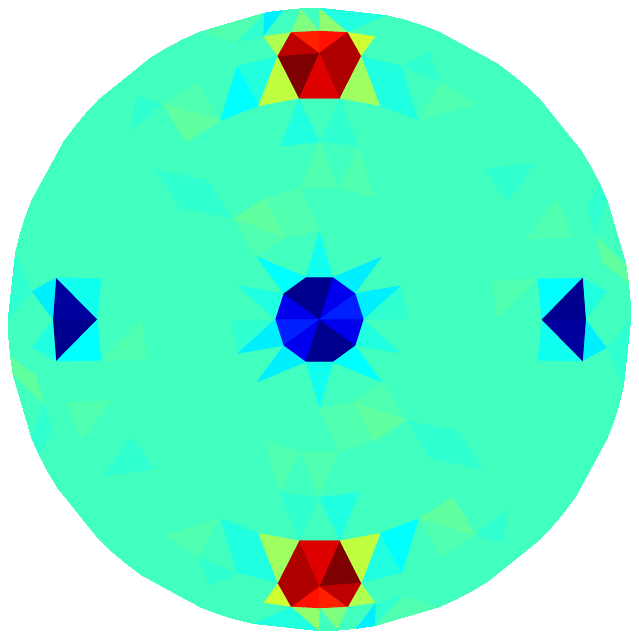}
\end{minipage}
\caption{Reconstructions with noise-free case at $\beta=0.1$. (a) by Algorithm 1; (b) by Algorithm 2. (All displayed by the unified colorbars from 1 to 8).}
\label{Figure 8}
\end{figure}

\begin{table}[htbp]\footnotesize
\centering
\caption{Comparison of image errors for the results displayed in Figure 8.}
\centering
\begin{tabular*}{13cm}{@{\extracolsep{\fill}}cccc}
\hline
Phantoms & Phantom A  & Phantom B  & Phantom C\\
\hline
Algorithm 1 &0.0272& 0.0324&0.0670\\
Algorithm 2 &0.0156& 0.0195&0.0535\\
\hline
\end{tabular*}
\end{table}

\subsection{Noise-contaminated case}
EIT image reconstruction is a typical ill-posed problem.
In the data acquisition processes, the measurement data always contain noise, and the noise will influence the quality of the reconstruction image.
The robustness of reconstruction technique with respect to noise is crucial.
Usually, the synthetic Gaussian noise was generated to simulate systematic and random errors existing in the data acquisition processes.
The measurement data with noise was indicated as $$U^{\delta}=U^{syn}(\sigma^{*})+\varepsilon\cdot n,$$
where $\varepsilon$ is the noise level, and $n$ is a noise vector obeying a Gaussian distribution with zero mean and standard deviation of 1.
Due to the high sensitivity of EIT problem to noise, we here consider two additive noise levels with $\varepsilon=0.1\%,0.3\%$, respectively.

We next investigate the influence of parameter $\beta$ for the noise-contaminated case.
The selected algorithmic parameters are listed in Table 4.
Figure 9 and  10 depict the change of reconstruction errors along with $\beta$ under different noise levels by Algorithm 1 and Algorithm 2, respectively.
It is found from Figure 9 and 10 that relative errors decrease at first and then increase as $\beta$ changes from small to large value.
Due to the presence of noise, there exists an trade-off between $l_{1}$-norm and $l_{2}$-norm.
It can be seen that $\beta$ shouldn't be too large and there exists optimum values, which may range from 0.05 to 0.5.
That is to say, the optimum for $\beta$ fall in between 0.05 and 0.5, which will lead to good performance and better robustness to noise for the test cases.
Furthermore, it can be found that $l_{2}$-norm is essential for the regularity to this inverse problem from comparing Figure 6(d) with Figure 9.

\begin{table}[htbp]\footnotesize
\centering
\caption{Algorithmic parameters with noise-contaminated case.}
\centering
\begin{tabular*}{13.5cm}{@{\extracolsep{\fill}}cccccccc}
\hline
\multirow{2}{*}{}  &\multicolumn{3}{c}{Algorithm 1} & & \multicolumn{3}{c}{Algorithm 2} \\
\cline{2-4}\cline{6-8}
Noise $\varepsilon$ &$\alpha_{0}$  & $q_{\alpha}$  & $\mu$  && $\alpha_{0}$  & $q_{\alpha}$  & $\mu$  \\
\hline
$0.1\%$ &$1\times10^{-4}$ & 0.6  & $1\times10^{-6}$ &&$1\times10^{-5}$ & 0.6  & $1\times10^{-6}$ \\
$0.3\%$ &$1\times10^{-3}$ & 0.6  & $1\times10^{-5}$ &&$1\times10^{-3}$ & 0.6  & $1\times10^{-4}$\\
\hline
\end{tabular*}
\end{table}

\begin{figure}[htbp]
\centering
\includegraphics[width=1.5in]{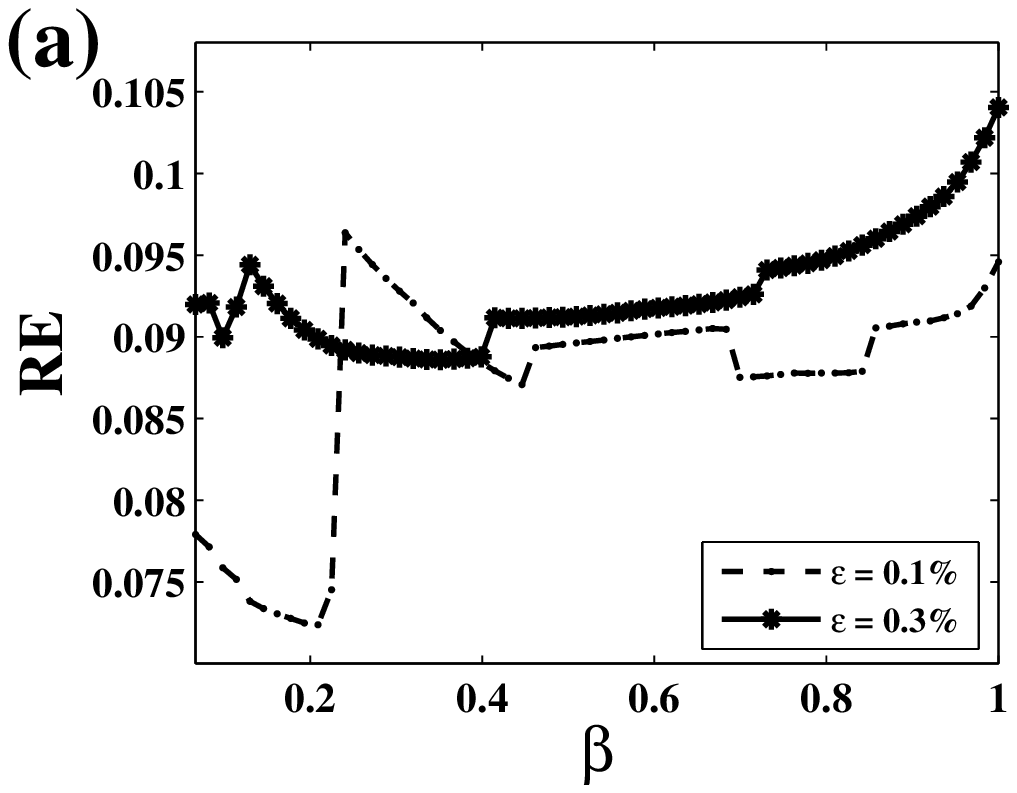}
\includegraphics[width=1.5in]{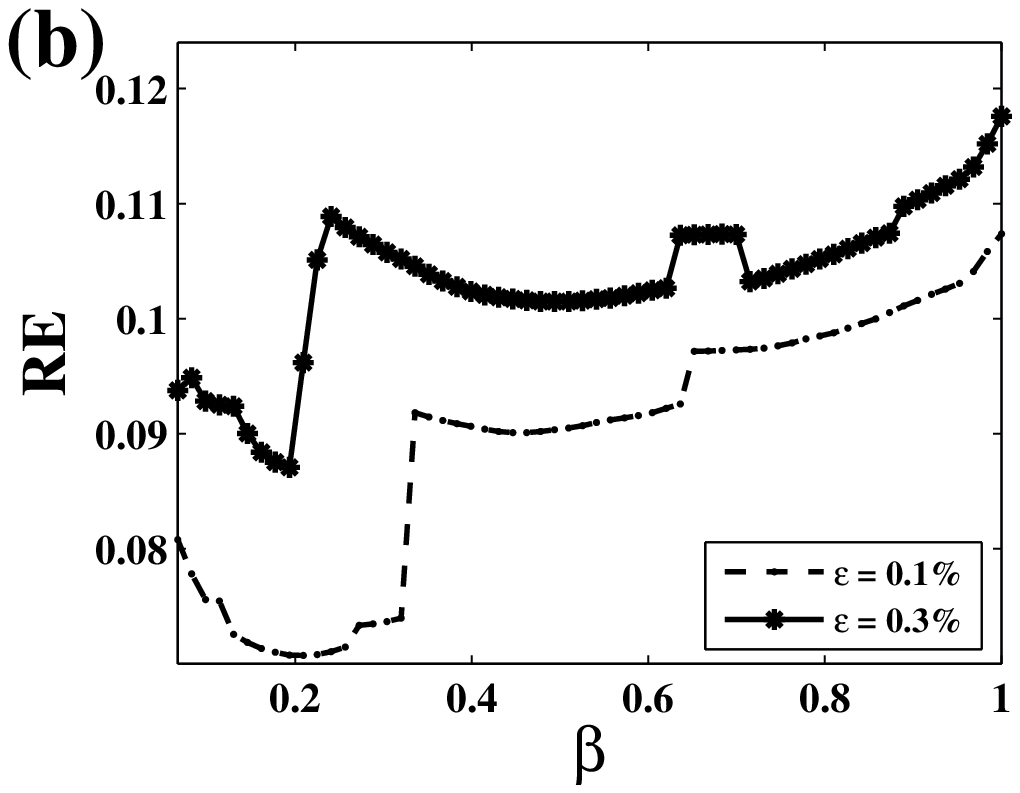}
\includegraphics[width=1.5in]{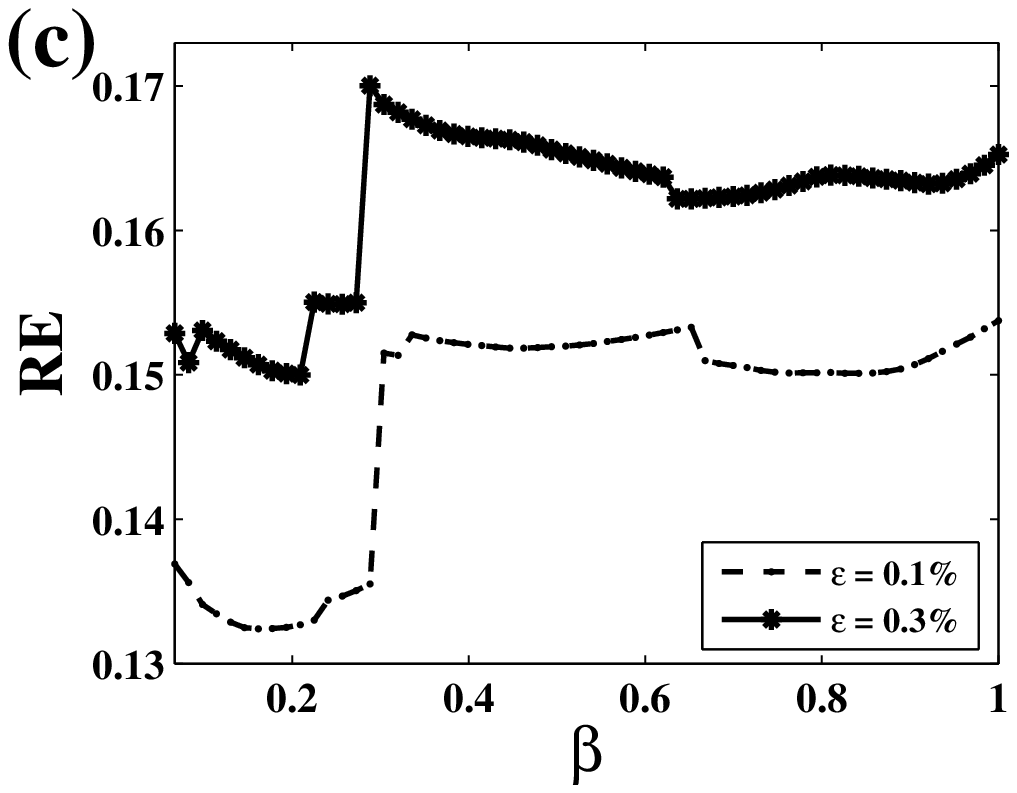}
\caption{Relative errors against parameter $\beta$ under different noise levels by Algorithm 1. (a) for Phantom A.  (b) for Phantom B.   (c) for Phantom C.}
\label{Figure 9}
\end{figure}

\begin{figure}[htbp]
\centering
\includegraphics[width=1.5in]{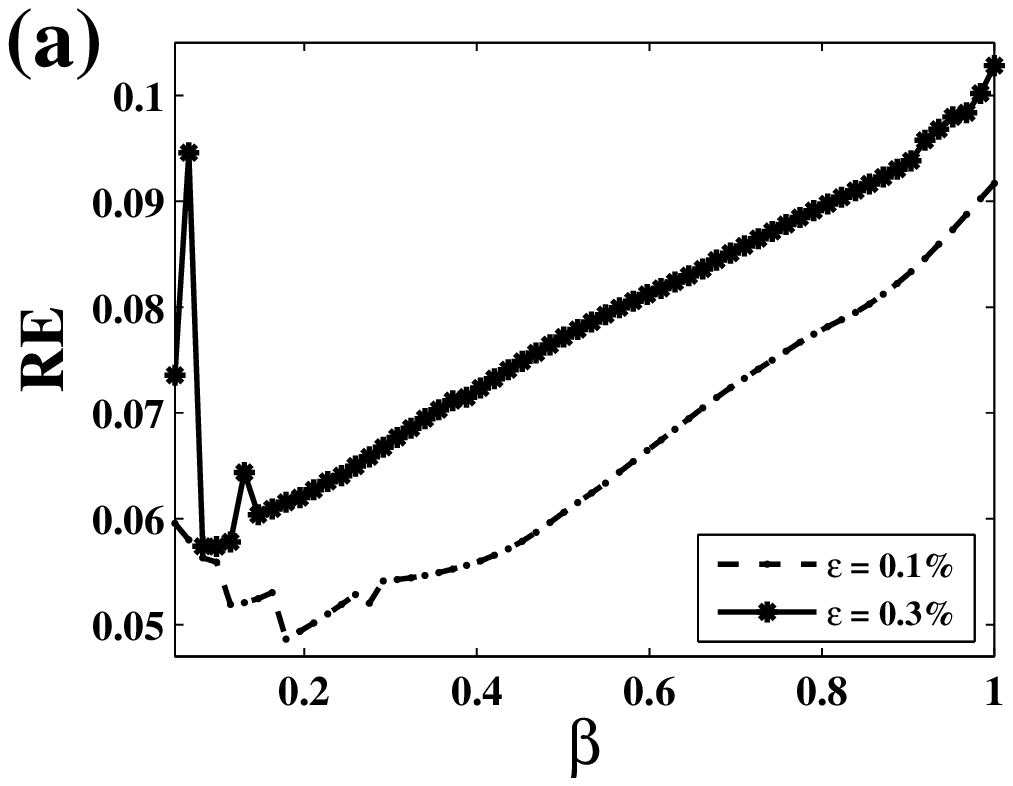}
\includegraphics[width=1.5in]{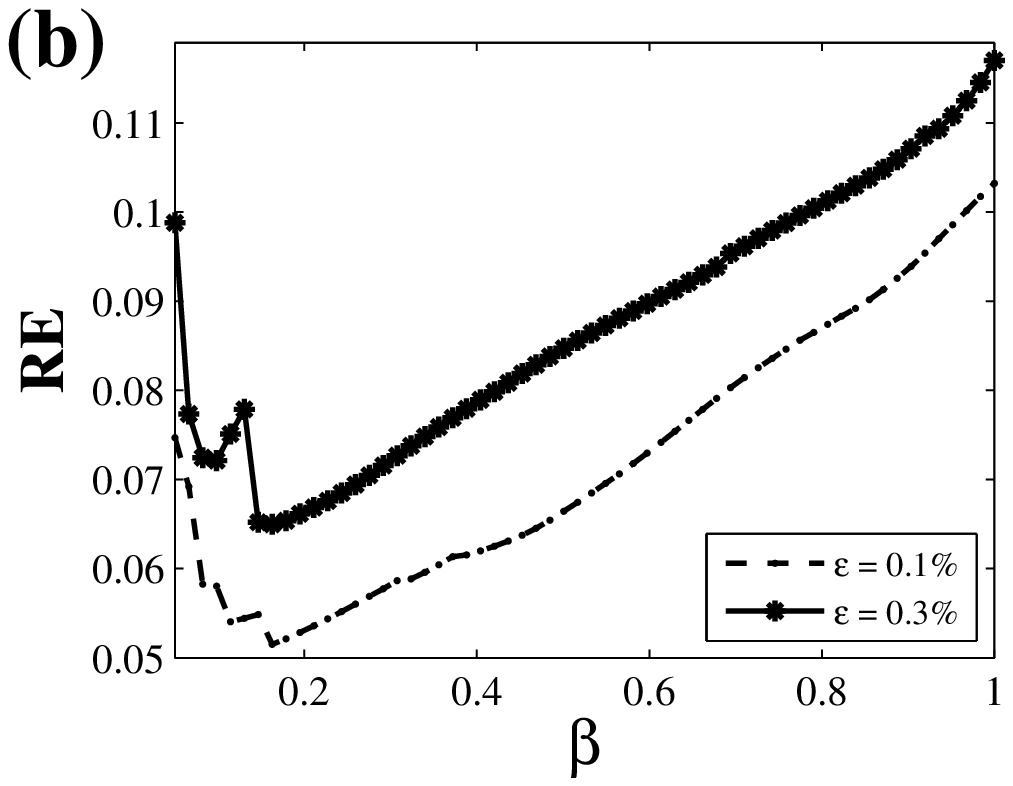}
\includegraphics[width=1.5in]{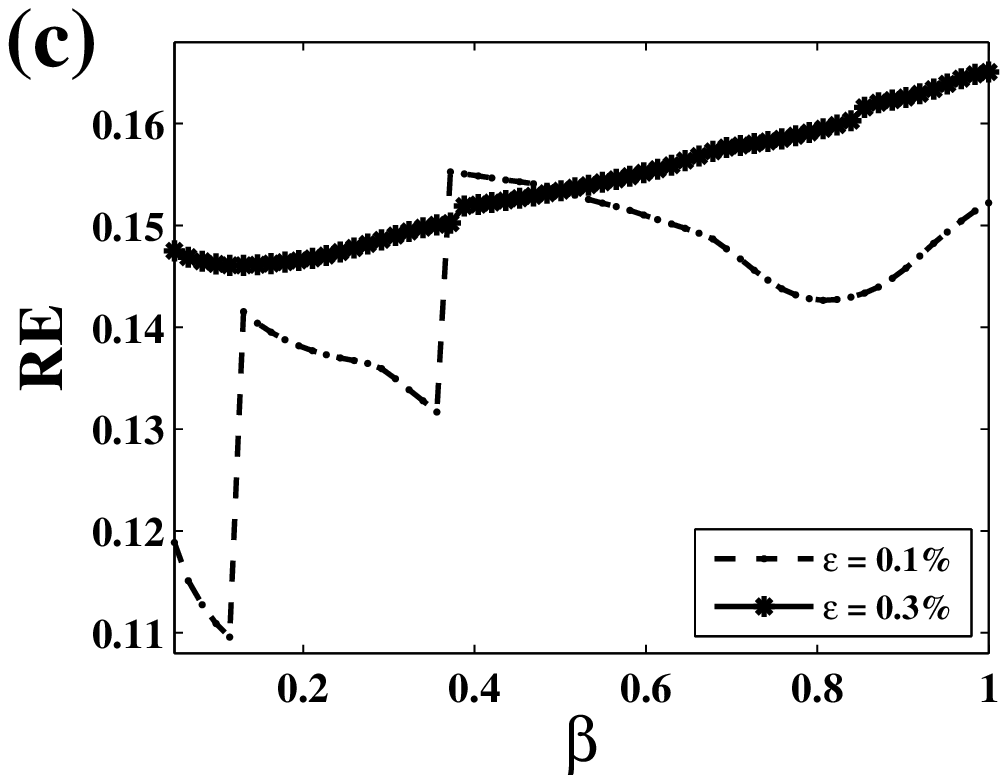}
\caption{Relative errors against parameter $\beta$ under different noise levels by Algorithm 2. (a) for Phantom A.  (b) for Phantom B.   (c) for Phantom C.}
\label{Figure 10}
\end{figure}

\begin{figure}[htbp]
\centering
\begin{minipage}[b]{0.15\textwidth}
\centering
\includegraphics[width=1.0in]{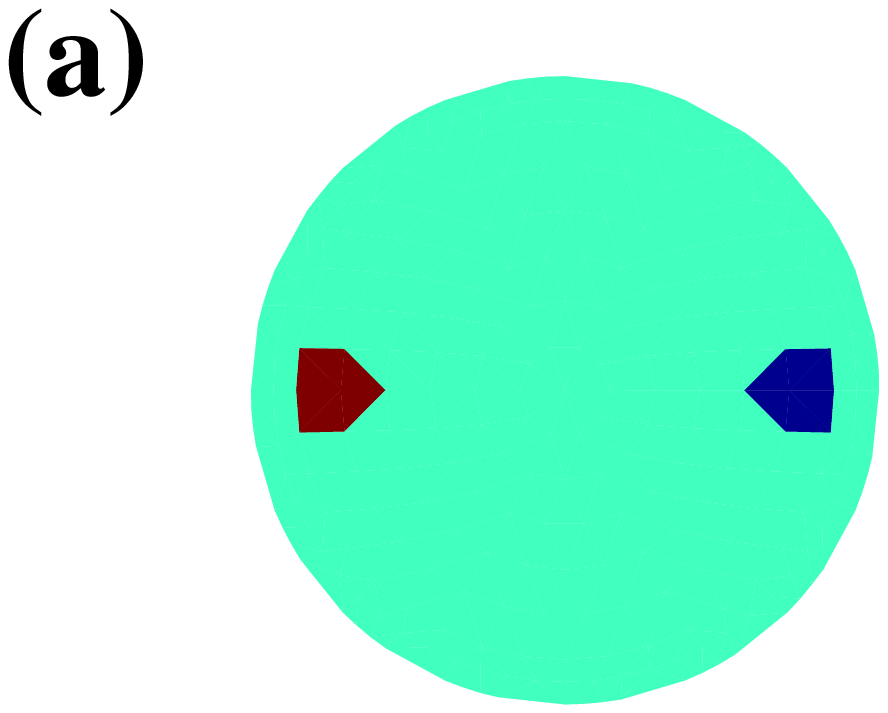}
\end{minipage}
\begin{minipage}[b]{0.15\textwidth}
\centering
\includegraphics[width=1.0in]{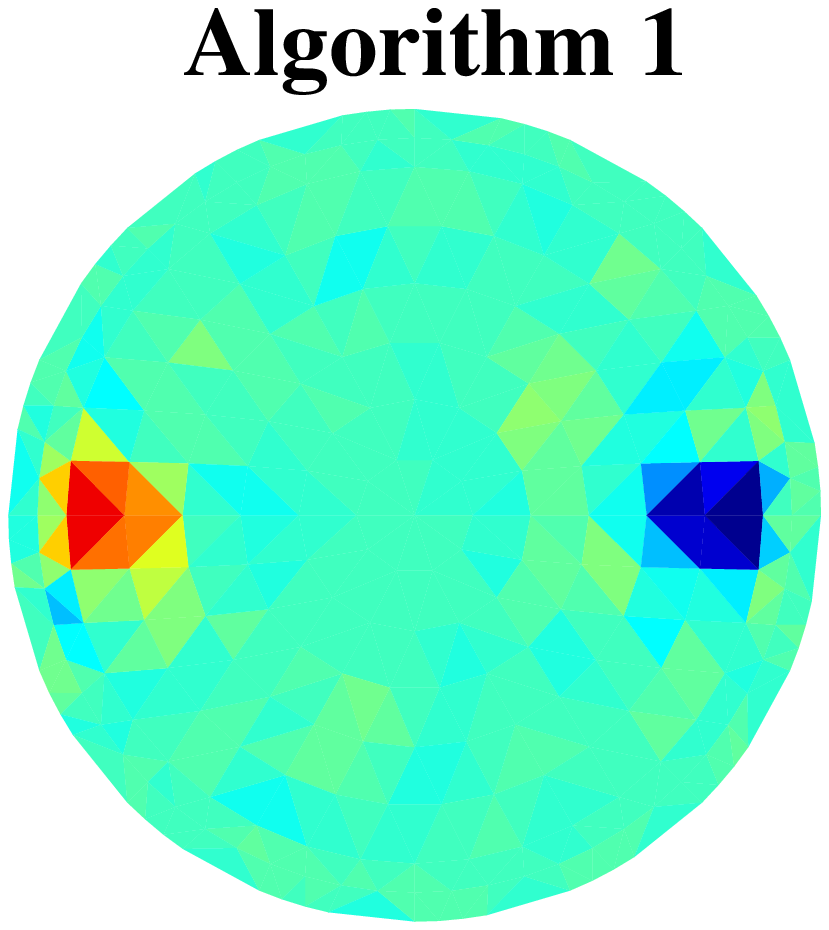}
\end{minipage}
\begin{minipage}[b]{0.15\textwidth}
\centering
\includegraphics[width=1.0in]{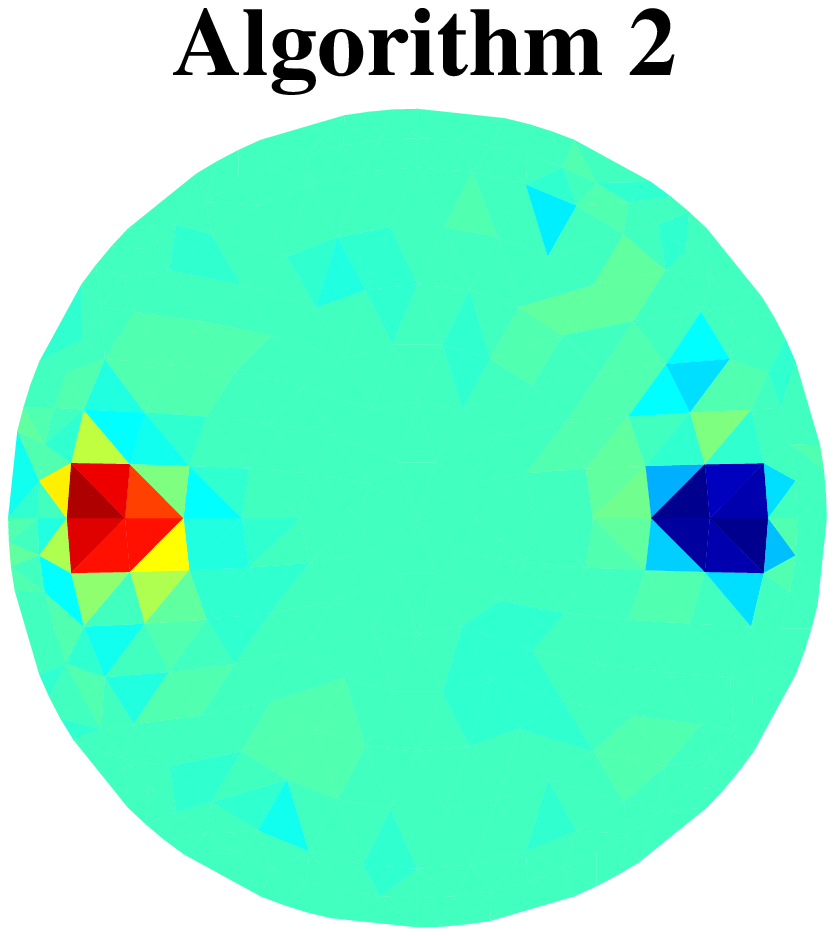}
\end{minipage}
\begin{minipage}[b]{0.15\textwidth}
\centering
\includegraphics[width=1.0in]{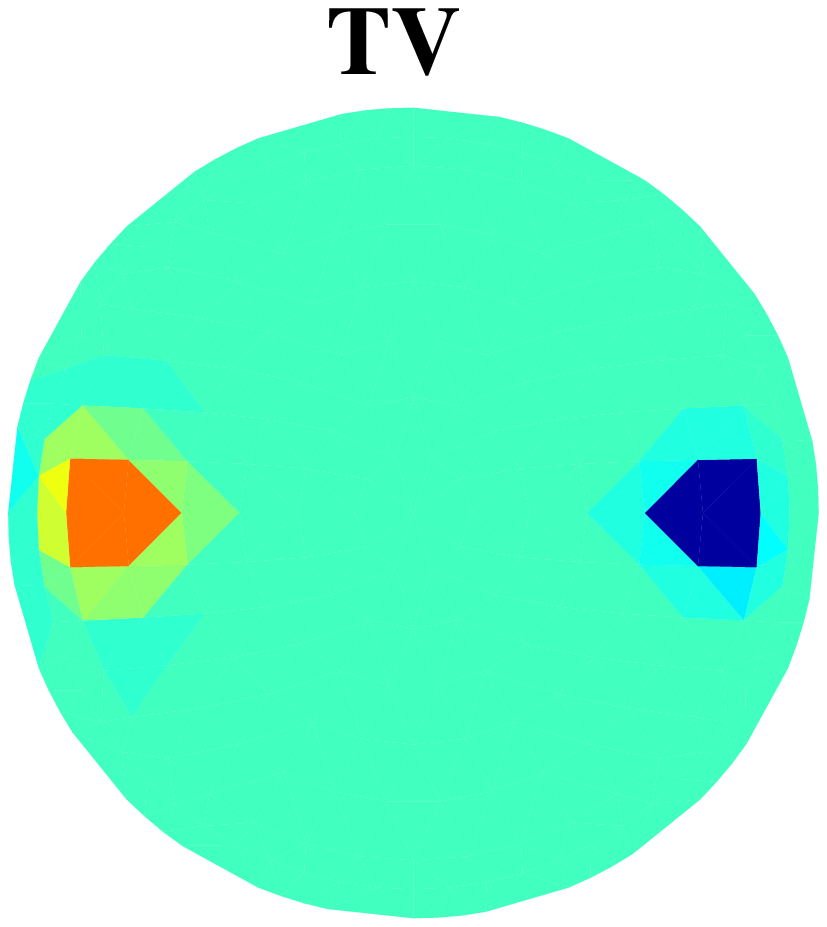}
\end{minipage}
\begin{minipage}[b]{0.15\textwidth}
\centering
\includegraphics[width=1.0in]{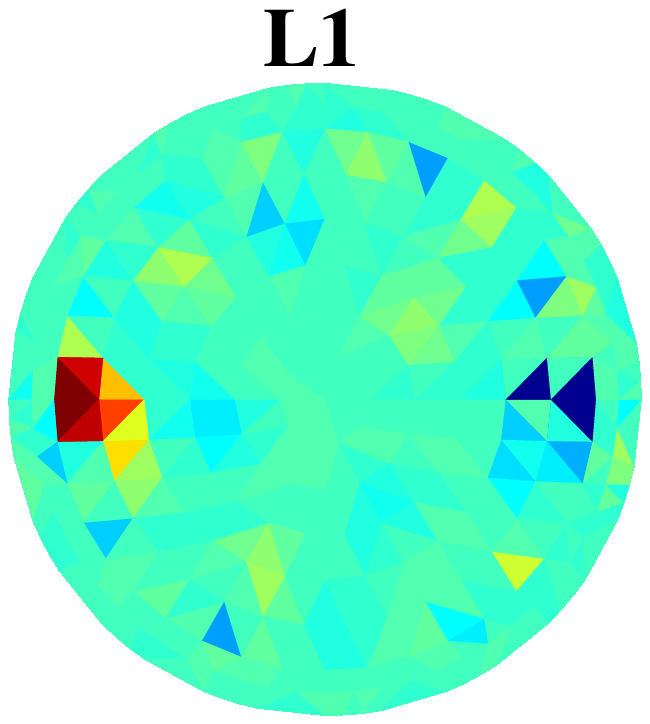}
\end{minipage}
\begin{minipage}[b]{0.15\textwidth}
\centering
\includegraphics[width=1.0in]{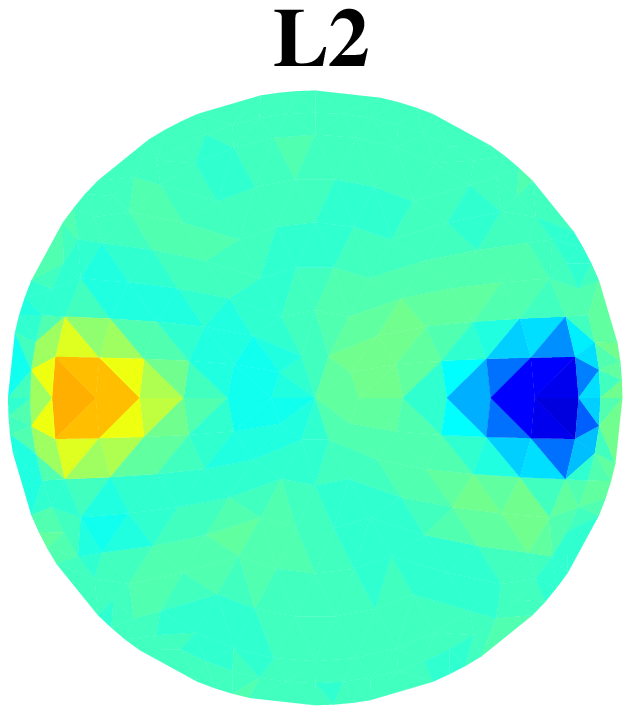}
\end{minipage}

\begin{minipage}[b]{0.15\textwidth}
\centering
\includegraphics[width=1.0in]{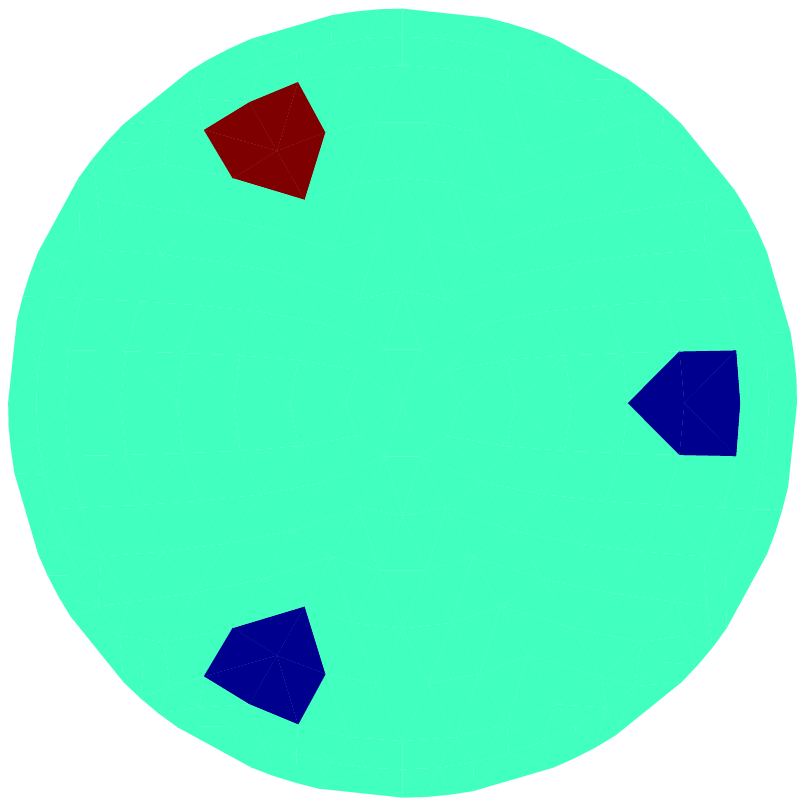}
\end{minipage}
\begin{minipage}[b]{0.15\textwidth}
\centering
\includegraphics[width=1.0in]{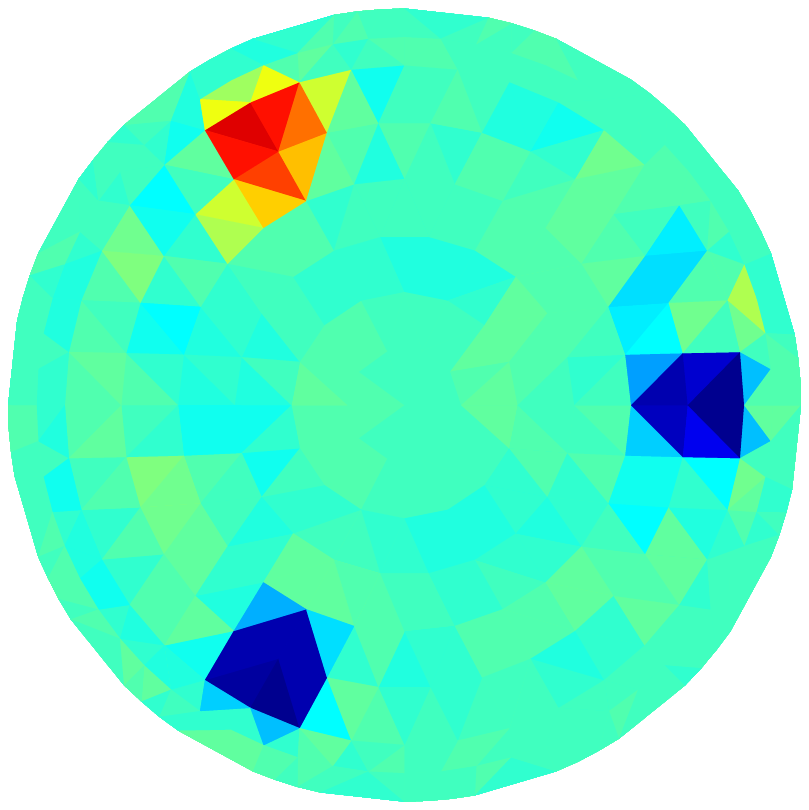}
\end{minipage}
\begin{minipage}[b]{0.15\textwidth}
\centering
\includegraphics[width=1.0in]{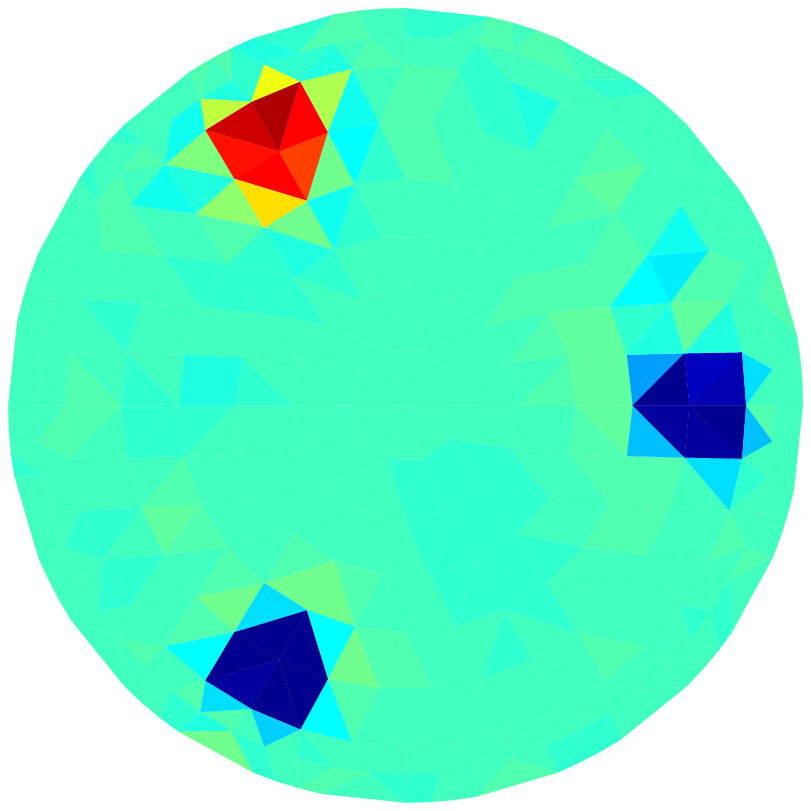}
\end{minipage}
\begin{minipage}[b]{0.15\textwidth}
\centering
\includegraphics[width=1.0in]{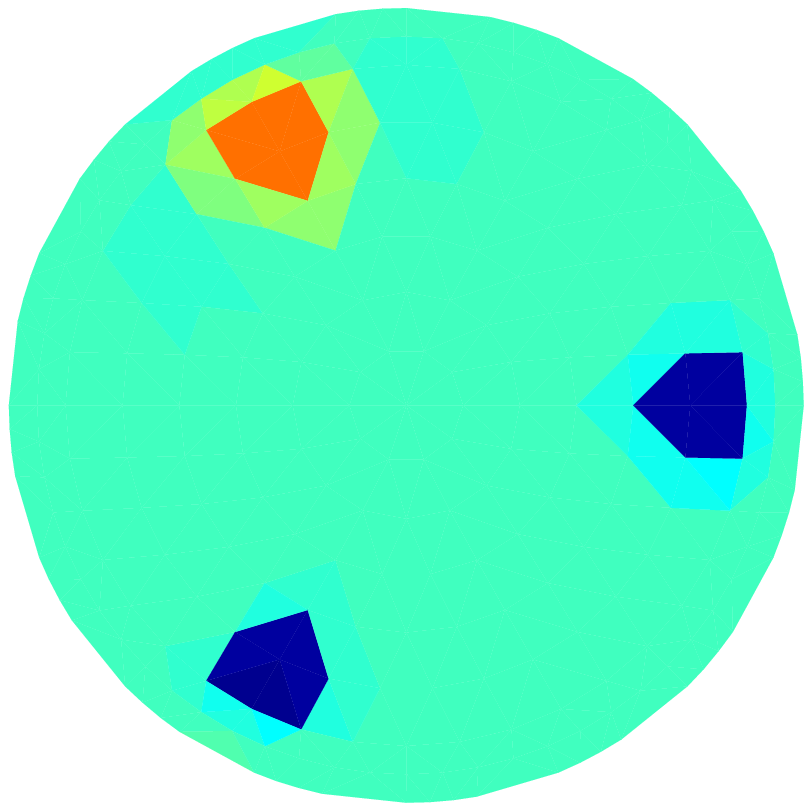}
\end{minipage}
\begin{minipage}[b]{0.15\textwidth}
\centering
\includegraphics[width=1.0in]{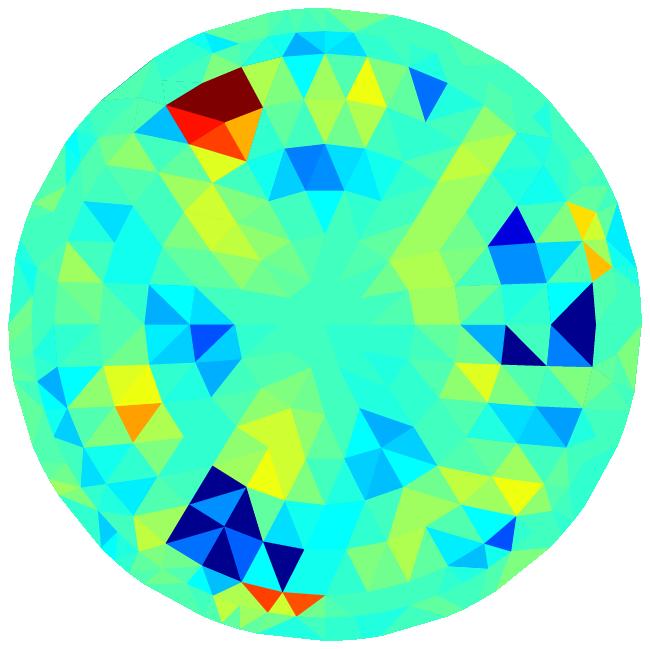}
\end{minipage}
\begin{minipage}[b]{0.15\textwidth}
\centering
\includegraphics[width=1.0in]{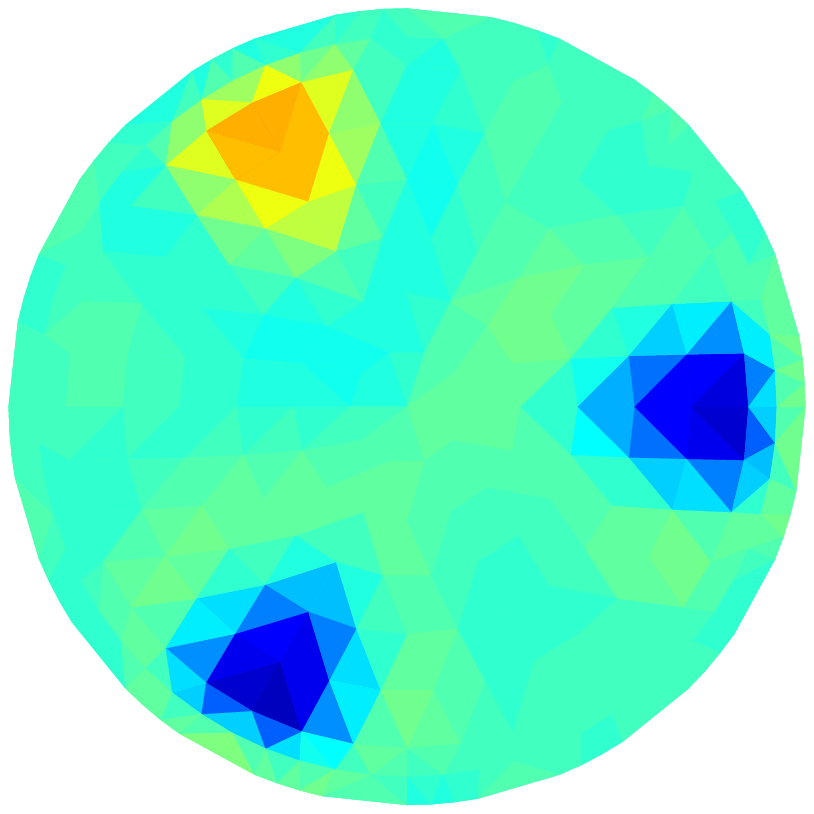}
\end{minipage}

\begin{minipage}[b]{0.15\textwidth}
\centering
\includegraphics[width=1.0in]{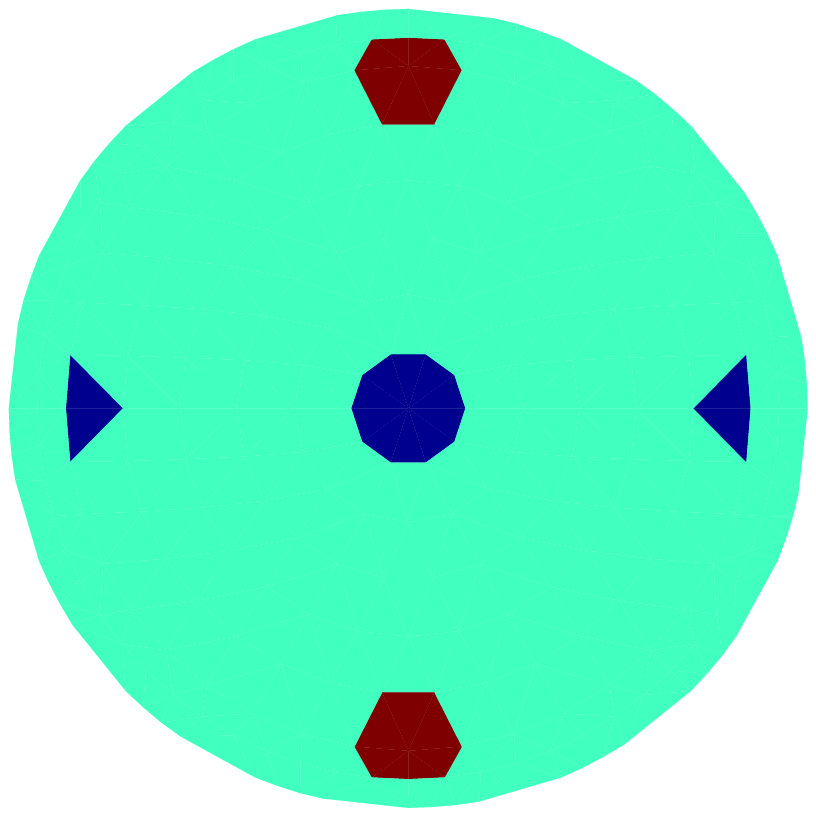}
\end{minipage}
\begin{minipage}[b]{0.15\textwidth}
\centering
\includegraphics[width=1.0in]{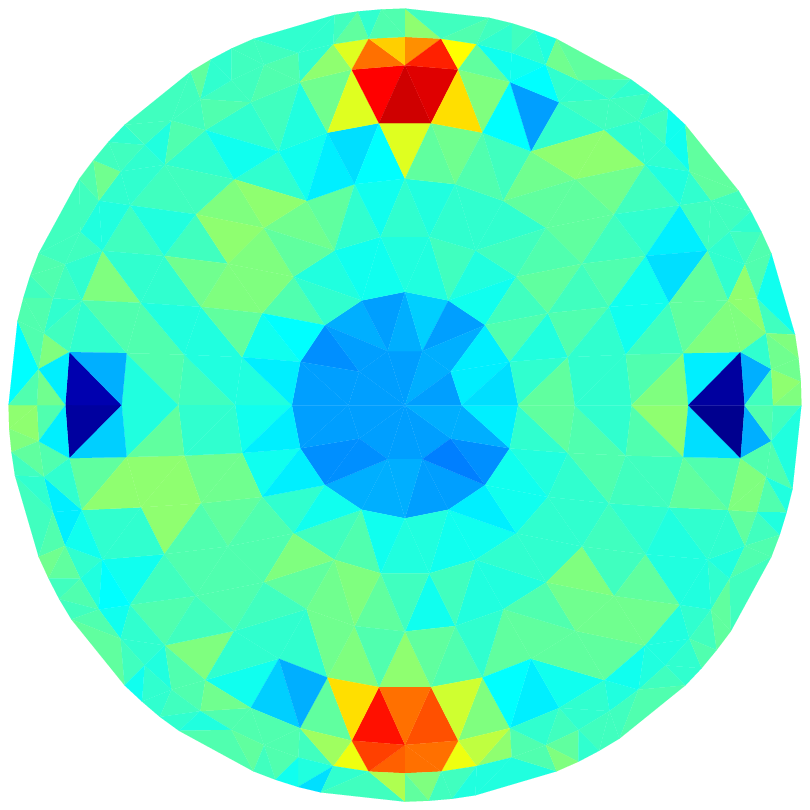}
\end{minipage}
\begin{minipage}[b]{0.15\textwidth}
\centering
\includegraphics[width=1.0in]{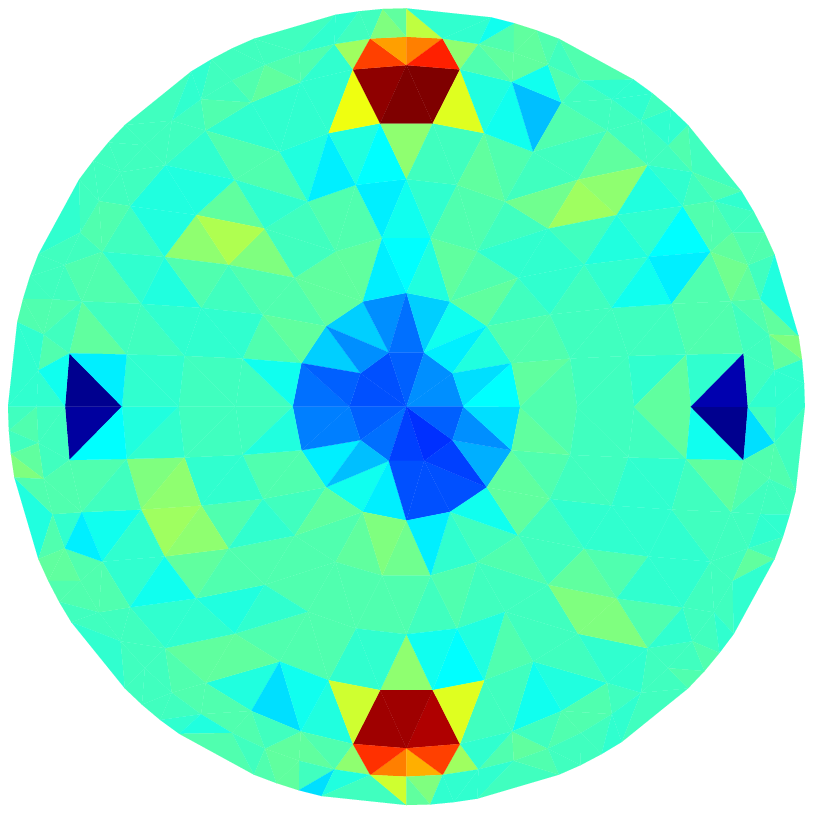}
\end{minipage}
\begin{minipage}[b]{0.15\textwidth}
\centering
\includegraphics[width=1.0in]{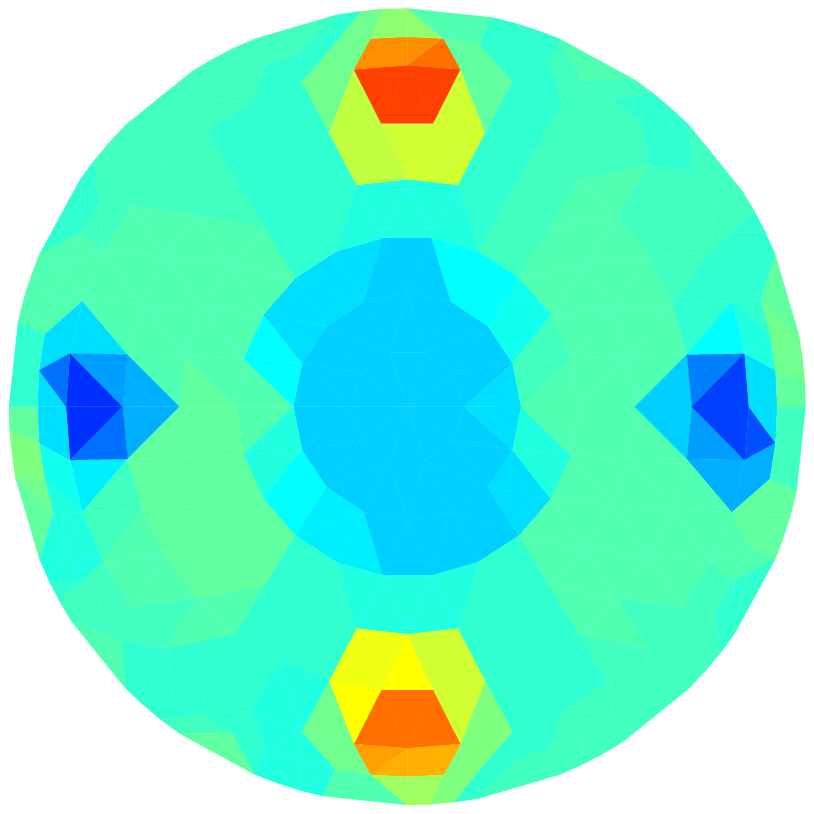}
\end{minipage}
\begin{minipage}[b]{0.15\textwidth}
\centering
\includegraphics[width=1.0in]{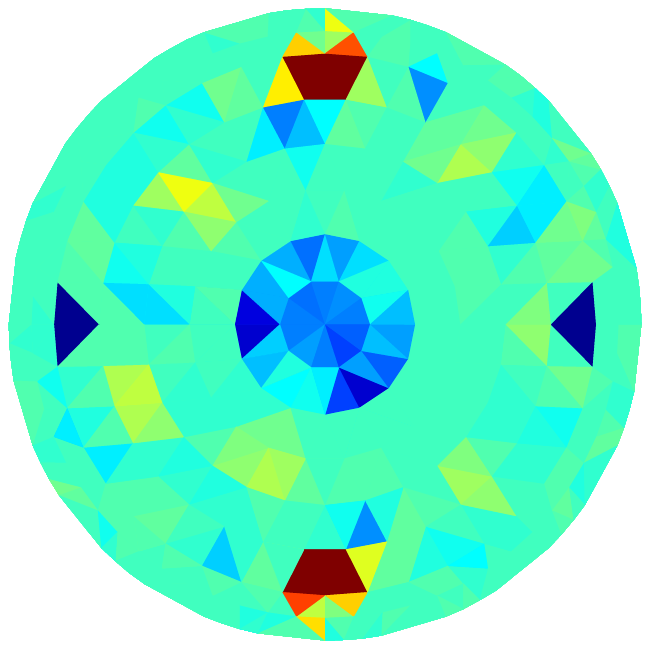}
\end{minipage}
\begin{minipage}[b]{0.15\textwidth}
\centering
\includegraphics[width=1.0in]{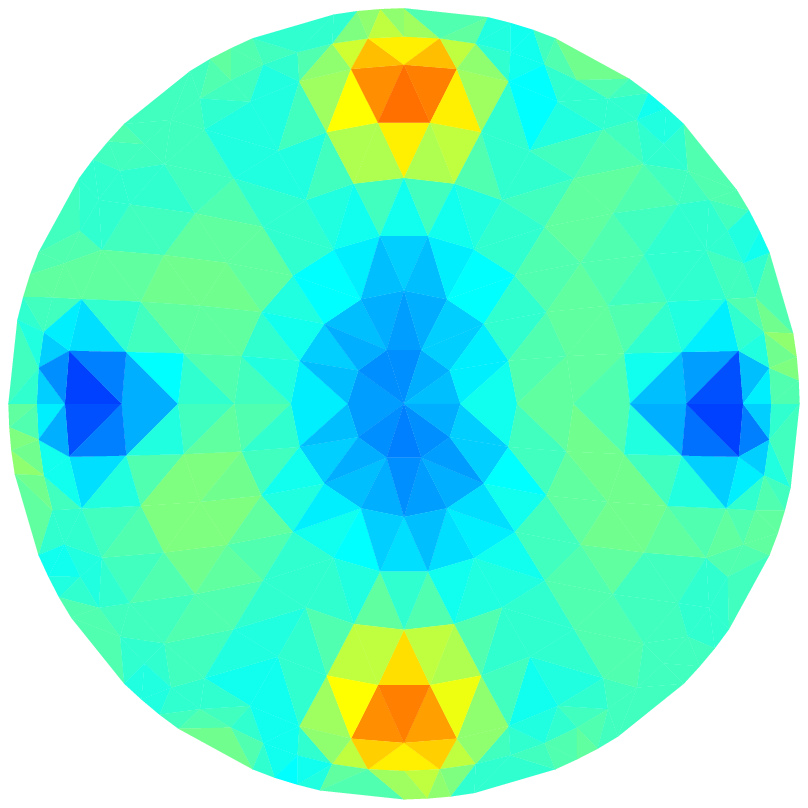}
\end{minipage}

\begin{minipage}[b]{0.15\textwidth}
\centering
\includegraphics[width=1.0in]{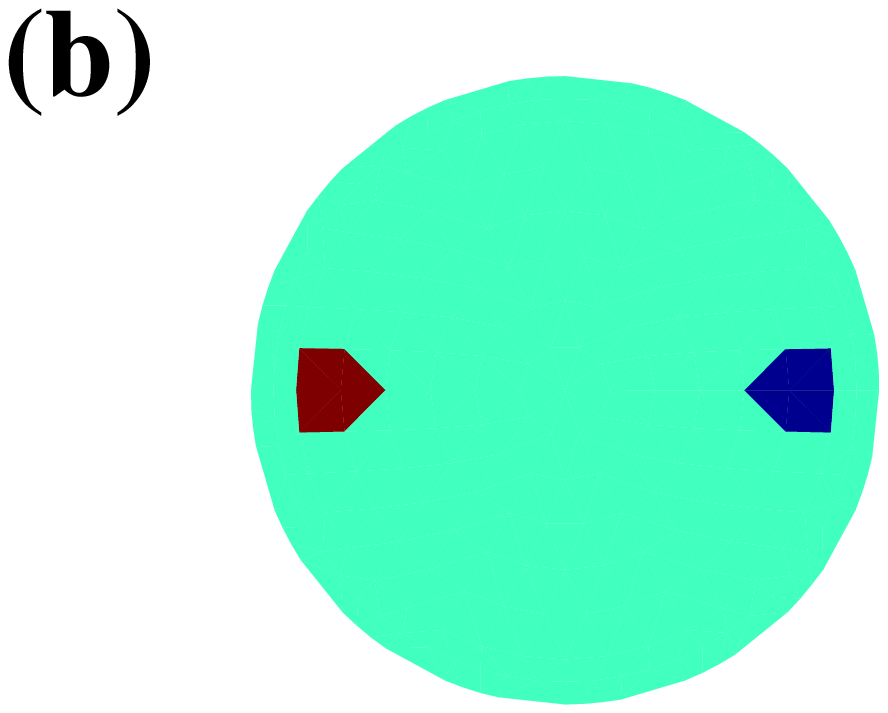}
\end{minipage}
\begin{minipage}[b]{0.15\textwidth}
\centering
\includegraphics[width=1.0in]{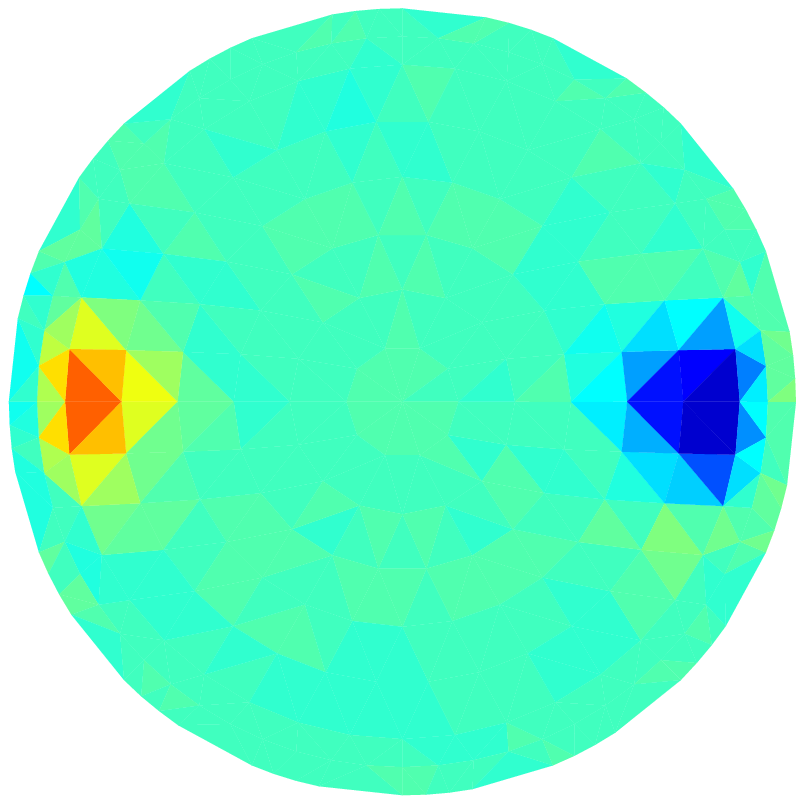}
\end{minipage}
\begin{minipage}[b]{0.15\textwidth}
\centering
\includegraphics[width=1.0in]{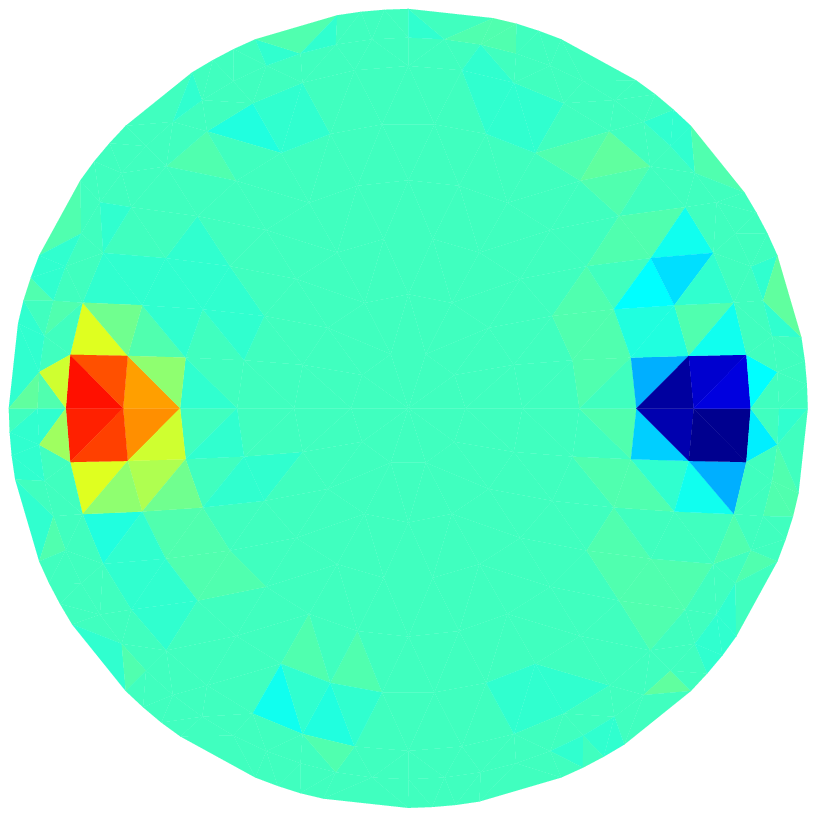}
\end{minipage}
\begin{minipage}[b]{0.15\textwidth}
\centering
\includegraphics[width=1.0in]{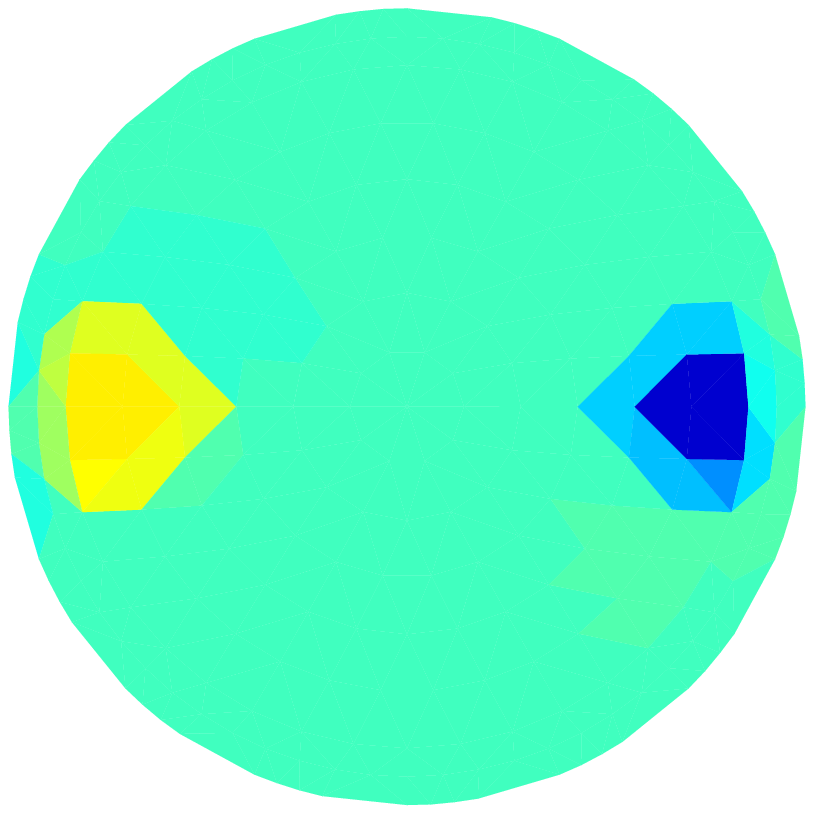}
\end{minipage}
\begin{minipage}[b]{0.15\textwidth}
\centering
\includegraphics[width=1.0in]{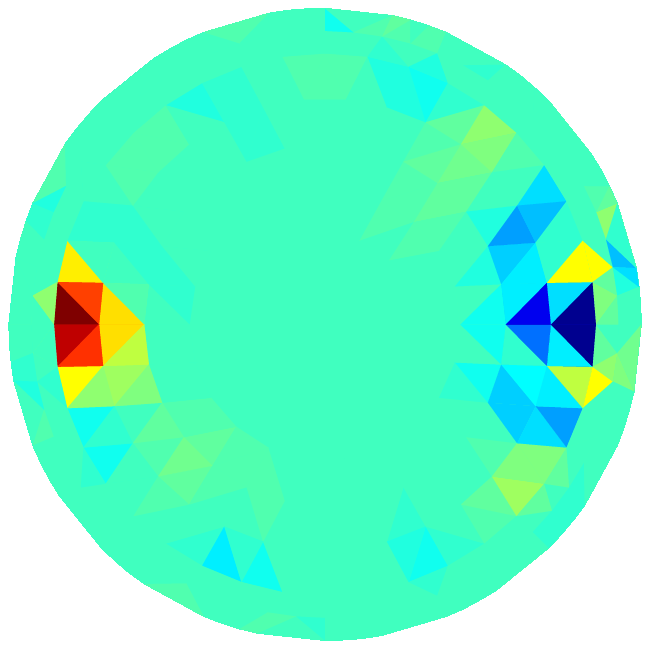}
\end{minipage}
\begin{minipage}[b]{0.15\textwidth}
\centering
\includegraphics[width=1.0in]{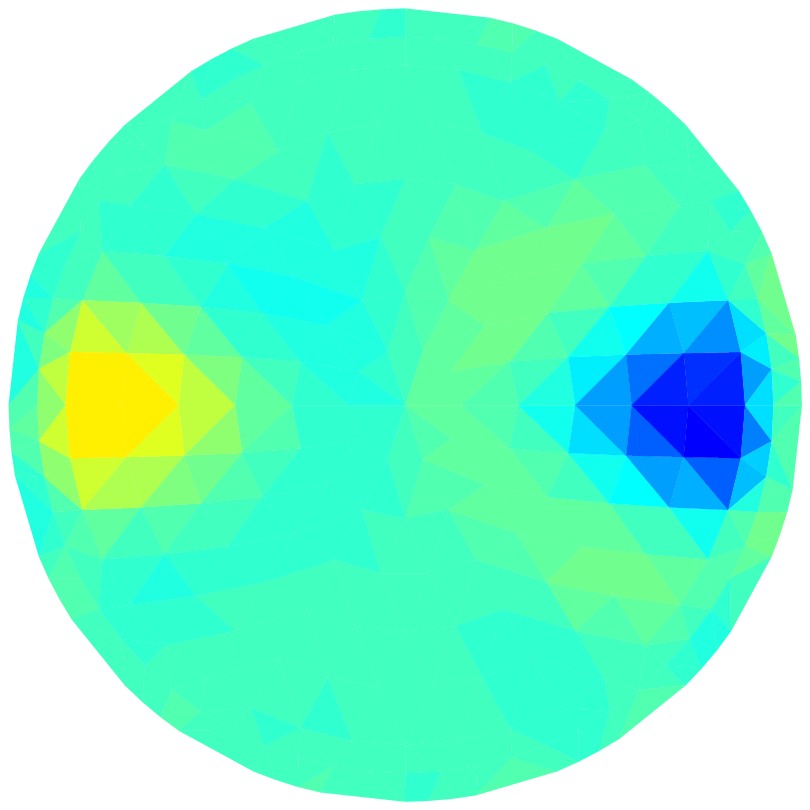}
\end{minipage}

\begin{minipage}[b]{0.15\textwidth}
\centering
\includegraphics[width=1.0in]{Model3}
\end{minipage}
\begin{minipage}[b]{0.15\textwidth}
\centering
\includegraphics[width=1.0in]{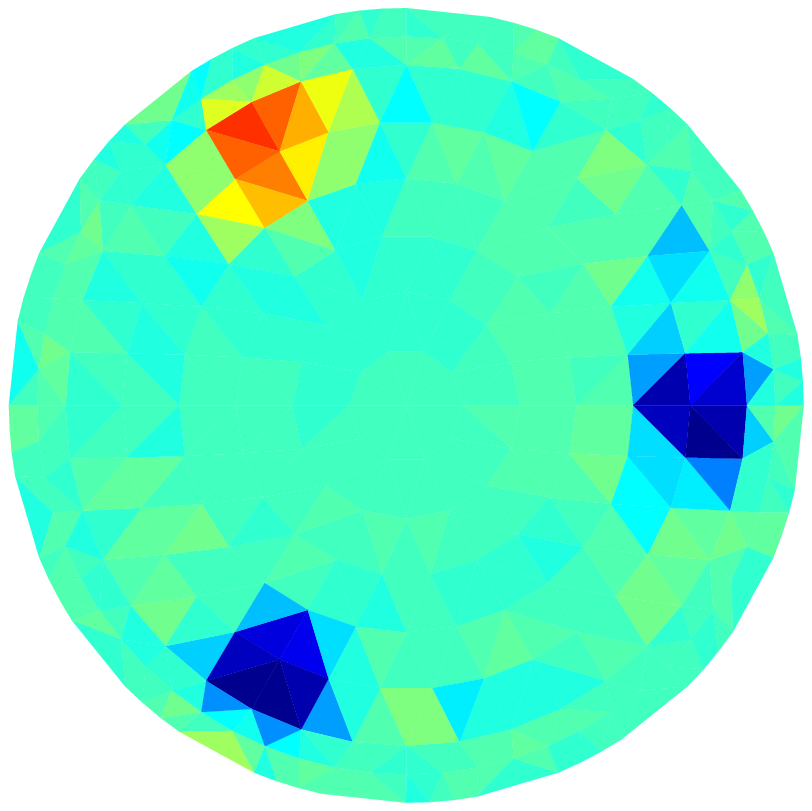}
\end{minipage}
\begin{minipage}[b]{0.15\textwidth}
\centering
\includegraphics[width=1.0in]{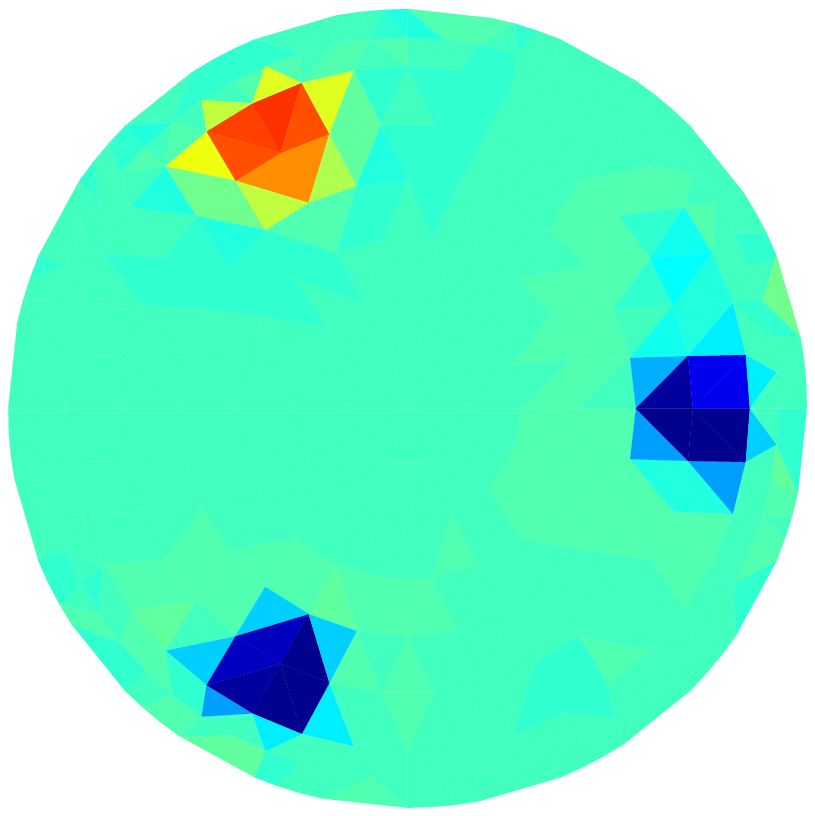}
\end{minipage}
\begin{minipage}[b]{0.15\textwidth}
\centering
\includegraphics[width=1.0in]{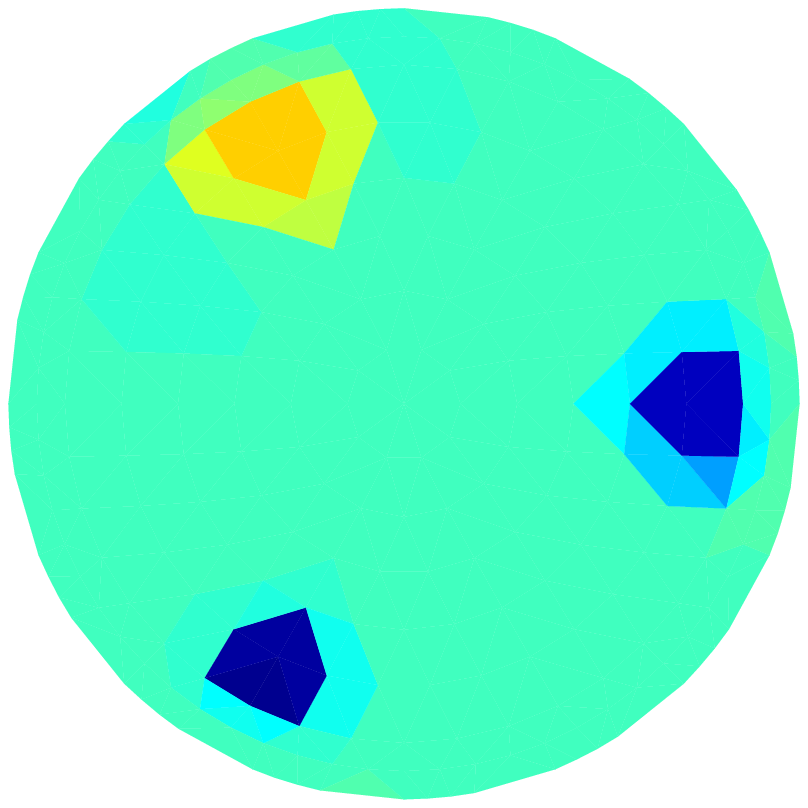}
\end{minipage}
\begin{minipage}[b]{0.15\textwidth}
\centering
\includegraphics[width=1.0in]{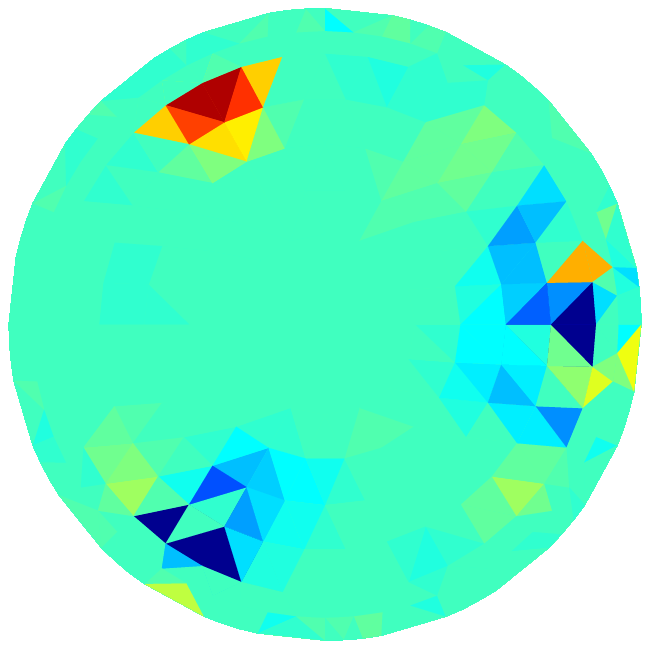}
\end{minipage}
\begin{minipage}[b]{0.15\textwidth}
\centering
\includegraphics[width=1.0in]{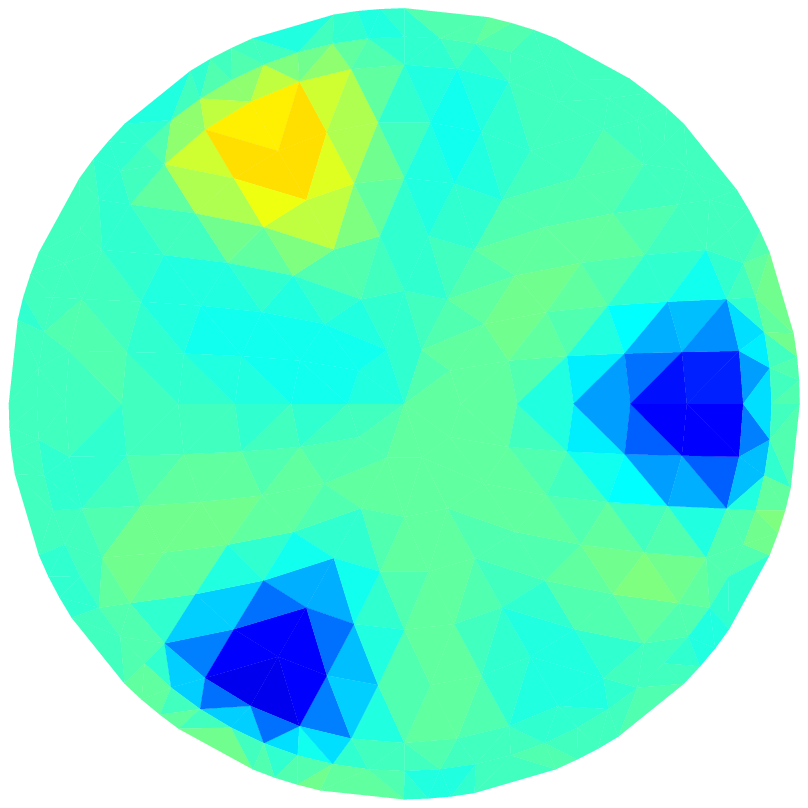}
\end{minipage}

\begin{minipage}[b]{0.15\textwidth}
\centering
\includegraphics[width=1.0in]{Model5}
\end{minipage}
\begin{minipage}[b]{0.15\textwidth}
\centering
\includegraphics[width=1.0in]{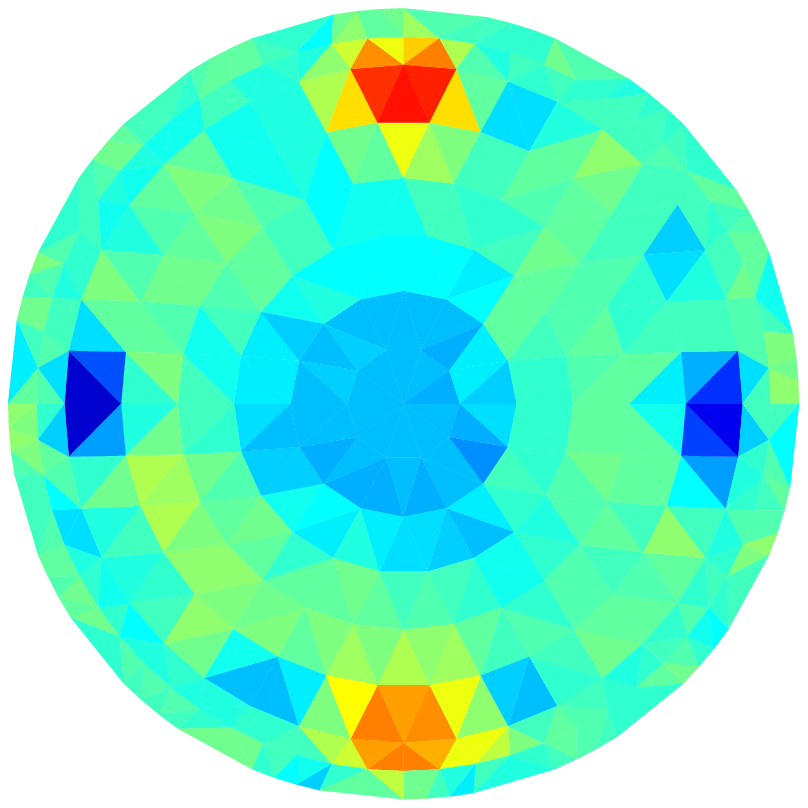}
\end{minipage}
\begin{minipage}[b]{0.15\textwidth}
\centering
\includegraphics[width=1.0in]{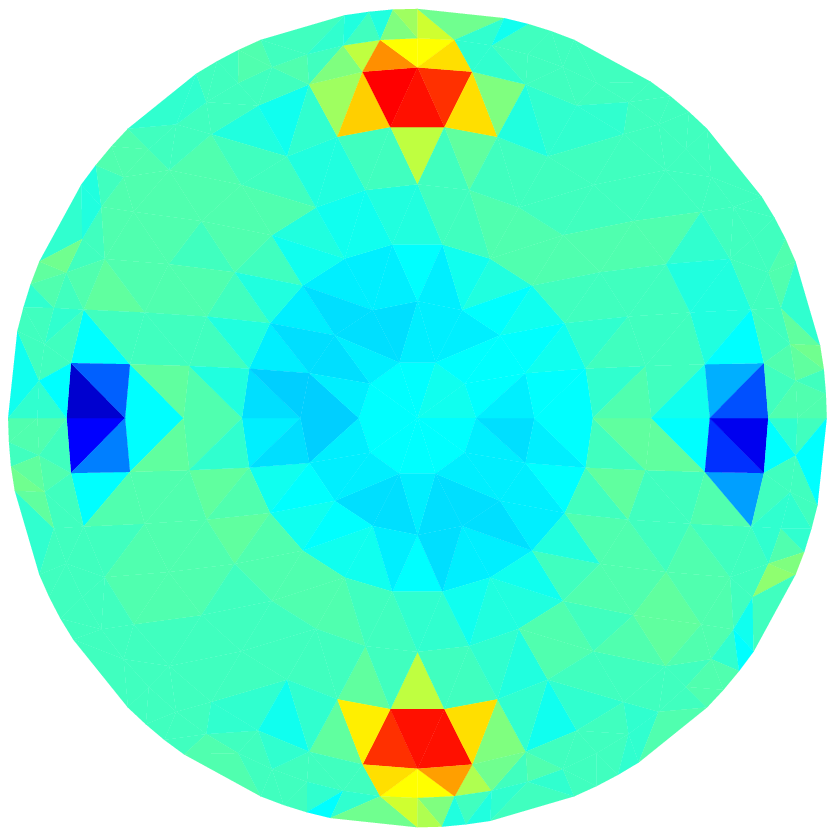}
\end{minipage}
\begin{minipage}[b]{0.15\textwidth}
\centering
\includegraphics[width=1.0in]{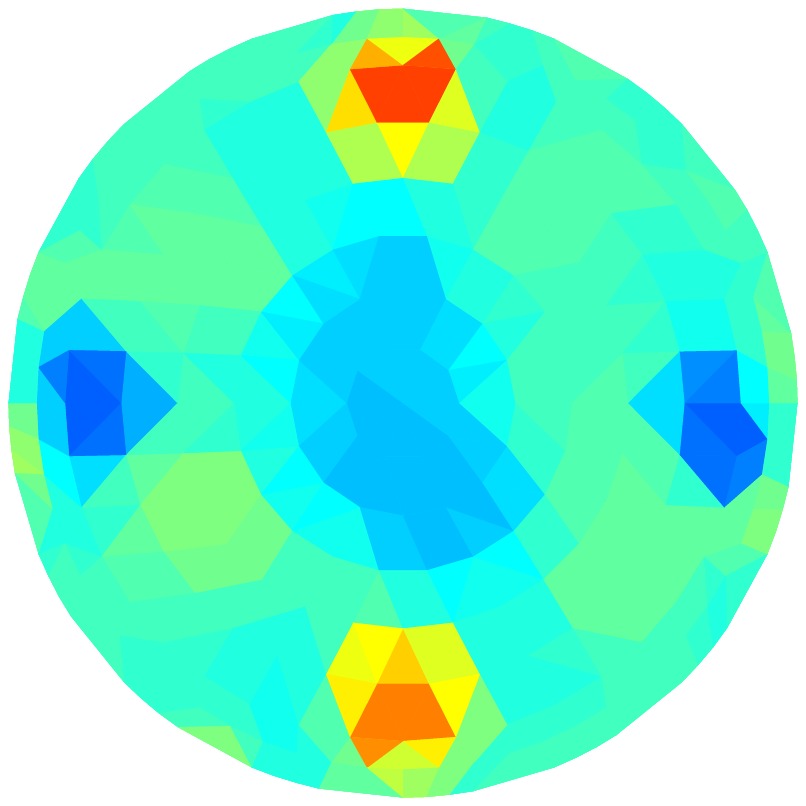}
\end{minipage}
\begin{minipage}[b]{0.15\textwidth}
\centering
\includegraphics[width=1.0in]{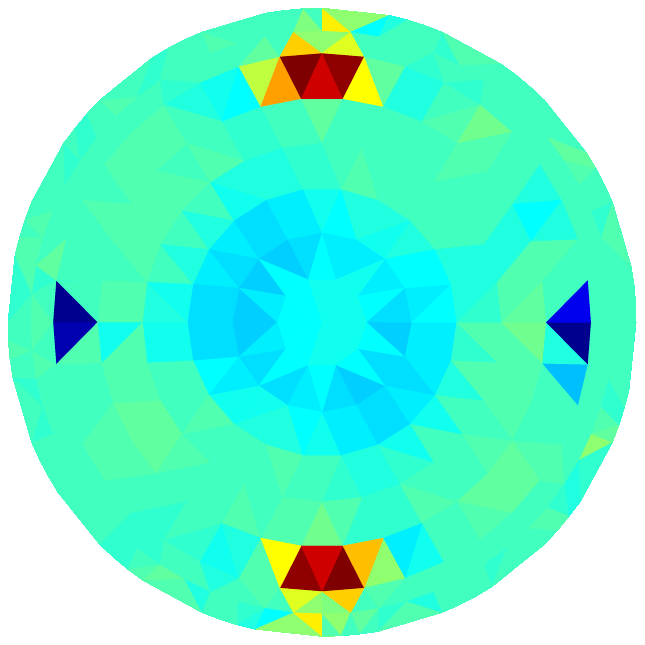}
\end{minipage}
\begin{minipage}[b]{0.15\textwidth}
\centering
\includegraphics[width=1.0in]{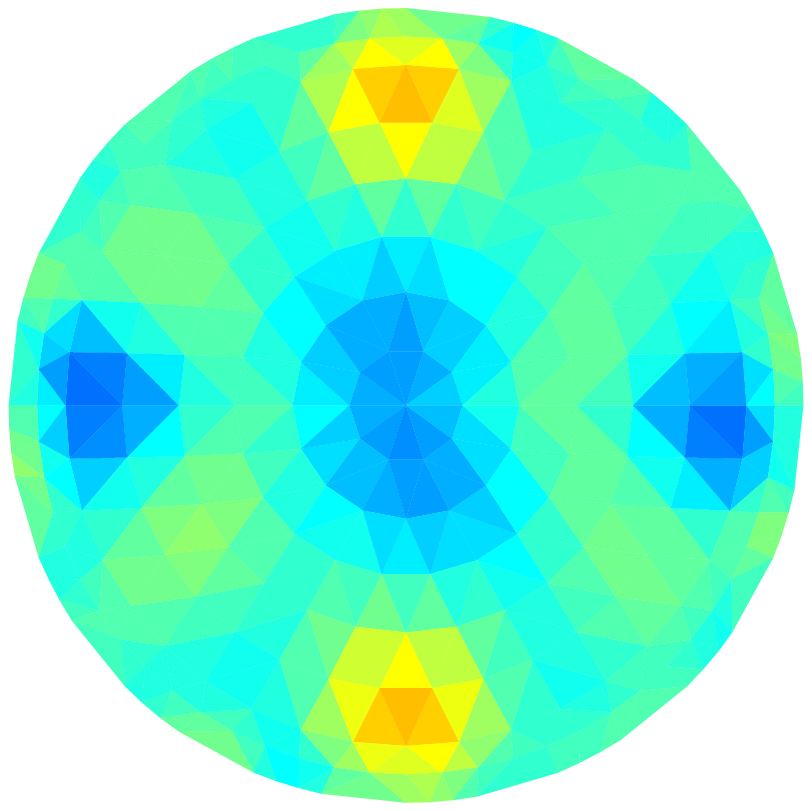}
\end{minipage}
\caption{Comparisons with TV regularization, $l_{1}$-norm regularization and $l_{2}$-norm regularization for three Phantoms.
(a) reconstructions with the noise level of $0.1\%$; (b) reconstructions with the noise level of $0.3\%$.
(All displayed by the unified colorbars from 1 to 8).}
\label{Figure 11}
\end{figure}

At the last simulation, the performance of our approach is compared with TV regularization, $l_{1}$-norm regularization and $l_{2}$-norm regularization.
Figure 11 illustrates the reconstructed images under different noise levels by Algorithm 1, Algorithm 2, TV using the lagged diffusivity method \cite{TV:2010}, $l_{1}$-norm regularization ($\beta=0$) and $l_{2}$-norm regularization ($\beta=1.0$).
It is visible that Algorithm 2 yields better reconstructions with sharper edges and less artifacts, but TV regularization appears "stair-case effect" and $l_{1}$-norm regularization appears more artifacts.
As can be expected, a little deterioration in the spatial resolution shows with the increase of noise level.
For each noise level, we compare the reconstruction errors shown in Table 5.
It is shown from Table 5 that our approach has more improvement in the accuracy, and especially Algorithm 2 is competitive and preferable in comparison with TV regularization.
At the same time, it can be seen that the reconstruction errors increase with the noise level.
Therefore simulation results indicate that our approach shows a satisfactory numerical performance and improves the imaging quality.

\section{Conclusion}
Electrical impedance tomography suffers from lower resolution.
The most important issue in the development of EIT reconstruction algorithms is to improve the imaging quality.
Using $l_{2}$-norm regularization, the reconstructed images often appear overly smooth.
This work investigates the effect of regularization on imaging quality, and proposes a novel regularization scheme constructed by a convex combination of the sparsity-promoting $l_{1}$-norm and $l_{2}$-norm for nonlinear EIT inverse problem.
The aim is to simultaneously encourage properties of sparsity and smoothness in our reconstructed images.
We then provide one simple and fast numerical iteration based on split Bregman technique for solving the joint inversion model.
Numerical simulations with synthetic data are carried out to evaluate the feasibility and effectiveness of our approach.
Results show that the proposed method have good performance on sharpening the edges of the inclusions and is robust with respect to data noise.
As a result, one promising algorithm is introduced for solving nonlinear EIT image reconstruction problem.

\begin{table}[htbp]
\centering
\caption{Comparison of image errors for the results displayed in Figure 11.}
\centering
\begin{tabular*}{13.5cm}{@{\extracolsep{\fill}}ccccccccc}
\hline
\multirow{2}{*}{}  &\multicolumn{2}{c}{Phantom A} & & \multicolumn{2}{c}{Phantom B} & &\multicolumn{2}{c}{Phantom C}\\
\cline{2-3}\cline{5-6}\cline{8-9}
noise &$0.1\%$ & $0.3\%$   &&$0.1\%$ & $0.3\%$ & & $0.1\%$ & $0.3\%$  \\
\hline
Algorithm 1 &0.0724 & 0.0886   && 0.0707 & 0.0870  &&0.1324 & 0.1500 \\
Algorithm 2 & 0.0486 & 0.0574  && 0.0515 & 0.0650  &&0.1096 & 0.1461 \\
TV & 0.0584 & 0.0846  && 0.0575 & 0.0819  &&0.1387 & 0.1579 \\
L1 & 0.0887 & 0.0842 && 0.1661 & 0.1252  &&0.1451 & 0.1503 \\
L2 & 0.0946 & 0.1040  && 0.1074 & 0.1176  &&0.1538 & 0.1652 \\
\hline
\end{tabular*}
\end{table}

\section*{Acknowledgement}
The work of Jing Wang is supported by the National Natural Science Foundation of China (NSFC) [grant number 11626092], Bo Han by NSFC [grant number 41474102], Wei Wang by NSFC [grant number 11401257].

\end{document}